\documentclass[superscriptaddress, amsmath, nofootinbib]{revtex4}
\usepackage{amsfonts}
\usepackage{graphicx}
\usepackage[latin1]{inputenc}
\usepackage{amsmath}
\usepackage{amssymb}
\usepackage[brazil]{babel}
\usepackage{indentfirst}
\usepackage{epsfig,float}

\begin{document}


\title{Basic notions of Poisson and symplectic geometry in local coordinates, with applications to Hamiltonian systems.}

\author{Alexei A. Deriglazov }
\email{alexei.deriglazov@ufjf.br} \affiliation{Depto. de Matem\'atica, ICE, Universidade Federal de Juiz de Fora,
MG, Brazil} \affiliation{Department of Physics, Tomsk State University, Lenin Prospekt 36, 634050, Tomsk, Russia}


\begin{abstract}
This work contains a brief and elementary exposition of the foundations of Poisson and symplectic geometries, with an emphasis on applications for Hamiltonian systems with second-class constraints. In particular, we clarify the geometric meaning of the Dirac bracket on a symplectic manifold and provide a proof of the Jacobi identity on a Poisson manifold.  A number of applications of the Dirac bracket are described: applications for the proof of the compatibility of a system consisting of differential and algebraic equations, as well as applications for the problem of reduction of a Hamiltonian system with known integrals of motion.
\end{abstract}

\maketitle 

\tableofcontents


\section{Introduction.}\label{SymNot}

In modern classical mechanics, equations of motion for the most of mechanical and field models can be obtained as extreme conditions for a suitably chosen variational problem. If we restrict ourselves to mechanical models, the resulting system of Euler-Lagrange equations in the general case contains differential second-order and first-order equations, as well as algebraic equations. The structure of this system becomes more transparent after the transition to the Hamiltonian formalism, which studies the equivalent system of equations, the latter no longer contains second-order equations. For the Euler-Lagrange system consisting only of second-order equations, the transition to the Hamiltonian formalism was formulated already at the dawn of the birth of classical mechanics. 
For the systems of a general form, the Hamiltonization procedure was developed by Dirac, and is known now as the Dirac formalism for constrained systems \cite{Dir_1950, GT, deriglazov2010classical}. 
In the Dirac formalism, the Hamiltonian systems naturally fall into three classes, depending on the structure of algebraic equations presented in the system.
According to the terminology adopted in \cite{GT}, they are called nonsingular, singular nondegenerate and singular degenerate theories. 

The study of these Hamiltonian systems  gave rise to a number of remarkable mathematical constructions. 
They are precisely the subject of investigation of Poisson and symplectic geometries \cite{Olver, Arn_1, Trof_1988, Per_1990, Marsden_1990,Cra_2021, Gui_1984, Vai_1994, Maslov_1993}. In particular, the geometry behind a singular nondegenerate theory could be summarized by the diagram (\ref{diag}), that clarifies geometric meaning of the famous Dirac bracket.  This will be explored in Sect. \ref{SymDiracB} to study the structure of a singular nondegenerate dynamical system. The geometric methods are widely used in current literature, in particular, for the study of massive spinning particles and bodies in external fields as well as in the analysis of propagation of light in dispersive media and in the gravitational field \cite{Kin_2021, Kas_2021, Flo_2021, Ner_2021, Kim_2022, Chen_2022, Ghez_2021, Dah_2021, Chak_2021, Heon_2021, Heon_2022, Dag_2022, Beck_2021, Giri_2022, Abd_2022, Peng_2022, Bub_2022, Ver_2022, Zhu_2022, Mik_2021, Tak_2022, Com_2021, Her_2021,  Meh_2021, Awo_2021, Saf_2020, Chu_2022, Aok_2022, Ota_2022}.

In the rest of this section we briefly describe the nonsingular and singular nondegenerate theories\footnote{Singular degenerate theories usually arise if we work within a manifestly covariant formalism, when basic variables of the theory transform linearly under the action of the Poincare group. Their description can be found in \cite{Dir_1950, GT, deriglazov2010classical}.}. Then we fix our notation and recall some basic notions of the theory of differentiable manifolds, that will be useful in what follows.

\subsection{Nonsingular theories.}  

This is the name of mechanical systems that in the Hamiltonian formulation can be described using only first-order differential equations (called Hamiltonian equations)  
\begin{eqnarray}\label{sim.0.01}
\dot q^a=\frac{\partial H}{\partial p_a}, \qquad 
\dot p_a=-\frac{\partial H}{\partial q^a},
\end{eqnarray}
where $H(q, p)$ is a given function and $\dot q^a=\frac{d}{d\tau}q^a$. 
The variables $q^a(\tau)$ describe position of the system, while $p_a(\tau)$ are related to the velocities, and in simple cases are just proportional to them. 
The equations show that the function $H(q, p)$, called the Hamiltonian, encodes in fact all the information about dynamics of the mechanical system. The equations can be written in a more compact form, if we introduce an operation assigning to every pair of functions $A(q, p)$ and $B(q, p)$ a new function,  denoted $\{A, B\}_P$, as follows:
\begin{eqnarray}\label{sim.0.02}
\{ A, B\}_P=\frac{\partial A}{\partial q^a}\frac{\partial B}{\partial p_a}-\frac{\partial B}{\partial q^a}\frac{\partial A}{\partial p_a}.
\end{eqnarray}
This is called the canonical Poisson bracket of $A$ and $B$. Then the Hamiltonian equations acquire the form 
\begin{eqnarray}\label{sim.0.03}
\dot z^i=\{z^i, H\}_P,   
\end{eqnarray}
where $z^i=(q^a, p_b)$, $i=1, 2, \ldots , 2n$.
The equations determine integral lines $z^i(\tau)$ of the vector field $\{z^i, H\}_P$ on $\mathbb{R}^{2n}$, created by function $H$. For smooth vector fields, the Cauchy problem, that is Eqs. (\ref{sim.0.01}) with the initial condidions $z^i(\tau_0)=z^i_0$, has unique solution in a vicinity of any point $z^i_0 \in\mathbb{R}^{2n}$. Formal solution to these equations in terms of power series is as follows \cite{deriglazov2010classical}
\begin{eqnarray}\label{sim.04}
z^i(\tau, z^j_0)=e^{\tau\{z^k_0, H(z^i_0)\}_P\frac{\partial}{\partial z^k_0}}z^i_0.   
\end{eqnarray}
The functions $z^i(\tau, z^j_0)$ depend on $2n$ arbitrary constants $z^j_0$, and hence represent a general solution to the system (\ref{sim.0.03}). 

In the Lagrangian formalism, an analogue of this formula is not known. So, Eq. (\ref{sim.04}) can be considered as the first example, showing the usefulness of the transition from the Lagrangian to the Hamiltonian description.

\subsection{Singular nondegenerate theories.}  

Consider the system consisting of differential and algebraic equations
\begin{eqnarray}
\dot z^i=\{z^i, H\}_P,  \label{sim.0.5} \qquad  \qquad  \qquad  \qquad \qquad ~  ~   \\
\Phi^\alpha(z^i)=0, \qquad ~ ~ \alpha =1, 2, \ldots ,  2p< 2n,   \label{sim.0.6}
\end{eqnarray}
where $H(z^i) $ and $\Phi^\alpha(z^i)$ are given functions. It is supposed that $\Phi^\alpha(z^i)$ are functionally independent functions\footnote{We recall that the functional independence of functions $\Phi^\alpha$ guarantee that the system (\ref{sim.0.6}) can be resolved with respect to $2p$ variables $z^\alpha$ among $z^i$, then $z^\alpha=f^\alpha(z^b)$ are parametric equations of the surface $\Phi^\alpha=0$.} (constraints), so the equations (\ref{sim.0.6}) determine $2n-2p$\,-dimensional surface $\mathbb{N}$. The system is called singular nondegenerate theory, if there are satisfied the following two conditions. The first condition is 
\begin{eqnarray}\label{sim.0.5.1}
\left.\det \{\Phi^\alpha, \Phi^\beta\}_P\right|_{\Phi^\alpha=0}\ne 0,
\end{eqnarray}
hence the name "nondegenerate system". In the Dirac formalism, functions with the property (\ref{sim.0.5.1}) are called second-class constraints. The second condition is that the functions $\{\Phi^\alpha, H\}_P(z^i)$ vanish on the surface $\mathbb{N}$
\begin{eqnarray}\label{sim.0.5.2}
\left.\{\Phi^\alpha, H\}_P\right|_{\Phi^\alpha=0}=0.
\end{eqnarray}
The two conditions guarantee the existence of solutions to the system (\ref{sim.0.5}), (\ref{sim.0.6}). To discuss this point, we adopt the following \par

\noindent {\bf Definition 1.1.} The system (\ref{sim.0.5}), (\ref{sim.0.6}) is called self-consistent, if through any point of the surface $\mathbb{N}$ passes a solution of the system.

For the self-consistent system, its formal solution can be written as in (\ref{sim.04}), it is sufficient to take the integration constants $z^i_0$ on the surface of constraints. 

Let us discuss the self-consistency of the system. Given point of the surface $\mathbb{N}$, there is unique solution of equations (\ref{sim.0.5}), that passes through this point. It will be a solution of the whole system, if it entirely lies on the surface: 
\begin{eqnarray}\label{sim.0.6.1}
\dot z^i=\{z^i, H\}_P \quad \mbox{and} \quad \Phi^\alpha(z^i(0))=0, \quad \mbox{implies} \quad \Phi^\alpha(z^i(\tau))=0 \quad \mbox{for all} \quad \tau.
\end{eqnarray}
This is a strong requirement, and equations (\ref{sim.0.5.1}) and (\ref{sim.0.5.2}) turn out to be the sufficient conditions for its fulfilment. In a physical context, the proof with use of special coordinates of $\mathbb{R}^{2n}$ was done in \cite{GT}. A more simple proof with use of Dirac bracket will be presented in Sect. \ref{SymDiracB}.  

An example of a self-consistent system like (\ref{sim.0.5}), (\ref{sim.0.6}) will be considered in Sect. \ref{SymDirac}, see Affirmation 7.5.

Here we discuss the necessity of the condition (\ref{sim.0.5.2}). \par 

\noindent {\bf Affirmation 1.1.} Consider the system (\ref{sim.0.5}), (\ref{sim.0.6}) with functionally independent functions $\Phi^\alpha$. 
Then  $\left.\{\Phi^\alpha, H\}_P\right|_{z^i(\tau)}=0$ for any solution $z^i(\tau)$, if any.  That is the algebraic equations $\{\Phi^\alpha, H\}_P=0$ are consequences of the system. \par

\noindent {\bf Proof.}  Let the system admits the solution $z^i(\tau)$. Then $\Phi^\alpha(z^i(\tau))=0$ for all $\tau$, this implies $\dot \Phi^\alpha(z^i(\tau))=0$. On other hand we get 
\begin{eqnarray}\label{sim.07}
0=\dot \Phi^\alpha=\left.\frac{\partial \Phi^\alpha}{\partial z^i}\right|_{z^i(\tau)}\dot z^i(\tau)=
\left.\frac{\partial \Phi^\alpha}{\partial z^i}\right|_{z(\tau)}\left.\{z^i, H\}_P\right|_{z(\tau)}=\left.\{\Phi^\alpha, H\}_P\right|_{z(\tau)}.
\end{eqnarray}
that is $\{\Phi^\alpha, H\}_P=0$ for any solution $z^i(\tau)$. $\blacksquare$ \par

\noindent {\bf Affirmation 1.2.} If the system (\ref{sim.0.5}), (\ref{sim.0.6}) with functionally independent functions $\Phi^\alpha$ is self-consistent, the conditions (\ref{sim.0.5.2}) hold. \par

\noindent {\bf Proof.}  Let $z^i_0$ be any point of the surface $\Phi^\alpha=0$. Due to the self-consistency, there is a solution $z^i(\tau)$ that passes through this point,  $z^i(0)=z^i_0$. As the equation$\{\Phi^\alpha, H\}_P=0$ is a consequence of the system, we 
have $\{\Phi^\alpha, H\}_P(z^i(\tau))=0$, in particular $\{\Phi^\alpha, H\}_P(z^i(0))=\{\Phi^\alpha, H\}_P(z^i_0)=0$, that is it vanishes at all points of the surface $\mathbb{N}$. $\blacksquare$

Consider the system (\ref{sim.0.5}), (\ref{sim.0.6}), and now suppose that some of the functions $\{\Phi^\alpha, H\}_P$ do not vanish identically on the surface $\mathbb{N}$. As we saw above, this means that the system is not a self-consistent. Then we can look for a sub-surface of $\mathbb{N}$ where the system could be a self-consistent. The procedure is as follows. We separate the functionally independent functions among $\{\Phi^\alpha, H\}_P$,  say $\Psi^1, \Psi^2, \ldots  , \Psi^k$. As the equations $\Psi^a=0$ are consequences of the system  (\ref{sim.0.5}), (\ref{sim.0.6}), we add them to the system, obtaining an equivalent system of equations. If the set $\Phi^\alpha, \Psi^a$ is composed of functionally independent functions, we repeat the procedure, analysing the functions $\{\Psi^a, H\}$, and so on. 
Since the number of functionally independent functions cannot be more than $2n$, the procedure will end at some step. If, in addition to this, the resulting set of functions satisfies the condition (\ref{sim.0.5.1}), we arrive at the self-consistent system of equations: $\dot z^i=\{z^i, H\}_P$,  $\Phi^\alpha(z^i)=0$, $\Psi^a=0, \ldots$ . 

It remains to discuss what happens if, at some stage, the extended system of algebraic equations consists of functionally dependent functions. Without loss of generality, we assume that the extended system is $\dot z^i=\{z^i, H\}_P$, $\Phi^\alpha=0$,  $\Psi\equiv\{\Phi^1, H\}=0$. By construction, it is equivalent to the original system, the function $\Psi(z^i)$ does not vanish identically on  $\mathbb{N}$, and the functions $\Phi^\alpha,  \Psi$ are functionally dependent. As $\Phi^\alpha$ are functionally independent, we present the equations $\Phi^\alpha(z^i)=0$ in the form $z^\alpha=f^\alpha(z^b)$, and substitute them into the expression for $\Psi(z^i)$, obtaining the system $\dot z^i=\{z^i, H\}_P$, $z^\alpha-f^\alpha(z^b)=0$,  $\Psi(z^b, f^\alpha(z^b))=0$, which is equivalent to (\ref{sim.0.5}), (\ref{sim.0.6}). The function $\Psi(z^b, f^\alpha(z^b))\ne 0$ identically. On other hand, it does not depend on $z^b$ (otherwise we could write it in the form like $z^1=\psi(z^2, z^3, \ldots)$, then the functions $z^\alpha-f^\alpha(z^b)$, $z^1-\psi(z^2, z^3, \ldots)$ are functionally independent). So, the only possibility is $\Psi=c=const\ne 0$. This means that the system (\ref{sim.0.5}), (\ref{sim.0.6}) contains the equation $c=0$, where $c\ne0$. Hence the system is contradictory and has no solutions at all. 

It should be noted that the outlined procedure for obtaining a self-consistent system lies at the corner of the Dirac method \cite{Dir_1950}.

Since all trajectories of the system (\ref{sim.0.5}), (\ref{sim.0.6}) lie on the surface $\Phi^\alpha(z^i)=0$ with coordinates, say $z^b$, a number of questions naturally arise. Can equations for independent variables $z^b$ be written in the form of a Hamiltonian system like (\ref{sim.0.5})? What are the Hamiltonian $H(z^b)$ and the bracket $\{ {}, {} \}_{{\mathbb N}}$ in these equations, and how they should be constructed?  Is the new bracket a kind of restriction of the original one to ${\mathbb N}$? The answers to these questions will be given in Sect. \ref{SymDirac}. In particular, we will show that the new bracket is a restriction of the Dirac bracket to ${\mathbb N}$, and not a restriction of the original bracket.

\subsection{Smooth manifolds.}\label{SymNot_C}

{\bf Notation.} Latin indices from the middle of alphabet are used to represent coordinates $z^k$ of a manifold $\mathbb{M}_n$ and run from $1$ to $n$. If coordinates are divided on two groups, we write $y^k=(y^\alpha, y^b)$, that is Greek indices from the beginning of alphabet are used to represent one group, while Latin indices from the beginning of alphabet represent another group. The notation like $U_i(z^j)$ means that we work with the functions $U_i(z^1, z^2, \ldots , z^n)$, where $i=1, 2, \ldots , n$. The notation like $\partial_i A(z^k)|_{z^k\rightarrow f^k(y^j)}$ means that in the expression $(\partial A(z^k)/\partial z^i)$ the symbols $z^k$ should be replaced on the functions $f^k(y^j)$. We often denote the inverse matrix $\omega^{-1}$ as $\tilde\omega$. We use the standard convention of summing over repeated indices. Since we are working in local coordinates, all statements should be understood locally, that is, they are true in some vicinity of the point in question.

{\bf Definition 1.2.} Vector space $\mathbb{V}=\{ \vec V, \vec U, \ldots \}$ is called the Lie algebra, if on $\mathbb{V}$ is defined the bilinear mapping $[ {}, {} ]: \mathbb{V}\times \mathbb{V}\rightarrow \mathbb{V}$ (called the Lie bracket), with the properties  
\begin{eqnarray}
\qquad  \quad  ~~[\vec V, \vec U ]=-[\vec U, \vec V] \quad \qquad \qquad \qquad \qquad \qquad \qquad \qquad\mbox{(antisymmetric)}, \qquad   \qquad \label{sim.0.1.1}\\ 
\quad [\vec V, [\vec U, \vec W]]+[\vec U, [\vec W, \vec V]]+[\vec W, [\vec V, \vec U]]\equiv [\vec V, [\vec U, \vec W]]+ cycle=0 \quad \mbox{(Jacobi identity)}. \qquad \qquad \label{sim.0.1.2}
\end{eqnarray}
\noindent Due to the bilinearity, all properties of the Lie bracket are encoded in the Lie brackets of basic vectors $T^i$: $[\vec V, \vec U ]=V_i U_j [T^i, T^j]$. Since $[T^i, T^j]=\vec W\in\mathbb{V}$, we can expand $W$ in the basis $T^i$, obtaining  
\begin{eqnarray}\label{sim.0.1}
[ T^i, T^j]=c^{ij}{}_k T^k,
\end{eqnarray}
where the numbers $c^{ij}{}_k$ are called the structure constants of the algebra in the basis $T^i$. The conditions (\ref{sim.0.1.1}) and (\ref{sim.0.1.2})  are satisfied, if the structure constants obey ({\sc Exercise})
\begin{eqnarray}\label{sim.0.2}
c^{ij}{}_k=-c^{ji}{}_k , \qquad \qquad c^{ij}{}_ac^{ak}{}_b+cycle(i, j, k)=0.
\end{eqnarray}

{\sc Example 1.1.} For the three-dimensional vector space with elements $V=v_i T^i$, $i=1, 2, 3$, let us define $[ T^i, T^j]=\epsilon^{ijk}T^k$, where $\epsilon^{ijk}$ is the Levi-Chivita symbol with $\epsilon^{123}=1$. It can be verified, that the set $c^{ij}{}_k\equiv\epsilon^{ijk}$ has the properties (\ref{sim.0.2}), so the vector space turn into a Lie algebra. It is called the Lie algebra of three-dimensional group of rotations, see Sect. 1.2 in \cite{deriglazov2010classical} for details.

Let $\mathbb{M}_n=\{z, y, \ldots \}$ be $n$\,-dimensional manifold, and $\mathbb{F}_{\mathbb{M}}=\{A, B, \ldots \}$ be space of scalar  functions on $\mathbb{M}_n$, that is the mappings $A: \mathbb{M}_n \rightarrow \mathbb{R}$. Let $z^i$ be local coordinates on  $\mathbb{M}_n$,  that is we have an isomorphism $z\in\mathbb{M}_n\rightarrow z^i(z)\in\mathbb{R}^n$. If $z'^i$ is another coordinate system, we have the relations
\begin{eqnarray}\label{sim.1}
z'^{i}=\varphi^i(z^j), \qquad \qquad z^j=\tilde \varphi^j(z'^i), \qquad \tilde \varphi^i(\varphi^j(z^k))=z^i. 
\end{eqnarray}
Let in the coordinates $z^i$ and $z'^i$  the mapping $A$ is represented by the functions $A(z^i): \mathbb{R}^n \rightarrow \mathbb{R}$ and $A'(z'^i): \mathbb{R}^n \rightarrow \mathbb{R}$. They are related by
\begin{eqnarray}\label{sim.2}
A'(z'^i)=\left. A(z^j)\right|_{z^j\rightarrow\tilde\varphi^j(z'^i)}\equiv A(\tilde\varphi^j(z'^i)). 
\end{eqnarray}
We call (\ref{sim.2}) the transformation law of a scalar function in the passage from $z^i$ to $z'^i$. In certain abuse of terminology, we often said "scalar function $A(z^i)$"  ~ instead of "the function $A(z^i)$ is representative of a scalar function 
$A: \mathbb{M}_n\rightarrow\mathbb{R}$ in the coordinates $z^i$". 

{\sc Example 1.2.}  Scalar function of a coordinate. Given coordinate system $z^i$, define the scalar function $A^1: z\rightarrow z^1$, where  $z^1$ is the first coordinate of the point $z$ in the system $z^i$. In the coordinates $z^i$ the mapping is represented by the following function: $A^1(z^1, z^2, \ldots , z^n)=z^1$. In the coordinates $z'^i=\varphi^i(z^j)$ it is represented by $A'^1(z'^1, z'^2, \ldots , z'^n)=\left. z^1\right|_{z^j\rightarrow\tilde\varphi^j(z'^i)}=\tilde\varphi^1(z'^1, z'^2, \ldots , z'^n)$. 

We often write $z'^i(z^j)$ instead of $\varphi^i(z^j)$, 
$z^j(z'^i)$ instead of $\tilde\varphi^j(z'^i)$, and use the notation $z'^i\equiv z^{i'}$. In the latter case, $i'$ and $i$, when they appear in the same expression, are considered as two different indexes. For instance, in these notations the scalar function of $z^1$\,-coordinate in the system $z^{i'}=z^{i'}(z^j)$ is represented by the function $z^1(z^{i'})$. 

{\sc Exercise 1.1.} Observe that (\ref{sim.1}) implies, that derivatives of the transition functions $\varphi$ and $\tilde\varphi$ form the inverse matrices 
\begin{eqnarray}\label{sim.2.0}
\left.\frac{\partial\tilde\varphi^i}{\partial z'^k}\right|_{z'\rightarrow\varphi(z)} \frac{\partial\varphi^k}{\partial z^j} =\delta^i{}_j \quad \mbox{or, in short notation} \quad 
\frac{\partial z^i}{\partial z^{k'}} \frac{\partial z^{k'}}{\partial z^j} =\delta^i{}_j.
\end{eqnarray}

Given curve $z^i(\tau)\in \mathbb{M}_n$ with $z^i(0)=z^i_0$, the numbers $V^i=\dot z^i(0)$ are called components (coordinates) of tangent vector to the curve at the point $z^i_0$. If $V^{i'}$ are components of the tangent vector in the coordinates $z^{i'}$, we have the relation $V^{i'}=\left.\frac{\partial z^{i'}}{\partial z^i}\right|_{z_0}V^i$. The set of tangent vectors at $z_0$ is $n$\,-dimensional vector space denoted ${\mathbb T}_{{\mathbb M}_n}(z_0)$. 

We said that we have a vector field $\vec V(z)$ on $\mathbb{M}_n$, if in each coordinate system $z^j$ it is defined the set of functions $V^i(z^j)$ with the transformation law
\begin{eqnarray}\label{sim.3.1}
V^{i'}(z^{j'})=\left.\frac{\partial z^{i'}}{\partial z^i}V^i(z^k)\right|_{z\rightarrow z(z')}. 
\end{eqnarray}
The space of all vector fields on $\mathbb{M}_n$ is denoted $\mathbb{T}_{\mathbb{M}_n}$. In the tensor analysis, $\vec V(z)$ is called the contravariant vector field. 

We said that we have covariant vector field $U(z)$ on $\mathbb{M}_n$, if in each coordinate system $z^j$ it is defined the set of functions $U_i(z^j)$ with the transformation law
\begin{eqnarray}\label{sim.3}
U_{i'}(z^{j'})=\left.\frac{\partial z^i}{\partial z^{i'}}U_i(z^k)\right|_{z\rightarrow z(z')}.
\end{eqnarray}
Gradient of a scalar function $A$ is an example of the covariant vector field. Its components are $U_i=\partial_i A$.  

{\sc Exercise 1.2.} Let $A(z^i)=\frac12[(z^1)^2+(z^2)^2+(z^3)^2]$ represent a scalar function in the coordinates $z^i$. Then in the coordinates $z^{i'}$, defined by (\ref{sim.1}), it is represented by $A'(z^{i'})=\frac12[(\tilde\varphi^1(z^{i'}))^2+(\tilde\varphi^2(z^{i'}))^2+(\tilde\varphi^3(z^{i'}))^2]$. Gradients of these functions are $U_i=z^i$ and $U_{i'}=\tilde\varphi^j(z^{i'})\frac{\partial\tilde\varphi^j(z^{i'})}{\partial z^{i'}}$. Confirm that the two gradients are related by Eq. (\ref{sim.3}). 

Similarly to this, contravariant tensor of second-rank is a quantity with the transformation law
\begin{eqnarray}\label{sim.4}
\omega^{i'j'}(z^{k'})=\left.\frac{\partial z^{i'}}{\partial z^i}\frac{\partial z^{j'}}{\partial z^j}\omega^{ij}(z^m)\right|_{z\rightarrow z(z')}, 
\end{eqnarray} 
and so on.

\noindent {\sc Exercise 1.3.} Contraction of $\omega$ with covariant vector field $grad ~ A$ gives a quantity with the 
components $V^i=\omega^{ij}\partial_j A$.  Confirm that $\vec V$ is a contravariant vector field. 

Integral line of the vector field $V^i(z^k)$ on $\mathbb{M}_n$ is a solution $z^i(\tau)$ to the system $\frac{d z^i(\tau)}{d\tau}=V^i(z^k(\tau))$.  We assume that $V^i(z^k)$ is a smooth field, so through each point of the manifold passes unique integral line of $\vec V$.

{\bf Submanifold of $\mathbb{M}_n$}.  The $k$\,-dimensional submanifold $\mathbb{N}_k^{\vec c}\in\mathbb{M}_n$ is often defined as a constant-level surface of a set of functionally independent scalar functions  $\Phi^\alpha(z)$
\begin{eqnarray}\label{sim.2.1}
\mathbb{N}_k^{\vec c}=\{ z\in \mathbb{M}_n, ~  \Phi^\alpha(z^k)=c^\alpha, ~  \alpha=1, 2, \ldots n-k\}, 
\end{eqnarray}
where $c^\alpha$ are given numbers. 

We recall that scalar functions $\Phi^\alpha(z)$, $\alpha=1, 2, \ldots , n-k$ are called functionally independent, if for their 
representatives $\Phi^\alpha(z^i)$ in the coordinates $z^i$ we have: $rank ~(\partial_i \Phi^\alpha)=n-k$. This implies, that covariant vectors $V_{(\alpha)}$ with coordinates $V_{(\alpha)i}=\partial_i \Phi^\alpha$ are linearly independent. The equations $\Phi^\alpha(z^i)=c^\alpha$  for the functionally independent functions  can be 
resolved: $z^\alpha=f^\alpha(z^a)$, $a=1, 2, \ldots , k$.  So the coordinates $z^i$ are naturally divided on two groups: $(z^\alpha, z^a)$, and $z^a$, $a=1, 2, \ldots , k$,  can be taken as local coordinates of the submanifold $\mathbb{N}_k^{\vec c}$. Below we always assume that the coordinates have been grouped in this way, and $\det\frac{\partial \Phi^\alpha}{\partial z^\beta}\ne 0$. 

If we have only one function $\Phi(z)$, we assume that it has a non vanishing gradient, $rank ~(\partial_i \Phi)=1$. 

Taking $c^\alpha=0$ in (\ref{sim.2.1}), we have the surface of level zero  
\begin{eqnarray}\label{sim.2.1.0}
\mathbb{N}_k=\{ z\in \mathbb{M}_n, ~  \Phi^\alpha(z^k)=0,   ~ ~ \alpha=1, 2, \ldots n-k\}. 
\end{eqnarray}

Let us introduce the notions that will be useful in discussing the Frobenius theorem (see Appendix C). \label{tang}

For the curve $z^i(\tau)\subset{\mathbb N}_k\subset{\mathbb M}_n$, with $z^i(0)=z^i_0$, the tangent vector $V^i(z_0)=\frac{dz^i(0)}{d\tau}\in{\mathbb T}_{\mathbb M}(z_0)$ is called tangent vector to ${\mathbb N}_k$ at $z_0$. The set of all tangent vectors at $z_0$ is $k$\,-dimensional vector space denoted 
${\mathbb T}_{\mathbb N}(z_0)$. For any such vector holds\footnote{In three-dimensional Euclidean space this equality has simple geometric meaning:  vector $grad ~ F(x, y, z)$ in $\mathbb{R}^3$ is orthogonal to the surfaces of level $F(x, y, z)=c$ of the scalar function $F(x, y, z)$.} the equality $V^i\partial_i\Phi^\alpha(z_0)=0$. 

Vector field $V^i(z)$ on ${\mathbb M}_n$ is tangent to  ${\mathbb N}_k$, if any integral curve of $V^i(z)$ crossing ${\mathbb N}_k$, lie entirely in ${\mathbb N}_k$: $\Phi^\alpha(z^k(0))=0$ implies $\Phi^\alpha(z^k(\tau))=0$ for any $\tau$. Vector field $V^i(z)$ on ${\mathbb M}_n$ touches the surface  ${\mathbb N}_k$, if $V^i\partial_i\Phi^\alpha|_{z_0}=0$ for 
any $z_0\in{\mathbb N}_k$. The tangent field touches the surface. The converse is not true.

{\bf Foliation of $\mathbb{M}_n$}. The set $\{ \mathbb{N}_k^{\vec c}, ~ \vec c\in \mathbb{R}^{n-k}\}$ of the submanifolds (\ref{sim.2.1}) is called a foliation of $\mathbb{M}_n$, while $\mathbb{N}_k^{\vec c}$ are called leaves of the foliation.  Notice that submanifolds with different $\vec c$ do not intercept, and any\footnote{Recall that all our assertions hold locally.} $z\in\mathbb{M}_n$ lies in one 
of $\mathbb{N}_k^{\vec c}$.  

There are coordinates, naturally adapted with the foliation: $z^k\rightarrow y^k=(y^\alpha, y^a)$, with the transition functions $y^a=z^a$, $y^\alpha=\Phi^\alpha(z^\beta, z^b)$.  In these coordinates the sumanifolds $\mathbb{N}_k^{\vec c}$ look like hyperplanes:
\begin{eqnarray}\label{sim.2.2}
\mathbb{N}_k^{\vec c}=\{y^i\in\mathbb{M}_n, ~ y^\alpha=c^\alpha \}, 
\end{eqnarray} 
and $y^a=z^a$ can be taken as local coordinates of $\mathbb{N}_k^{\vec c}$. The useful identity is 
\begin{eqnarray}\label{sim.2.2.0}
\left. A(z^i(y^j))\right|_{y^\alpha=0}=\left. A(f^\alpha(z^a), z^a)\right|_{z^a\rightarrow y^a}. 
\end{eqnarray}

{\bf Lie bracket (commutator) of vector fields} is bilinear operation $[{~}, {~}]:  \mathbb{T}_{\mathbb{M}_n}\times\mathbb{T}_{\mathbb{M}_n}\rightarrow \mathbb{T}_{\mathbb{M}_n}$, that with each pair of vector fields  $\vec V$ and $\vec U$ of $\mathbb{T}_{\mathbb{M}_n}$ associates the vector field $[\vec V, \vec U]$ of $\mathbb{T}_{\mathbb{M}_n}$ according to the rule
\begin{eqnarray}\label{sim.2.3}
[\vec V, \vec U]^i=V^j\partial_j U^i-U^j\partial_j V^i. 
\end{eqnarray} 
The quantity $[\vec V, \vec U]^i$ is indeed a vector field, which can be verified by direct computation. We have $[\vec V, \vec U]^i=V^{j'}\partial_{j'}(\frac{\partial z^i}{\partial z^{i'}}U^{i'})-(V\leftrightarrow U)=\frac{\partial z^i}{\partial z^{i'}} (V^{j'}\partial_{j'} U^{i'}-(V\leftrightarrow U))=\frac{\partial z^i}{\partial z^{i'}}[\vec V', \vec U']^{i'}$, in agreement with Eq. (\ref{sim.3.1}). The Lie bracket has the properties (\ref{sim.0.1.1}) and (\ref{sim.0.1.2}), and turns the space of vector fields into infinito-dimensional Lie algebra. 

Each vector field  determines a linear mapping $\vec V: \mathbb{F}_{\mathbb{M}_n}\rightarrow \mathbb{F}_{\mathbb{M}_n}$ on the space of scalar functions  according to the rule 
\begin{eqnarray}\label{sim.2.4}
\vec V: A \rightarrow \vec V(A)=V^i\partial_i A. 
\end{eqnarray} 
Notice that $\vec V(A)=0$ for all $A$ implies $V^i=0$.  Then the Lie bracket can be considered as a commutator of two differential operators
\begin{eqnarray}\label{sim.2.5}
[\vec V, \vec U](A)=\vec V(\vec U(A))- \vec U(\vec V(A)). 
\end{eqnarray} 
Using this formula, it is easy to confirm by direct computation the Jacobi identity (\ref{sim.0.1.2}) for the Lie bracket (\ref{sim.2.3}).

\subsection{The mapping of manifolds and induced mappings of tensor fields.}\label{I D}

Given two manifolds $\mathbb{N}_k=\{x^a\}$, $\mathbb{M}_n=\{z^i\}$, consider the functions $z^i=\phi^i(x^a)$. They determine the mapping 
\begin{eqnarray}\label{sim.2.6}
\phi: \mathbb{N}_k \rightarrow\mathbb{M}_n, \qquad x^a\rightarrow z^i=\phi^i(x^a)\equiv z^i(x^a). 
\end{eqnarray} 
If $\phi$ is an injective function: $rank ~ \frac{\partial\phi^i}{\partial x^a}=k$, the image of the mapping is $k$\,-dimensional submanifold of $\mathbb{M}_n$: $\mathbb{N}_k=\{z^i\in\mathbb{M}_n, ~ z^\alpha-f^\alpha(z^a)=0\}$, where the equalities $z^\alpha=f^\alpha(z^a)$ are obtained  excluding $x^a$ from the equations $z^i=\phi^i(x^a)$. In some cases \cite{Olver}, the manifold $\mathbb{N}_k$ can be identified with this submanifold of $\mathbb{M}_n$. 

Conversely, let $\mathbb{N}_k\subset\mathbb{M}_n$, then the parametric equations $z^\alpha=f^\alpha(z^b)$ of the submanifold (\ref{sim.2.1.0}) can be considered as determining the mapping of embedding 
\begin{eqnarray}\label{sim.2.7}
\eta: \mathbb{N}_k=\{z^b\} \rightarrow \mathbb{M}_n=\{z^i\},  \quad z^b\rightarrow z^i=(z^\alpha, z^b), \quad \mbox{where} \quad z^\alpha=f^\alpha(z^b). 
\end{eqnarray} 
Using the mapping (\ref{sim.2.6}), some geometric objects from one manifold can be transferred to another.  We start from the spaces of covariant  and contravariant tensors {\it at the points} $x_0$ and $z_0=\phi(x_0)$. Take, for definiteness, the second-rank tensors. Given  $U_{ij}(z_0)$, we can construct the  induced tensor  $U_{ab}(x_0)$ 
\begin{eqnarray}\label{sim.2.8}
\mathbb{T}^{(0, 2)}_{\mathbb{M}} \rightarrow \mathbb{T}^{(0, 2)}_{\mathbb{N}}, \qquad 
U_{ij}(z_0) \rightarrow U_{ab}(x_0)=\frac{\partial z^i(x_0)}{\partial x^a}\frac{\partial z^j(x_0)}{\partial x^b}U_{ij}(z_0). 
\end{eqnarray} 
Given  $V^{ab}(x_0)$, we can construct the  induced tensor  $U^{ij}(z_0)$ 
\begin{eqnarray}\label{sim.2.9}
\mathbb{T}^{(2, 0)}_{\mathbb{N}} \rightarrow \mathbb{T}^{(2, 0)}_{\mathbb{M}}, \qquad
V^{ab}(x_0) \rightarrow V^{ij}(z_0)=\frac{\partial z^i(x_0)}{\partial x^a}\frac{\partial z^j(x_0)}{\partial x^b}V^{ab}(x_0). 
\end{eqnarray} 
For the case of vector, the notion of induced mapping, 
\begin{eqnarray}\label{sim.2.10}
V^{i}(z_0)=\frac{\partial z^i(x_0)}{\partial x^a}V^{a}(x_0), 
\end{eqnarray} 
is consistent with the notion of a tangent vector: if $V^{a}=\frac{dx^a}{d\tau}$ is tangent vector to the curve $x^a(\tau)$, then $V^{i}$, given by (\ref{sim.2.10}), is tangent vector to the 
image $z^i(x^a(\tau))$ 
\begin{eqnarray}\label{sim.2.11}
V^{i}(z(\tau))=\frac{d}{d\tau}z^i(x^a(\tau)). 
\end{eqnarray} 

Concerning the fields on the manifolds, $\phi$ naturally  induces the mappings of scalar functions and of covariant tensor fields. For the functions the induced mapping 
\begin{eqnarray}\label{sim.2.12}
\phi^*: \mathbb{F}_{\mathbb{M}_n}\rightarrow \mathbb{F}_{\mathbb{N}_k}, \qquad A(z^i)\rightarrow \bar A(x^a)=A(z^i(x^a)), 
\end{eqnarray} 
is just the composition: $\bar A=A \circ \phi$. For the covariant tensor fields we have
\begin{eqnarray}\label{sim.2.13}
\phi_*: \mathbb{T}^{(0, 2)}_{\mathbb{M}_n}\rightarrow \mathbb{T}^{(0, 2)}_{\mathbb{N}_k}, \qquad
U_{ij}(z^i) \rightarrow U_{ab}(x^a)=\frac{\partial z^i}{\partial x^a}\frac{\partial z^j}{\partial x^b}U_{ij}(z^i(x^a)).  
\end{eqnarray} 
Notice that the contravariant tensor fields can not be transferred to another manifold (submanifold) in this manner. As we will see in the next section, the Poisson structure on ${\mathbb{M}_n}$ is determined namely by secod-rank contravariant tensor. Hence it can not be directly transferred on a submanifold. This turns out to be possible in special case of Casimir submanifolds (see Sect. \ref{IV B}), and leads to the Dirac bracket (see Sect. \ref{VI B}). 

Note also that if two contravariant fields on the manifolds ${\mathbb{M}_n}$ and ${\mathbb{N}_k}$ are given, we can of course compare them using the mapping (\ref{sim.2.6}), see Eq. (\ref{sim.15.3}) below as an example.

\section{Poisson manifold.}\label{SymPoi}
Let on the space of functions $\mathbb{F}_{\mathbb{M}}$ is defined a bilinear mapping $\{ {}, {}\}: \mathbb{F}_{\mathbb{M}}\times \mathbb{F}_{\mathbb{M}}\rightarrow \mathbb{F}_{\mathbb{M}}$ (called the Poisson bracket), with the properties 
\begin{eqnarray}\label{sim.5}
\qquad \qquad \qquad \quad   \{A, B\}=-\{B, A\} \qquad \mbox{(antisymmetric)}, \label{sim.A5}\\ 
\qquad \qquad ~\{A, \{B, C\}\}+ cycle=0 \qquad \mbox{(Jacobi identity)}, \label{sim.B5}\\
\qquad \{A, BC\}=\{A, B\}C+\{A, C\}B \qquad \mbox{(Leibnitz rule)}.  \label{sim.C5}
 \end{eqnarray}
\noindent  When $\mathbb{F}_{\mathbb{M}}$ is equipped with the Poisson bracket, the manifold $\mathbb{M}_n$ is called the Poisson manifold. Comparing  (\ref{sim.A5}) and (\ref{sim.B5}) with (\ref{sim.0.1.1}) and (\ref{sim.0.1.2}), we see that the infinite-dimensional vector space $\mathbb{F}_{\mathbb{M}}$ is equipped with the structure of a Lie algebra.  \par

\noindent {\sc Exercise 2.1.} Show that constant functions, $A(z)=c$ for any $z$, have vanishing brackets (commute) with all other functions. 

One of the ways to define the Poisson structure on $\mathbb{M}_n$ is as follows. \par

\noindent {\bf Affirmation 2.1.}  Let $\omega^{ij}(z)$ be the contravariant tensor of second rank on $\mathbb{M}_n$. The mapping 
\begin{eqnarray}\label{sim.5}
\{A, B\}=\partial_i A ~ \omega^{ij} ~ \partial_j B, 
\end{eqnarray}
determines the Poisson bracket, if the tensor $\omega$ obeys the properties 
\begin{eqnarray}
\qquad \qquad  \qquad  \omega^{ij}=-\omega^{ji} \qquad \mbox{(antisymmetric)}, \label{sim.D5} \\
\qquad \qquad \qquad    \omega^{ip}\partial_p\omega^{jk}+ cycle(i, j, k)=0. \qquad\label{sim.E5}
\end{eqnarray}
In particular, each numeric antisymmetric matrix determines a Poisson bracket. We call $\omega$ the Poisson tensor. \par 

\noindent {\bf Proof.} First, we note that (\ref{sim.D5}) implies (\ref{sim.A5}). Second, the mapping (\ref{sim.5}), being combination of derivatives, is bilinear and automatically obeys the Leibnitz rule. To complete the proof, we need to show that (\ref{sim.E5}) is equivalent to (\ref{sim.B5}). Using (\ref{sim.E5}), by direct calculation we get 
\begin{eqnarray}\label{sim.6}
\{A, \{B, C\}\}+ cycle(A, B, C)=\partial_i A\partial_j B\partial_k C\omega^{ip}\partial_p\omega^{jk}+cycle(A, B, C)+\omega^{ip}\omega^{jk}\partial_p\left[\partial_i A\partial_j B\partial_k C\right]+cycle(A, B, C)
\end{eqnarray}
By direct calculation, we can show also that in the first term on r.h.s. the $cycle(A, B, C)$ is equivalent to $cycle(i, j, k)$. So we write the previous equality as
\begin{eqnarray}\label{sim.8}
\{A, \{B, C\}\}+ cycle(A, B, C)=\partial_i A\partial_j B\partial_k C\left[\omega^{ip}\partial_p\omega^{jk}+cycle(i, j, k)\right]+   \cr  
\omega^{ip}\omega^{jk}\partial_p\left[\partial_i A\partial_j B\partial_k C+\partial_j A\partial_k B\partial_i C+\partial_k A\partial_i B\partial_j C \right]. 
\end{eqnarray}
The second line in this equality identically vanishes due to symmetry properties of this term. Indeed, we write the first term of the line as follows:
\begin{eqnarray}\label{sim.9}
\omega^{ip}\omega^{jk}\partial_p\left[\partial_i A\partial_j B\partial_k C\right]=\omega^{ip}\omega^{jk}\partial_p\partial_i \left[A\partial_j B\partial_k C\right]-\omega^{ip}\omega^{jk}\partial_p\left[A\partial_i \partial_j B\partial_k C+ A\partial_j B\partial_i\partial_k C\right]= \cr
\omega^{ip}\omega^{jk}\partial_p\left[A\partial_k(\partial_i B\partial_j C-[i\leftrightarrow j])\right]. \qquad \qquad \qquad \qquad 
\end{eqnarray}
The two remaining terms of the line we write as
\begin{eqnarray}\label{sim.10}
\omega^{ip}\omega^{jk}\partial_p\left[\partial_j A\partial_k B\partial_i C+\partial_k A\partial_i B\partial_j C\right]
=\omega^{ip}\omega^{jk}\partial_p\left[\partial_kA(\partial_i B\partial_j C-[i\leftrightarrow j])\right]= \cr 
\omega^{ip}\omega^{jk}\partial_p\partial_k\left[A(\partial_i B\partial_j C-[i\leftrightarrow j])\right]-\omega^{ip}\omega^{jk}\partial_p\left[A\partial_k(\partial_i B\partial_j C-[i\leftrightarrow j])\right]. 
\end{eqnarray}
The last terms in (\ref{sim.9}) and (\ref{sim.10}) cancel each other, while the first term in (\ref{sim.10}) iz zero, being the trace of the product of symmetric 
$D^{ij}\equiv \omega^{ip}\omega^{jk}\partial_p\partial_k$ and antisymmetric $E_{ij}\equiv A(\partial_i B\partial_j C-[i\leftrightarrow j])$ quantities. Thus we have obtained the identity
\begin{eqnarray}\label{sim.11} 
\omega^{ip}\omega^{jk}\partial_p\left[\partial_i A\partial_j B\partial_k C+\partial_j A\partial_k B\partial_i C+\partial_k A\partial_i B\partial_j C \right]=0. 
\end{eqnarray}
Taking this into account in (\ref{sim.8}), we see the equivalence of the conditions (\ref{sim.B5}) and (\ref{sim.E5}). $\blacksquare$  \par

\noindent {\bf Affirmation 2.2.} Let the bracket (\ref{sim.5}) obeys the Jacobi identity in the coordinates $z^i$. Then the Jacobi identity is satisfied in any other coordinates. 

This is an immediate consequence of tensor character of involved quantities. Indeed, the bracket $\{A, B\}=\partial_i A \omega^{ij} \partial_j B$ is a contraction of three tensors and so is a scalar function under diffeomorphisms. Then the same is true for $\{A, \{ B, C\}\}$. Let us denote l.h.s. of the Jacobi identity as $D(z)$. Then the Jacobi identity is the coordinate-independent statement that the scalar function $D(z)$ identically vanishes for all $z\in\mathbb{M}_n$. This can be verified also by direct computations, see Appendix A. As a consequence, the left hand side of Eq. (\ref{sim.E5}) is a tensor of third rank\footnote{This is a non trivial affirmation, since $\partial_p\omega^{jk}$ is not a covariant object.}.

For the scalar functions of coordinates (see Example 1.2), the Poisson bracket (\ref{sim.5}) reads
\begin{eqnarray}\label{sim.11.4}
\{z^i, z^j\}=\omega^{ij}. 
\end{eqnarray}
In classical mechanics these equalities are known as fundamental brackets of the coordinates. Observe that the identity (\ref{sim.E5}) can be written as follows: $\{z^i, \{z^j, z^k\}+cycle (i,j,k)=0$.

The bracket (\ref{sim.5}) is called nondegenerate if $\det\omega\ne 0$, and degenerate when $\det\omega=0$. Examples will be presented below: (\ref{sim.15}) is nondegenerate while (\ref{sim.16.1}) is degenerate. The structure of the matrix $\omega$ depends on its rank, and becomes clear in the so called canonical coordinates specified by the following theorem: \par 

\noindent {\bf Generalized Darboux theorem.} Let $rank~\omega=2k$ at the point $z^i\in\mathbb{M}_n$. Then there are local coordinates 
$z^{i'}=(z^{\beta'}, z^{a'})$, $a'=1, 2. \ldots , 2k$, $\beta'=1, 2, \ldots , p=n-2k$ such that  $\omega$ in some vicinity of $z^i$ has the form: 
\begin{eqnarray}\label{dar.1.0}
\omega^{i'j'}=\left(
\begin{array}{ccc}
0_{p\times p}  & 0 & 0 \\
0  & 0_{k\times k} & 1_{k\times k} \\
0  & -1_{k\times k} & 0_{k\times k} 
\end{array}\right),
\end{eqnarray}
or
\begin{eqnarray}\label{dar.1}
\omega^{a'b'}=\left(
\begin{array}{cc}
0 & 1 \\
-1 & 0
\end{array}
\right), \qquad \omega^{\beta' j'}= \omega^{j'\beta'}=0, \quad \mbox{where} \quad j'=1, 2, \ldots , n.   
\end{eqnarray}

Proof is given in Appendix B. 
We recall that determinant of any odd-dimensional antisymmetric matrix vanishes, this implies that $rank ~ \omega$ is necessary an even number, as it is written above. Let us further denote $z^{a'}=(q^1, q^2, \ldots , q^k, p_1, p_2, \ldots , p_k)$.
Then, in terms of fundamental brackets, the equalities  (\ref{dar.1}) can be written as follows:
\begin{eqnarray}\label{dar.2}
\{q^{a'}, p_{b'}\}=\delta^{a'}{}_{b'}, \qquad  \{q^{a'}, q^{b'}\}=0,  \qquad \{p_{a'}, p_{b'}\}=0,     \qquad \{z^{j'}, z^{\beta'}\}= 0. 
\end{eqnarray}

\section{Hamiltonian dynamical systems on a Poisson manifold.}\label{DynPoi}
\subsection{Hamiltonian vector fields.} 

Using the Poisson structure (\ref{sim.5}), with each function $H(z^i)\in\mathbb{F}_{\mathbb{M}}$ we can associate the contravariant vector field $X^i_H\equiv \omega^{ij}\partial_j H=\{z^i, H\}\in\mathbb{T}_{\mathbb{M}} $. That is we have the mapping
\begin{eqnarray}\label{sim.12.0}
\omega: \mathbb{F}_{\mathbb{M}}\rightarrow \mathbb{T}_{\mathbb{M}}, \qquad 
\omega: H \rightarrow [\omega(H)]^i=\omega^{ij}\partial_j H, \quad \mbox{we also denote} \quad \omega(H)\equiv \vec X_H\in\mathbb{T}_{\mathbb{M}}. 
\end{eqnarray}
$\vec X_H$ is called the Hamiltonian vector field of the function $H$.
Then 
\begin{eqnarray}\label{sim.12}
\dot z^i=\{z^i, H\}\equiv\omega^{ij}\partial_j H, 
\end{eqnarray}
are called Hamiltonian equations, the scalar function $H$ is called the Hamiltonian.  Solutions $z^i(\tau)$ of the equations are called integral lines of the vector field $\{z^i, H\}$ created by $H$ on $\mathbb{M}_n$. We assume that $\vec X_H$ is a smooth vector field, so the Cauchy problem for (\ref{sim.12}) has unique solution in a vicinity of  any point  of $\mathbb{M}_n$. $\vec X_H$ at each point is tangent vector to the integral line that passes through this point. 

Let $z^k(\tau)$ be integral line of $\vec X_A$ and $B$ be scalar function. Then we can write
\begin{eqnarray}\label{sim.12.001}
\frac{d}{d\tau}B(z^k(\tau))=\left. \{B, A\}\right|_{z(\tau)}.
\end{eqnarray}
Using this equality, and the fact that integral lines pass through each point of $\mathbb{M}_n$, it is easy to prove the three affirmations presented below. They will be repeatedly used (and sometimes rephrased, see Sect. V) in our subsequent considerations.

\noindent {\bf Affirmation 3.1.}  Integral line of $\vec X_H$ entirely lies on one of the surfaces $H(z^k)=c=const$. In classical mechanics it is just the energy conservation law. 


Denote $\vec X_{(j)}$ the Hamiltonian vector field associated with scalar function of the coordinate $z^j$. Its components are $X_{(j)}^i=\omega^{ik}\partial_k z^j=\omega^{ij}$. Hence the Poisson matrix can be considered\footnote{Notice that it is an example of coordinate-dependent statement.} as composed of the columns $\vec X_{(j)}$
\begin{eqnarray}\label{sim.12.00}
\omega=(\vec X_{(1)}, \vec X_{(2)}, \ldots , \vec X_{(n)}).
\end{eqnarray}
According to the Affirmation 3.1, integral lines of the vector $\vec X_{(j)}$ lie on the hyperplanes $z^j=const$. 

\noindent {\bf Affirmation 3.2.} Given scalar functions $H$ and $Q^\alpha$, $\alpha=1, 2, \ldots , n-k$, the following two conditions are equivalent:

\noindent {\bf (A)}  Integral lines of $\vec X_H$ lie in the submanifolds $\mathbb{N}_k^{\vec c}=
\{z\in\mathbb{M}_n, \quad Q^\alpha=c^\alpha, H=c\}$.

\noindent {\bf (B)} All $Q^\alpha$ commute with $H$: $\{Q^\alpha, H\}=0$, for all  $z\in\mathbb{M}_n$.

In classical mechanics the quantities $Q^\alpha$ are called first integrals (or conserved charges) of the system.

%

\noindent {\bf Affirmation 3.3.} Let $A^\alpha$, $\alpha=1, 2, \ldots , n-k$ be functionally independent scalar functions, and denote $\vec V_{(\alpha)}$ the Hamiltonian field of $A^\alpha$. The following two conditions are equivalent:

\noindent {\bf (A)} Integral lines of each $\vec V_{(\beta)}$ lie in the submanifolds $\mathbb{N}_k^{\vec c}=
\{z\in\mathbb{M}_n, \quad A^\alpha=c^\alpha \}$.

\noindent {\bf (B)} $\{A^\alpha, A^\beta \}=0$ on $\mathbb{M}_n$.

\subsection{Lie bracket and Poisson bracket.} 

Consider the spaces of scalar functions and of vector fields on $\mathbb{M}_n$, which are the infinito-dimensional Lie algebras: $\mathbb{F}_{\mathbb{M}}=\{ A, B, \ldots , ~ \{{~}, {~}\}~\}$ and $\mathbb{T}_{\mathbb{M}}=\{ \vec V, \vec U, \ldots , ~ [{~}, {~}]~\}$. 

\noindent {\bf Affirmation 3.3.} The mapping (\ref{sim.12.0}) respects the Lie products of $\mathbb{F}_{\mathbb{M}}$ and $\mathbb{T}_{\mathbb{M}}$:
\begin{eqnarray}\label{sim.12.3}
\omega(\{ A, B \})=-[\omega(A), \omega(B)], \quad \mbox{or, equivalently} \quad \vec X_{\{A, B\}} =-[\vec X_A, \vec X_B ].  
\end{eqnarray}
According to the last equality,  the Hamiltonian vector fields form a subalgebra of the Lie algebra $\mathbb{T}_{\mathbb{M}}$. \par

\noindent {\bf Proof.}  Using the vector notation (\ref{sim.2.4}), we can present the Poisson bracket as follows:
\begin{eqnarray}\label{sim.12.4}
\{ A, B \}=-\vec X_A(B). 
\end{eqnarray}
The equality (\ref{sim.12.3}) is the Jacobi identity (\ref{sim.C5}), rewritten in the vector notations. Indeed 
\begin{eqnarray}\label{sim.12.5}
\{\{ A, B\}, C \}=\{ A, \{B, C \}\}-\{ B, \{A , C \}\}, \quad \mbox{or} \quad \vec X_{\{A, B\}}(C)=\vec X_A(\vec X_B(C))-\vec X_B(\vec X_A(C)), 
\end{eqnarray}
for all $C$, which is just (\ref{sim.12.3}).  $\blacksquare$

We also note that in the vector notation, the Jacobi identity (\ref{sim.E5}) states that Hamiltonian fields of coordinates form the closed algebra
\begin{eqnarray}\label{sim.12.6}
[\vec X_{(i)},  \vec X_{(j)}]=c_{(i)(j)}{}^{(k)}\vec X_{(k)}, 
\end{eqnarray}
with the structure functions $c_{(i)(j)}{}^{(k)}=-\partial_k\omega^{ij}$.  

{\sc Exercise 3.1.} Show that $\{Q, H\}=const$ implies $[\vec X_Q, \vec X_H]=0$.

\subsection{Two basic examples of Poisson structures.} \par 

\noindent {\bf 1.} Consider the space $\mathbb{R}^{2n}$, denote its coordinates $z^i=(q^1, q^2, \ldots, q^n, p_1, p_2, \ldots, p_n)\equiv (q^a, p_b)$, $a, b=1, 2, \ldots,  n$, and take the matrix composed from four $n\times n$ blocks as follows:
\begin{eqnarray}\label{sim.14}
\omega^{ij}=\left(
\begin{array}{cc}
0 & 1 \\
-1 & 0
\end{array}
\right).
\end{eqnarray}
In all other coordinate systems $z^{i'}$, we define components of the matrix $\omega^{i'j'}$ according to Eq. (\ref{sim.4}). Then $\omega$ is the contravariant tensor of second rank, which (in the system $z$) determines the Poisson structure on $\mathbb{R}^{2n}$ according to Eq. (\ref{sim.5}):  
\begin{eqnarray}\label{sim.15}
\{ A, B\}_P=\frac{\partial A}{\partial q^a}\frac{\partial B}{\partial p_a}-\frac{\partial B}{\partial q^a}\frac{\partial A}{\partial p_a}, \quad \mbox{and fundamental (nonvanishing) brackets are:} \quad \{q^a, p_b \}_P=\delta^a{}_b.  
\end{eqnarray}
As $\omega$ is the numeric matrix, the condition (\ref{sim.E5}) is satisfied in the coordinate system $(q^a, p_b)$. According to Affirmation 2.2, it is then satisfied in all other coordinates. Given Hamiltonian function $H$, the Hamiltonian equations acquire the following form:
\begin{eqnarray}\label{sim.15.1}
\dot q^a=\{ q^a, H \}_P=\frac{\partial H}{\partial p_a}, \qquad \qquad \dot p_a=\{ p_a, H \}_P=-\frac{\partial H}{\partial q^a}.
\end{eqnarray}
It is known (see Sect. 2.9 in \cite{deriglazov2010classical}) that they  follow from the variational problem for the functional 
\begin{eqnarray}\label{sim.16}
S_H: (q^a(\tau), p_a(\tau) )\rightarrow \mathbb{R}; \qquad S_H=\int_{\tau_1}^{\tau_2} d\tau \left[p_a\dot q^a-H(q^a, p_b)\right].
\end{eqnarray}
In classical mechanics, $\mathbb{R}^{2n}$ equipped with the coordinates $(q^a, p_b)$ is called the phase space, the bracket (\ref{sim.15}) is called the canonical Poisson bracket, while the functional $S_H$ is called the Hamiltonian action.  \par 

\noindent {\bf 2.} Given manifold $\mathbb{M}_n$, let $c^{ij}{}_k$ be structure constants of  an $n$\,-dimensional  Lie algebra. We define $\omega^{ij}(z)=c^{ij}{}_k z^k$. Then the equalities (\ref{sim.0.2}) imply  (\ref{sim.D5}) and (\ref{sim.E5}), so the tensor $\omega^{ij}$ determines a Poisson structure on $\mathbb{M}_n$. The corresponding bracket 
\begin{eqnarray}\label{sim.16.0}
\{ A, B \}_{LP}=\partial_i A c^{ij}{}_k z^k\partial_j B, \quad \mbox{fundamental brackets:} \quad \{z^i, z^j\}_{LP}=c^{ij}{}_k z^k, 
\end{eqnarray}
is called the Lie-Poisson bracket. 
In particular, the Lie algebra of rotations determines the Lie-Poisson bracket on $\mathbb{R}^3$
\begin{eqnarray}\label{sim.16.1}
\omega^{ij}=\epsilon^{ijk}z^k. 
\end{eqnarray}
Let $B^i$ are coordinates of a constant vector ${\bf B}\in\mathbb{R}^3$. Taking $H=z^i B^i$ as the Hamiltonian, we obtain the Hamiltonian equations (called the equations of precession) 
\begin{eqnarray}\label{sim.16.2}
\dot z^i=\epsilon^{ijk}B^j z^k, \quad \mbox{or} \quad \dot{\bf z}={\bf B}\times {\bf z},
\end{eqnarray}
where ${\bf B}\times{\bf z}$ is the usual vector product in $\mathbb{R}^3$. For any solution $\bf z(\tau)$, the end of this vector lies in a plane 
perpendicular to ${\bf B}$, and describes a circle around ${\bf B}$, with an angular velocity equal to the magnitude $|{\bf B}|$ of this vector. The compass needle in the earth's magnetic field moves just according to this law.

\subsection{Poisson mapping and Poisson submanifold.} 

Here we discuss the mappings which are compatible with Poisson brackets of the involved manifolds. Intuitively, such a mapping turns the bracket of one manifold into the bracket of another. As an instructive example, we first consider 
the manifolds with the brackets (\ref{sim.14}) and (\ref{sim.16.0}).  Introduce the mapping 
\begin{eqnarray}\label{sim.15.1}
\phi: \mathbb{R}^{2n} \rightarrow \mathbb{M}^{n}, \quad (q^a, p_b) \rightarrow  z^a=\phi^a(q, p)=-c^{ab}{}_c p_b q^c.  
\end{eqnarray}
Computing the canonical Poisson bracket (\ref{sim.15}) of the functions $\phi^a(q, p)$, 
we obtain a remarkable relation between the two brackets ({\sc Exercise}): 
\begin{eqnarray}\label{sim.15.2}
\{\phi^a(q, p), \phi^b(q, p)\}_P=c^{ab}{}_c\phi^c(q, p), 
\quad \mbox{or}  \quad \{\phi^a(q, p), \phi^b(q, p)\}_P=\left.\{z^a,  z^b\}_{LP}\right|_{ z\rightarrow\phi(q, p)}, \quad \mbox{or}  \quad 
\end{eqnarray}
\begin{eqnarray}\label{sim.15.3}
\partial_i\phi^a\omega^{ij}\partial_j\phi^b=\left.\omega^{ab}(z^c)\right|_{ z\rightarrow\phi(q, p)},  \quad \mbox{where} \quad \omega^{ab}(z^c)=c^{ab}{}_cz^c. 
\end{eqnarray}
The relation (\ref{sim.15.3}) shows that Poisson structures  $\omega^{ij}$ and $\omega^{ab}$ are related by the tensor-like law (\ref{sim.4}). 
The relations (\ref{sim.15.2}) show that Poisson brackets of the special functions $\phi^a$ on $\mathbb{R}^{2n}$ are the same as fundamental Lie-Poisson brackets (\ref{sim.16.0}) of the manifold $\mathbb{R}^{n}$.    
We can made these relations to hold for an arbitrary scalar functions, by using  the induced mapping  between the functions $A(z^a)$ of $\mathbb{M}^{n}$ and $\bar A(q^a, p_b)$ of $\mathbb{R}^{2n}$
\begin{eqnarray}\label{sim.15.4}
\phi^*: A(z^a) \rightarrow \bar A(q^a, p_b)\equiv \phi^*(A)(q, p)=A(\phi^a(q, p)).
\end{eqnarray}
This implies the following relation between the Poisson and Lie-Poisson brackets  ({\sc Exercise}):
\begin{eqnarray}\label{sim.15.5}
\{\phi^*(A), \phi^*(B)\}_P=\phi^*\left(\{A, B\}_{LP}\right).  
\end{eqnarray}

Formalizing this example, we arrive at the notion of a Poisson mapping. 

\noindent {\bf Definition 3.1.} Consider the Poisson manifolds $\mathbb{N}_k=\{x^a, ~ \{{}, {}\}_{\mathbb{N}}\}$ and $\mathbb{M}_n=\{z^i, ~ \{ {},  {}\}_{\mathbb{M}}\}$. The mapping (\ref{sim.2.6}) is called Poisson mapping if the induced mapping (\ref{sim.2.12}) preserves the Poisson brackets
\begin{eqnarray}\label{sim15.8}
\{ \phi^*(A), \phi^*(B)\}_{\mathbb{N}}=\phi^*\left(\{ A,   B\}_{\mathbb{M}}\right). 
\end{eqnarray} 
This allows us to compare Poisson brackets of $\mathbb{M}$ and $\mathbb{N}$: given two functions $A$ and $B$ of $\mathbb{M}$ and their images $\bar A$ and $\bar B$, we can compare the bracket $\{ \bar A,  \bar B\}_{\mathbb{N}}$ with the image of scalar function $\{A,  B\}_{\mathbb{M}}$, that is with $\phi^*\left(\{A,  B\}_{\mathbb{M}}\right)$. If they coincide, we have the mapping (\ref{sim.2.8}) that respects the Poisson structures of the manifolds. The mapping (\ref{sim.15.1}) is an example of Poisson mapping of the canonical Poisson manifold on the Lie-Poisson manifold. 

{\bf Poisson submanifold of the Poisson manifold.} Let the Poisson manifold $\mathbb{N}_k$ be a submanifold of Poisson manufold $\mathbb{M}_n$, determined by functionally independent set of scalar functions $\Phi^\beta(z^i)$ of $\mathbb{M}_n$ 
\begin{eqnarray}\label{sim15.11}
\mathbb{N}_k=\{ z^i\in\mathbb{M}; \quad \Phi^\beta(z^i)=0 \}. 
\end{eqnarray}
Solving $\Phi^\beta(z^i)=0$, we obtain the parametric equations $z^\beta=f^\beta(z^a)$, and take $z^a$ as local  coordinates of $\mathbb{N}_k$. Any scalar function $A(z^i)$ on $\mathbb{M}_n$ is defined, in particular, at the points of $\mathbb{N}_k$, and hence we can consider the restriction of $A(z^i)$ on $\mathbb{N}_k$
\begin{eqnarray}\label{sim15.12}
\eta^*: \mathbb{F}_{\mathbb{M}} \rightarrow \mathbb{F}_{\mathbb{N}}, \qquad A(z^\beta, z^a) \rightarrow \bar A(z^a)=A(f^\beta(z^a), z^a),
\end{eqnarray}
The Poisson manifold $\mathbb{N}_k$ is called Poisson submanifold of $\mathbb{M}_n$, if the mapping $\eta^*$ turn the bracket of $\mathbb{M}$ into the 
bracket of $\mathbb{N}$: 
\begin{eqnarray}\label{sim15.10}
\eta^*\left(\{A, B\}_M\right)=\{\bar A, \bar B \}_N. 
\end{eqnarray} 
or 
\begin{eqnarray}\label{sim15.13}
\left. \{ A(z^\beta, z^a), B(z^\beta, z^a) \}_{\mathbb{M}}\right|_{z^\beta \rightarrow f^\beta(z^a)}=\{ A(f^\beta(z^a), z^a), B(f^\beta(z^a), z^a) \}_{\mathbb{N}}. 
\end{eqnarray} 
Various examples of Poisson mappings and Poisson submanifolds will appear in the analysis of dynamical systems in Sect. \ref{FirstInt2} .  
Notice that $\eta^*$ determined by (\ref{sim15.12}) is the mapping induced by  the embedding mapping (\ref{sim.2.7}).

\section{Degenerate Poisson manifold.}\label{DegPoi}

The affirmations discussed above are equally valid for nondegenerate and degenerate manifolds. Now we consider some characteristic properties of a Poisson manifold with degenerate Poisson bracket. Nondegenerate Poisson manifolds will be discussed in Sect. \ref{SymSym}.   

\subsection{Casimir functions.}\label{IV A} 

Poisson manifold with a degenerate bracket has the following property:  in the space $\mathbb{F}_{\mathbb{M}}$ there is a set of functionally independent functions, that have null brackets (commute) with all functions of $\mathbb{F}_{\mathbb{M}}$. They are called 
the Casimir functions. This allows to construct a remarcable foliation of the manifold with the leaves determined by the Casimir functions. \par

\noindent {\bf Affirmation 4.1.} Let $K_\beta$, $\beta=1, 2, \ldots , p$, are $p$ functionally independent Casimir functions of a Poisson manifold $\mathbb{M}_n$. Then $\omega$ is degenerated, and $rank ~ \omega\le n-p$. \par

\noindent {\bf Proof.} $K_\beta$ commutes with any function, in particular, we can write $\{z^i, K_\beta\}=0$, or $\omega^{ij}\partial_j K_\beta=0$. The latter equation means that $\omega$ admits at least  $p$ independent null-vectors $\vec V_{(\beta)}$, so $rank ~ \omega\le n-p$. $\blacksquare$ \par

\noindent {\bf Affirmation 4.2.}  Consider Poisson manifold with $rank ~ \omega=n-p$. Then there are exactly $p$ functionally independent Casimir functions:
\begin{eqnarray}\label{dar.3}
\{z^i, K_\beta \}=0, \quad \mbox{or} \quad  \vec X_{K_\beta}=0,   \quad i=1, 2, \ldots , n, \quad \beta=1, 2, \ldots , p. 
\end{eqnarray} 

\noindent {\bf Proof.} First, let us consider the particular case of $2n+1$\,-dimensional Poisson manifold with $rank ~ \omega=2n$. According the Darboux theorem, there are canonical coordinates $z^{i'}$ such that one of them coomutes with all others, say $z^{1'}$ commute with all coordinates, $\{z^{i'}, z^{1'}\}\equiv \omega^{i'1'}= 0$. Let us define a scalar function as follows: at the point $z\in\mathbb{M}_{2n+1}$, its value coincides with the value of first coordinate of this point in the canonical system: $K(z)=z^{1'}$. In the canonical coordinates this function is represented by $K'(z^{1'}, z^{2'}, \ldots , z^{(2n+1)'})=z^{1'}$. Then, according to Eqs. (\ref{sim.2}) and (\ref{sim.1}),  in the original coordinates it is represented 
by $K(z^i)=z^{1'}(z^1, z^2, \ldots , z^{2n+1})$. Let us confirm that $K(z)$ is the Casimir function. Using the transformation laws (\ref{sim.2}), (\ref{sim.3}) and (\ref{sim.4}), we obtain
\begin{eqnarray}\label{dar.4}
\{z^i, K(z)\}=\omega^{ij}(z)\partial_jK(z)=\left.\frac{\partial z^i}{\partial z^{i'}}\right|_{z'(z)}\omega^{i'j'}\left.\frac{\partial z^j}{\partial z^{j'}}\right|_{z'(z)}\frac{\partial z^{k'}}{\partial z^j}\left.\frac{\partial K'(z')}{\partial z^{k'}}\right|_{z'(z)}= \cr
\left.\left[\frac{\partial z^i}{\partial z^{i'}}\omega^{i'j'}\frac{\partial z^{1'}}{\partial z^{j'}}\right]\right|_{z'(z)}=\left.\frac{\partial z^i}{\partial z^{i'}}\right|_{z'(z)}\omega^{i'1'}=0.
\end{eqnarray} 
Let us return to the general case with $rank ~\omega=n-p$ Casimir functions. According to Eq. (\ref{dar.2}), in the Darboux coordinates the functions of $z^{\beta'}$ are Casimir functions. As the complete set of functionally independent Casimir functions,  we can take the coordinates $z^{\beta'}$ themselves. More than $p$ functionally independent Casimir functions would be in contradiction with Affirmation 4.1.   $\blacksquare$

Consider the foliation with the leaves determined by the Casimir 
functions, $\mathbb{N}_{n-p}^{\vec c}=\{ z\in\mathbb{M}_n, \quad  K_\beta(z^i)=c_\beta \}$. Then equations (\ref{dar.3}) have the following remarkable interpretation: for any function $A\in{\mathbb F}_{{\mathbb M}_n}$, the Hamiltonian vector field $X^i_A=\omega^{ij}\partial_j A$  is tangent to the hypersurfaces $\mathbb{N}_{n-p}^{\vec c}$, that is its integral lines lie in $\mathbb{N}_{n-p}^{\vec c}$.  Indeed, let $z^i(\tau)$ be an integral line of $X_A^i$. We get: $\frac{d}{d\tau}K_\beta(z^i(\tau))=\left.X_A^i\partial_i K_\beta(z^i)\right|_{z(\tau)}=\left. \{ K_\beta, A\}\right|_{z(\tau)}=0$. 
Then $K_\beta(z^i(\tau))=c_\beta=\mbox{const}$, that is $z^i(\tau)$ lies on one of the surfaces, so $\vec X_A\in\mathbb{T}_{\mathbb{N}}$.  

{\sc Exercise 4.1.} Observe: $K=z^i z^i$ is the Casimir function of (\ref{sim.16.1}).

\subsection{Induced bracket on the Casimir submanifold.}\label{IV B} 

Consider degenerate Poisson manifold ${\mathbb M}_n=\{ z^i; ~ \{A, B\}=\partial_i A  \omega^{ij}  \partial_j B, ~ rank ~ \omega=n-p ~ \}$, and let $K_\beta(z^i)$ is a subset of Casimir functions (we can take either all functionally independent Casimirs, $\beta=1, 2, \ldots , p$, or some part of them). 
Consider the submanifold determined by $K_\beta$
\begin{eqnarray}\label{dar.5}
\mathbb{N}=\{ z^i\in\mathbb{M}_n, \quad  K_\beta(z^i)=0 \}.
\end{eqnarray} 
For shorteness, we call $\mathbb{N}$ the Casimir submanifold. We will show that Poisson bracket  on $\mathbb{M}_n$ can be used to construct a natural Poisson bracket $\{ {}, {} \}_{\mathbb N}$ on $\mathbb{N}$.

{\bf Induced bracket in special coordinates.} Since the functions $K_\beta(z^i)$ are functionally independent, we can take the coordinate system where they turn into a part of coordinates, say $\tilde z^i=(\tilde z^\beta=K_\beta, \tilde z^a)$. On the surface $\mathbb{N}$ we have $\tilde z^\beta=0$, so $\tilde z^a$ are the coordinates of $\mathbb{N}$. The Poisson tensor $\tilde\omega^{ij}=\{\tilde z^i, \tilde z^j\}$ of $\mathbb{M}_n$ in these coordinates has the following special form: $\{\tilde z^a, \tilde z^b\}=\omega^{ab}(\tilde z^\beta, \tilde z^c)$, $\{\tilde z^\beta, \tilde z^i\}=\{ K_\beta, \tilde z^i\}=0$, for any $i$. Since $\tilde\omega^{ij}$ obeys (\ref{sim.D5}) and (\ref{sim.E5}) for any value of the coordinates $\tilde z^\beta$, we get $\omega^{ab}(0, z^c)=-\tilde\omega^{ba}(0, z^c)$ and 
\begin{eqnarray}\label{ar.1}
\tilde\omega^{ap}(\tilde z^\beta, \tilde z^c)\partial_p\tilde\omega^{bc}(\tilde z^\beta, \tilde z^c)+cycle (a, b, c)=0, \quad \mbox{or} \quad 
\{\tilde z^a, \{\tilde z^b, \tilde z^c\}\}+cycle (a, b, c)=0, 
\end{eqnarray} 
where the index $p$ runs over both $\beta$ and $a$ subsets. But, observing that 
\begin{eqnarray}\label{ar.2}
\{\tilde z^a, \{\tilde z^b, \tilde z^c\}\}=\{\tilde z^a, \tilde\omega^{bc}(\tilde z^\beta, \tilde z^c)\}=\{\tilde z^a, \tilde z^\beta\}\partial_\beta\tilde\omega^{bc}(\tilde z^\beta, \tilde z^c)+\{\tilde z^a, \tilde z^d\}\partial_d\tilde\omega^{bc}(\tilde z^\beta, \tilde z^c)=
\omega^{ad}(\tilde z^\beta, \tilde z^c)\partial_d\omega^{bc}(\tilde z^\beta, \tilde z^c), 
\end{eqnarray} 
we can write the first equality in (\ref{ar.1}) as follows:
\begin{eqnarray}\label{ar.3}
\omega^{ad}(\tilde z^\beta, \tilde z^c)\partial_d\omega^{bc}(\tilde z^\beta, \tilde z^c)+cycle (a, b, c)=0.  
\end{eqnarray} 
As it is true for any value of the coordinates $\tilde z^\beta$, we can take $\tilde z^\beta=0$, then
\begin{eqnarray}\label{ar.4}
\tilde\omega^{ad}(0, \tilde z^c)\partial_d\tilde\omega^{bc}(0, \tilde z^c)+cycle (a, b, c)=0.  
\end{eqnarray} 
Let us define the matrix with elements 
\begin{eqnarray}\label{ar.5}
\bar\omega^{ab}(\tilde z^c)\equiv\tilde\omega^{ad}(0, \tilde z^c),  
\end{eqnarray}
in the coordinates $\tilde z^c$. In any other coordinate system on $\mathbb{N}$, say $\tilde z^{a'}$, we define the elements $\bar\omega^{a'b'}$ according the rule (\ref{sim.4}). Then $\bar\omega$ is a tensor of $\mathbb{N}$. According to our computations, it obeys to Eqs. (\ref{sim.D5}) and (\ref{sim.E5}), and thus determines a Poisson bracket $\{A, B\}_{\mathbb N}=\partial_i A(\tilde z^c)  \bar\omega^{ab}(\tilde z^c) \partial_j B(\tilde z^c)$ on ${\mathbb N}$.

{\bf Induced bracket in the original coordinates.} Let us solve the same problem in the original coordinates, divided in two subsets, $z^i=(z^\alpha,  z^a)$,   such that $\det\frac{\partial K_\beta}{\partial z^\alpha}\ne 0$. Notice that in this case the Eq. (\ref{dar.3}) reads 
\begin{eqnarray}\label{dar.4.01}
\omega^{i\alpha}\partial_\alpha K_\beta+\omega^{i b}\partial_b K_\beta=0.  
\end{eqnarray} 
Denoting $\partial_\alpha K_\beta=K_{\alpha\beta}$, this allows us to restore the whole $\omega^{ij}(z^k)$ from the known block $\omega^{ab}(z^k)$ as follows:
\begin{eqnarray}\label{dar.4.1}
\omega^{a\alpha}=-\omega^{ab}\partial_b K_\gamma(K^{-1})^{\gamma\alpha}, \qquad \omega^{\alpha\beta}=-\omega^{\alpha b}\partial_b K_\gamma(K^{-1})^{\gamma\beta}.
\end{eqnarray} 
Geometric interpretation of these relations will be duscussed in Sect. \ref{PresCas}.

It is instructive to obtain the induced bracket in the original coordinates in a manner, independent of calculations made in the previous subsection.

\noindent {\bf Affirmation 4.3.} Let $K_\beta(z^\alpha, z^b)$  are Casimir functions, and $z^\alpha=f^\alpha(z^b)$ is a solution to the equations $K_\beta(z^\alpha, z^b)=0$. Then \par
\noindent a) $z^\alpha-f^\alpha(z^b)$ are Casimir functions; \par 
\noindent b) The Poisson tensor of $\mathbb{M}_n$ satisfies the identity
\begin{eqnarray}\label{dar.4.2}
\omega^{i\alpha}=\omega^{ib}\partial_b f^\alpha. 
\end{eqnarray} 
{\bf Proof.}  a) Contracting the expression $\{z^i, z^\alpha-f^\alpha\}=\omega^{i\alpha}-\omega^{i a}\partial_a f^\alpha$ with $\partial_\alpha K_\beta$ and using (\ref{dar.4.01}),  we obtain $\omega^{i\alpha}\partial_\alpha K_\beta-\omega^{i a}\partial_a f^\alpha\partial_\alpha K_\beta=-\omega^{i a}(\partial_a K_\beta+\partial_a f^\alpha\partial_\alpha K_\beta)=-\omega^{i a}\partial_a K_\beta(f^\alpha, z^a)=0$ since $K_\beta(f^\alpha, z^a)\equiv 0$. As $\partial_\alpha K_\beta$ is an invertible matrix, the equality $\{z^i, z^\alpha-f^\alpha\}\partial_\alpha K_\beta=0$ implies $\{z^i, z^\alpha-f^\alpha\}=0$.  \par

b) Let $K_\alpha$ are Casimir functions. According to Item  a),  $z^\alpha-f^\alpha(z^b)$ also are the Casimir functions, then $\{z^i, z^\alpha-f^\alpha\}=0$, or $\omega^{i\alpha}=\omega^{i b}\partial_b f^\alpha$.  $\blacksquare$

\noindent {\bf Affirmation 4.4.}  For any function $B(z^i)$ and Casimir functions $z^\beta-f^\beta(z^b)$, there is the identity
\begin{eqnarray}\label{dar.5.1}
\left.\omega^{ip}\partial_p B\right|_{z^\beta=f^\beta(z^c)}=\omega^{ia}(z^c, f^\beta(z^c))\partial_a B(z^c, f^\beta(z^c)),
\end{eqnarray}
where $p=(1, 2, \ldots ,  n)$, while $a=(1, 2, \ldots ,  n-p)$. Note the geometric interpretation of this equality: if two functions $B$ and $B'$ 
of $\mathbb{F}_{\mathbb{M}}$ coincide on $\mathbb{N}$,  their Hamiltonian vector fields also coincide on $\mathbb{N}$: $B|_{\mathbb{N}}=B'|_{\mathbb{N}}$ implies $\vec X_B|_{\mathbb{N}}=\vec X_{B'}|_{\mathbb{N}}$. \par 

\noindent {\bf Proof.} Let us write
\begin{eqnarray}\label{dar.5.2}
\left.\omega^{ip}\partial_p B\right|_{z^\beta=f^\beta(z^c)}=\left.\omega^{ia}\partial_a B\right|_{z^\beta=f^\beta(z^c)}+\left.\omega^{i\beta}\partial_\beta B\right|_{z^\beta=f^\beta(z^c)}.
\end{eqnarray}
Using the identity (\ref{dar.4.2}) we have $\omega^{i\beta}\partial_\beta B=\omega^{id}\partial_d f^\beta(z^c)\partial_\beta B(z^c, z^\beta)\equiv\omega^{id}\left[\partial_d B(z^c, f^\beta(z^c))-\left.\partial_d B(z^c, z^\beta)\right|_{z^\beta=f^\beta(z^c)}\right]$. Using this expression for the term $\omega^{i\beta}\partial_\beta B$ in (\ref{dar.5.2}), we arrive at the desired identity (\ref{dar.5.1}). $\blacksquare$

We are ready to construct the induced Poisson structure. We take $z^a$ as local coordinates 
of $\mathbb{N}$, and using $\omega^{ab}$\,-block of $\omega^{ij}$, introduce the antisymmetric matrix 
\begin{eqnarray}\label{dar.6}
\bar\omega^{ab}(z^c)=\omega^{ab}(f^\beta(z^c), z^c).  
\end{eqnarray} 
Let us confirm that $\bar\omega$ obeys the condition (\ref{sim.E5}). We write the condition (\ref{sim.E5}), satisfied for $\omega^{ij}$, taking  the indices $i, j, k$ equals to $a, b, c$,  and substitute $z^\beta=f^\beta(z^c)$. This gives us the identity  
\begin{eqnarray}\label{dar.7}
\left.\omega^{ap}\partial_p\omega^{bc}\right|_{z^\beta=f^\beta(z^c)}+ cycle(a, b, c)=0. 
\end{eqnarray}
Using the identity (\ref{dar.5.1}), we immediately obtain
\begin{eqnarray}\label{dar.8}
\omega^{ad}(f^\beta(z^c), z^c)\partial_d\omega^{bc}(f^\beta(z^c), z^c)+cycle(a, b, c)=0,
\end{eqnarray}
which is just the Jacobi identity for the tensor $\bar\omega^{ab}$. Thus the bracket 
\begin{eqnarray}\label{dar.9}
\{ A(z^a), B(z^a)\}=\partial_a A\bar\omega^{ab}\partial_b B, 
\end{eqnarray} 
defined on $\mathbb{N}$, obeys the Jacobi identity.  In any other coordinate system on $\mathbb{N}$, say $z^{a'}$, we define the components $\bar\omega^{a'b'}$ according the rule (\ref{sim.4}): 
\begin{eqnarray}\label{dar.9.1}
\bar\omega^{a'b'}=\left.\partial_az^{a'}\partial_bz^{b'}\bar\omega^{ab}\right|_{z^a(z^{a'})}. 
\end{eqnarray} 
Then $\bar\omega$ is a tensor of $\mathbb{N}$, while the expression (\ref{dar.9}) is a scalar function, as it should be for the Poisson bracket. 

Let us confirm, that the obtained bracket does not depend on the coordinates of ${\mathbb M}_n$ chosen for its construction. Let $z^i=\phi^i(z^{j'})$ be transition functions between two coordinate systems. For a point of $\mathbb{N}$, this implies the following relation between its local coordinates $z^a$ and $z^{a'}$:
\begin{eqnarray}\label{dar.9.1.0}
z^a=\phi^a(f^{\alpha'}(z^{a'}), z^{a'}).
\end{eqnarray} 
Using these finctions in the expression (\ref{dar.9.1}), we obtain the components $\bar\omega^{a'b'}(z^{a'})$ of the tensor $\bar\omega^{ab}(z^a)$ in the coordinates $z^{a'}$. On the other hand, using the Poisson tensor $\omega^{i'j'}$ in coordinates $z^{i'}$, we could construct the matrix $\hat\omega^{a'b'}(f^{\alpha'}(z^{a'}), z^{a'})$ according the rule (\ref{dar.6}). The task is to show that $\hat\omega^{a'b'}$ coincides with $\bar\omega^{a'b'}$.

As the functions $\hat\omega^{a'b'}(f^{\beta'}(z^{a'}), z^{a'})$ are components of the tensor $\omega^{i'j'}$ of $\mathbb{M}$, we use the transformation law (\ref{sim.4}), and write
\begin{eqnarray}\label{dar.9.2}
\hat\omega^{a'b'}(f^{\beta'}(z^{a'}), z^{a'})=
\left.\hat\omega^{a'b'}(z^{\beta'}, z^{a'} )\right|_{\mathbb{N}}=
\left.\left.\partial_kz^{a'}\omega^{kp}\partial_pz^{b'}\right|_{z^i(z^{i'})}\right|_{\mathbb{N}}. 
\end{eqnarray} 
In the last expression we have a quantity $D(z^i)$, and need to replace the coordinates $z^i$ by the transition functions  $z^i(z^{i'})$ taken at the point of $\mathbb{N}$. Equivalently, we can first to restrict $D$ on $\mathbb{N}$, replacing $z^\beta$ on $f^\beta(z^a)$,
and then to replace $z^a$ on its expression (\ref{dar.9.1.0}) through coordinates $z^{a'}$.  Making this, and then using the identity (\ref{dar.5.1}), we obtain
\begin{eqnarray}\label{dar.9.3}
\hat\omega^{a'b'}(f^{\beta'}(z^{a'}), z^{a'})=
\left.\left.\partial_kz^{a'}(z^i)\omega^{kp}(z^i)\partial_pz^{b'}(z^i)\right|_{z^\beta\rightarrow f^\beta(z^a)}\right|_{z^a(z^{a'})}= \cr
\left.\partial_az^{a'}(z^\beta(z^a), z^a)\omega^{ab}(z^\beta(z^a), z^a)\partial_bz^{b'}(z^\beta(z^a), z^a)\right|_{z^a(z^{a'})}. \qquad \quad 
\end{eqnarray} 
Comparing this expression with (\ref{dar.9.1}), we arrive at the desired result: $\hat\omega=\bar\omega$.   \par

{\sc Exercise 4.2.} Confirm that the Poisson manifold $\mathbb{N}$ is the Poisson submanifold of $\mathbb{M}$ in the sense of definition (\ref{sim15.10}).  

Consider the Poisson manifold $\mathbb{M}_n$ with $rank~\omega^{ij}=n-p$, and let  the submanifold (\ref{dar.5}) is determined by a complete set of $p$ functionally independent Casimir functions. Then the induced Poisson structure is non degenerate: $\det \bar\omega^{ab}\ne 0$. To see this, suppose an opposite, $\det \bar\omega^{ab}= 0$, and let $z^i$ be canonical coordinates of $\mathbb{M}$. Then  $\bar\omega$ is a numeric degenerate matrix, so it has a numeric null-vector, $\bar\omega^{ab} c_b=0$. As a consequence, the function $A(z^i)=z^ac_a$ commutes with all coordinates (\ref{dar.2}), and hence is a Casimir function of $\mathbb{M}$. It depends only on the variables $z^a$, so it is functionally independent of the Casimir functions $z^\beta-f^\beta (z^a)=0$. This is in contradiction with the condition $rank~\omega^{ij}=n-p$, so $\det \bar\omega^{ab}\ne 0$.

Let us resume the obtained results. Let $\mathbb{M}$ be Poisson manifold with degenerate Poisson bracket $\omega$. Then on the submanifold $\mathbb{N}\in \mathbb{M}$, determined by any set of functionally independent Casimir functions, there exists the Poisson bracket $\bar\omega$, such that the Poisson manifold $\mathbb{N}$ turns into the Poisson submanifold of $\mathbb{M}$. In the original coordinates, divided on two groups  according to the structure of Casimir functions (\ref{dar.5}), $z^i=(z^\beta, z^a)$, elements of the matrix $\bar\omega$ coincide with fundamental brackets of coordinates $z^a$ restricted on $\mathbb{N}$: 
\begin{eqnarray}\label{dar.9.4}
\bar\omega^{ab}=\left.\{z^a, z^b\}_{\mathbb M}\right|_{z^\beta\rightarrow f^\beta(z^a)}.
\end{eqnarray}

\subsection{Restriction of Hamiltonian dynamics to the Casimir submanifold.}\label{IV C} 

The degeneracy of a Poisson structure implies that integral lines of {\it any} Hamiltonian system on this manifold have special properties: any solution started in a Casimir submanifold entirely lies in it. So the dynamics can be consistently restricted on the submanifold, and the resulting equations are still Hamiltonian.  To discuss these properties, we will need the notion of an invariant submanifold. \par 

\noindent {\bf Definition 4.1.} The submanifold $\mathbb{N}\in\mathbb{M}_n$ is called an invariant submanifold of the Hamiltonian $H$, if any trajectory of (\ref{sim.12})  that starts in $\mathbb{N}$, entirely lies in $\mathbb{N}$
\begin{eqnarray}\label{sim.16.3}
z^i(0)\in\mathbb{N}, \quad \rightarrow \quad z^i(\tau)\in\mathbb{N} \quad \mbox{for any} \quad \tau. 
\end{eqnarray}
The observation made in Sect.\ref{IV A} now can be rephrased as follows. 

\noindent {\bf Affirmation 4.5.}  A Casimir submanifold of $\mathbb{M}_n=\{z^i, ~ \{{}, {}\}\}$ is invariant submanifold of any Hamiltonian $H\in \mathbb{F}_{\mathbb{M}}$. \par

\noindent {\bf Affirmation 4.6.}  Solutions to the Hamiltonian equations: $\dot z^i=\omega^{ij}\partial_j H$, that belong to the Casimir submanifold (\ref{dar.5}), obey the Hamiltonian equations 
\begin{eqnarray}\label{sim.9.3.2} 
\dot z^a=\bar\omega ^{ab}\partial_b H(z^c, f^\alpha(z^c)), 
\end{eqnarray}
where $z^a$ are local coordinates,  and the Poisson tensor $\bar\omega^{ab}$ of $\mathbb{N}$ is the restriction of $\omega^{ij}$ on $\mathbb{N}$
\begin{eqnarray}\label{sim.9.3.3}  
\bar\omega^{ab}=\omega^{ab}(z^c, f^\alpha(z^c)). 
\end{eqnarray}

\noindent {\bf Proof.}  According to the Affirmation 4.5, we can add the algebraic equations $z^\beta=f^\beta(z^a)$ to the system (\ref{sim.12}), thus obtaining a consistent equations with solutions living on $\mathbb{N}$. In the equations for $\dot z^a$ we substitute $z^\beta=f^\beta(z^a)$, and using the identity (\ref{dar.5.1}), obtain the closed system (\ref{sim.9.3.2}) and (\ref{sim.9.3.3}) for determining $z^a$. Then the equations for $\dot z^\beta$ can be omitted from the system.  Jacobi identity for $\bar\omega$ has been confirmed above. $\blacksquare$ \par

\section{Integrals of motion of a Hamiltonian system.}\label{FirstInt}
\subsection{Basic notions.}

Let $z(\tau)$ is a solution to Hamiltonian equations (\ref{sim.12}).  For any function $Q(z)$ we have: $\dot Q(z(\tau))=\left.\{Q(z), H(z)\}\right|_{z(\tau)}$. In other words, functions $Q(z(\tau))$ follow the Hamiltonian dynamics together with $z(\tau)$. The function $Q(z)$ (with non-vanishing gradient) is called the integral of motion, if it preserves its value along the trajectories of 
(\ref{sim.12}): $Q(z(\tau))=const$, or $\dot Q(z(\tau))=0$. Note that the value of $Q(z(\tau))$ can vary from one trajectory to another. 

\noindent {\bf Affirmation 5.1.} $Q(z)$ is an integral of motion of the system (\ref{sim.12}), if and only if its bracket with $H$ vanishes 
\begin{eqnarray}\label{sim.13}
\{Q, H\}=0. 
\end{eqnarray}

Since $\{H, H \}=0$, the Hamiltonian himself is an example of the integral of motion. So, any Hamiltonian system admits at least one integral of motion. The Casimir functions obey to Eq. (\ref{sim.13}) for any $H$, so they represent integrals of motion of any Hamiltonian system on a given manifold. As a consequence, a Hamiltonian system on the manifold $\mathbb{M}_n$ with $rank~\omega=n-p$ has at least $p+1$ integrals of motion. \par

\noindent {\sc Exercise 5.1.} a) Confirm the Affirmation 5.1. (hint: take into account, that integral lines of (\ref{sim.12}) cover all the manifold). \par
\noindent b) Observe: if $Q_1$ and $Q_2$ are integrals of motion, then $c_1Q_1+c_2Q_2$, $f(Q_1)$ and $\{Q_1, Q_2 \}$ are the integrals of motion as well. The integral of motion $\{Q_1, Q_2 \}$ may be functionally independent of $Q_1$ and $Q_2$. 

The integrals of motion $Q_\alpha$ can be used to construct the surfaces of level in $\mathbb{M}_n$. Considering the Hamiltonian equations on the surfaces, it turns out to be possible to reduce the number of differential equations that we need to solve. This method, called the reduction procedure, is based on the following affirmations. 

\noindent {\bf Affirmation 5.2.} Let $Q_\alpha (z)$, $\alpha=1, 2, \ldots , p$ are functionally independent integrals of motion of $H$. Then $\mathbb{N}_{{\bf c}}=\{ z\in\mathbb{M}_n, \quad Q_\alpha(z)=c_\alpha=const \}$ is $n-p$\,-dimensional invariant submanifold of $H$. \par 

Indeed, given solution with $z(0)\in\mathbb{N}_{{\bf c}}$, that is $Q_\alpha (z(0))=c_\alpha$, we have $Q_\alpha(z(\tau))=Q_\alpha (z(0))=c_\alpha$ for any $\tau$,  so  the trajectory $z(\tau)$ entirely lies in $\mathbb{N}_{{\bf c}}$.  The manifolds $\mathbb{N}_{{\bf c}}$ and $\mathbb{N}_{{\bf d}}$ with ${{\bf c}}\ne {{\bf d}}$ do not intercept. So the Poisson manifold $\mathbb{M}_n$ is covered by $p$\,-parametric foliation of the invariant submanifolds $\mathbb{N}_{{\bf c}}$. 

Since the Casimir function is an integral of motion of any Hamiltonian, the Affirmation 5.2 implies, once again, the geometric interpretation \label{int} of Eq. (\ref{dar.3}): integral lines of all Hamiltonian vector fields of $\mathbb{M}_n$ lie on the surfaces of Casimir functions.   

\noindent {\bf Affirmation 5.3.} Let the Hamiltonian system 
\begin{eqnarray}\label{sim.16.4}
\dot z^i=\{z^i, H\},  
\end{eqnarray}
admits $p$ functionally independent integrals of motion $Q_\alpha(z)=c_\alpha$. We present them in the form $z^\alpha=f^\alpha(z^b, c_\alpha)$. Then the system of $n$ differential equations (\ref{sim.16.4}) is equivalent to the system 
\begin{eqnarray}\label{sim.16.5}
\dot z^b=\{z^b, H\}\equiv h^b(z^c, z^\alpha),  \qquad z^\alpha=f^\alpha(z^b, c_\alpha),
\end{eqnarray}
composed of $n-p$ differential and $p$ algebraic equations. \par 

\noindent {\bf Proof.} Adding the consequences $z^\alpha=f^\alpha(z^b, c_\alpha)$ to the equations (\ref{sim.16.4}), we write the resulting equivalent system as follows 
\begin{eqnarray}\label{sim.16.6}
\dot z^\alpha=\{z^\alpha, H\}, \qquad  
\dot z^b=\{z^b, H\}, \qquad  z^\alpha=f^\alpha(z^b, c_\alpha). 
\end{eqnarray}
To prove the equivalence of (\ref{sim.16.5}) and (\ref{sim.16.6}), we need to show that the equation $\dot z^\alpha=\{z^\alpha, H\}$ is a  consequence of the system (\ref{sim.16.5}). Let $z^\alpha(\tau)$, $z^b(\tau)$ be a solution to (\ref{sim.16.5}). Computing derivative of the identity $z^\alpha(\tau)\equiv f^\alpha(z^b(\tau), c_\alpha)$ we have $\dot z^\alpha(\tau)=\left.\partial_b f^\alpha(z^b, c_\alpha)\right|_{z\rightarrow z(\tau)}\dot z^b=
\left.\partial_b f^\alpha(z^b, c_\alpha)\{z^b, H\}\right|_{z\rightarrow z(\tau)}=\left.\{f^\alpha, H\}\right|_{z\rightarrow z(\tau)}=\left.\{z^\alpha, H\}\right|_{z\rightarrow z(\tau)}$. On the last step we used (\ref {sim.13}). Hence, the equation $\dot z^\alpha=\{z^\alpha, H\}$ is satisfied by any solution 
to the system (\ref{sim.16.5}). $\blacksquare$ \par

{\sc Example 5.1.} Using the reduction procedure, any two-dimensional Hamiltonian system can be completely integrated, that is solving the differential equations is reduced to evaluation of an integral. Indeed, consider the system $\dot x=\{x, H(x, y)\}\equiv h(x, y)$, $\dot y=\{y, H(x, y)\}$. We assume that $grad ~ H\ne 0$ (otherwise $H=const$ and the system is immediately integrated). Let $y=f(x, c)$ is a solution to the equation $H(x, y)=c$. As $H$ is an integral of motion, we use the Affirmation 5.3 to present the original system 
in the equivalent form: $\dot x= h(x, y)$, $y=f(x, c)$. Replacing  $y$ on $f(x, c)$ in the differential equation, the latter can be immediately integrated. The general solution $x(\tau, c, d)$, $y(\tau, c, d)$ in an implicit form is as follows:
\begin{eqnarray}\label{sim.16.7}
\int \frac{dx}{h(x, f(x, c))}=\tau+d, \qquad y=f(x). 
\end{eqnarray}
There is a kind of multi-dimensional generalization of this example, see Affirmation B1 in Appendix \ref{App B}. 

\subsection{Hamiltonian reduction to an invariant submanifold.}\label{FirstInt2} 

As we saw above, when a dynamical system admits an invariant submanifold, its dynamics can be consistently restricted to the submanifold. Then it is natural to ask, whether the resulting equations form a Hamiltonian system? For instance, according to Affirmation 5.2, we can add the algebraic equations\footnote{Without loss of generality, we have taken $c_\alpha=0$.} $Q_\alpha(z)=0$ to the Hamiltonian system (\ref{sim.16.4}), thus obtaining a consistent equations with solutions living on the invariant submanifold $\mathbb{N}=\{ z\in\mathbb{M}_n , Q_\alpha(z)=0 \}$. Using Affirmation 5.3, we exclude $z^\alpha$ and obtain differential equations on the manifold $\mathbb{N}$ with the local coordinates $z^b$  
\begin{eqnarray}\label{sim.16.8} 
\dot z^b=h^b(z^c)\equiv\left.\{z^b, H\}\right|_{z^\alpha\rightarrow f^\alpha(z^c)}. 
\end{eqnarray}
They already don't know anything about the ambient space $\mathbb{M}_n$. Hence we ask, if the resulting equations represent a Hamiltonian system on $\mathbb{N}$?  That is we look for the Hamiltonian equations 
\begin{eqnarray}\label{sim.16.8.0} 
\dot z^b=\omega^{ba}(z^c)\partial_a \hat H(z^c), 
\end{eqnarray}
that could be equivalent to (\ref{sim.16.8}). 

Let us list some known cases of the Hamiltonian reduction. 

{\bf 1.} Reduction of non singular theory (\ref{sim.0.03}) to the surface of constant Hamiltonian gives a Hamiltonian system with time-dependent Hamiltonian. The method is known as Maupertuis principle, see \cite{deriglazov2010classical} for details. 

{\bf 2.} Hamiltonian reduction to the surface of Casimir functions, see Affirmation 4.6. The particular example is a Hamiltonian system with Dirac bracket, see Eq. (\ref{sim.31}) below. 

{\bf 3.} Hamiltonian reduction of non singular theory to the surface of first integrals $\Phi^\alpha$ with the property $\det\{\Phi^\alpha, \Phi^\beta\}\ne 0$, see Eq. (\ref{sim.34.5}) below. 

{\bf 4.} Singular non degenerate theory (\ref{sim.0.5})-(\ref{sim.0.5.2}) is equivalent to the theory of Item 2, see Affirmations 7.5 and 7.3 below. Hence it admits the Hamiltonian reduction to the surface of constraints. 

 {\bf 5.} According to Gitman-Tyutin theorem, singular degenerate theory admits Hamiltonian reduction to the surface of all constraints, see \cite{GT} for details.

\section{Symplectic manifold and Dirac bracket.}\label{SymSym}
\subsection{Basic notions.}
As we saw in Sect. \ref{SymPoi}, Poisson manifold can be defined choosing a contravariant tensor with the properties (\ref{sim.D5}) and (\ref{sim.E5}). 
Here we discuss another way, which works for the construction of  a nondegenerate Poisson structures on even-dimensional manifolds. 
Let on the even-dimensional manifold $\mathbb{M}_{2n}$ is defined the covariant tensor $\tilde\omega_{ij}(z^k)$ (called the symplectic form) with the properties 
\begin{eqnarray}
\qquad  \tilde\omega_{ij}=-\tilde\omega_{ji} \qquad \mbox{(antisymmetric)},  \label{sim.F5}\\
~ \qquad \det\tilde\omega\ne 0 \qquad \mbox{(nondegenerate)}, \label{sim.G5}\\ 
\qquad \partial_i\tilde\omega_{jk} + cycle=0 \qquad \mbox{(closed)}. \quad \label{sim.H5}
\end{eqnarray}
\noindent $\mathbb{M}_{2n}$ equipped with a symplectic form is called the symplectic manifold. 

We recall that determinant of any odd-dimensional matrix vanishes, so (\ref{sim.G5}) implies that we work on the even-dimensional manifold. Some properties of a symplectic form are in order. \par 

\noindent {\bf  Affirmation 6.1} The inverse matrix $\omega^{ij}$ of the matrix $\tilde\omega_{ij}$ obeys the properties (\ref{sim.D5}) and (\ref{sim.E5}). So it determines the Poisson structure (\ref{sim.5}) on  $\mathbb{M}_{2n}$. In other words, any symplectic manifold locally is a Poisson manifold. \par 

\noindent {\sc Exercise 6.1.} Prove that (\ref{sim.H5}) implies (\ref{sim.E5}). 

Conversely, take a Poisson manifold with the nondegenerated bracket, $\det\omega\ne 0$, and let $\tilde\omega$ is its inverse. Contracting the condition (\ref{sim.E5}) with $\tilde\omega_{ni} \tilde\omega_{mj} \tilde\omega_{pk}$ we immediately obtain (\ref{sim.H5}). We obtained\par 

\noindent {\bf Affirmation 6.2.} Poisson manifold with non degenerate bracket is a symplectic manifold. \par

\noindent {\bf  Darboux theorem.} In a vicinity of any point, there are coordinates $y^{k}$ where $\omega^{ij}(z^k)$ acquires the form 
\begin{eqnarray}\label{sim.17.0}
\omega'^{ij}(y^{k})=\left(
\begin{array}{cc}
0 & 1 \\
-1 & 0
\end{array}
\right), \quad \mbox{then} \quad   
\tilde\omega'_{ij}(y^k)=\left(
\begin{array}{cc}
0 & -1 \\
1 & 0
\end{array}
\right). 
\end{eqnarray} \par 

Proof is given in Appendix B. \par 

\noindent {\bf Poincare lemma.}  In a vicinity of any point,  the symplectic form $\tilde\omega$ can be presented through some covariant vector field $a_j$ as follows:  
\begin{eqnarray}\label{sim.17}
\tilde\omega_{ij}=\partial_i a_j-\partial_j a_i. 
\end{eqnarray}
(In the language of differential forms, this is formulated as follows: closed form is locally an exact form). Conversely, the tensor $\tilde\omega$, constructed from given $a_j$ according to Eq. (\ref{sim.17}), obeys the condition (\ref{sim.H5}). \par 

\noindent  {\bf Proof.}  According to Darboux theorem, there are coordinates $y^{i}=(x^1, \ldots , x^n, p^1, \ldots , p^n)$ where $\tilde\omega_{ij}$ acquires the canonical form (\ref{sim.17.0}). Let us identically rewrite it as follows: $\tilde\omega'_{ij}=\partial_{i} a'_{j}-\partial_{j} a'_{i}$, where $a'_{i}(y^k)=\frac 12(-p^1, \ldots , -p^n, x^1, \ldots , x^n)$. Returning to the original coordinates, we write $a_i(z^k)=\frac{\partial y^{j}}{\partial z^i}a'_{j}(y(z))$, where  $a'_{j}(y(z))=\frac 12(-p^1(z^k), \ldots , -p^n(z^k), x^1(z^k), \ldots , x^n(z^k))$.  This contravariant vector field satisfy the desired property: $\partial_i a_j-\partial_j a_i=\frac{\partial y^{n}}{\partial z^i}\frac{\partial y^{m}}{\partial z^i} \tilde\omega'_{nm}=\tilde\omega_{ij}$. 

The field $a_i$ can equally be obtained by direct integrations
\begin{eqnarray}\label{sim.17.1}
a_i=-\frac{1}{n-1}\sum_{j=1}^{n}\int\tilde\omega_{ij}(z^k)dz^j.  \qquad \qquad  \qquad  \qquad \qquad  \qquad \blacksquare
\end{eqnarray}

Due to the Poincare lemma, it is easy to construct examples of  closed and non-constant form $\tilde\omega$, and then the tensor $\omega$, that will automatically obey a rather complicated equation (\ref{sim.E5}). Note also that in the Darboux coordinates $y^k$, the Poisson bracket acquires the canonical form (\ref{sim.15}). 

Since any symplectic manifold is simultaneously a Poisson manifold, it has all the properties discussed in Sect. \ref{DynPoi}. In particular, we have the mapping
\begin{eqnarray}\label{sim.19}
\omega: \mathbb{F}_{\mathbb{M}}\rightarrow \mathbb{T}_{\mathbb{M}}, \quad \omega: A\rightarrow X_A^i=[\omega(A)]^i=\omega^{ij}\partial_j A,
\end{eqnarray}
and the basic relation between the Lie and Poisson brackets
\begin{eqnarray}\label{sim.20}
\omega(\{ A, B \})=-[\omega(A), \omega(B)], \quad \mbox{or, equivalently} \quad X_{\{A, B\}} =-[\vec X_A, \vec X_B ].  
\end{eqnarray}

Symplectic form can be used to determine the mapping $\tilde\omega:  \mathbb{T}_{\mathbb{M}}\times \mathbb{T}_{\mathbb{M}}\rightarrow \mathbb{F}_{\mathbb{M}}$ as follows
\begin{eqnarray}\label{sim.21}
\tilde\omega:  \vec X,  \vec Y \rightarrow \tilde\omega( \vec X,  \vec Y)\equiv\tilde\omega_{ij}X^iY^j, \qquad \mbox{then} \qquad \tilde\omega( \vec X,  \vec X)=0.  
\end{eqnarray}
Then the Poisson bracket can be considered as a composition\footnote{In the coordinate-free formulation of Poisson geometry, the equality  $\tilde\omega(\omega(A), \omega(B))=-\{ A, B \}$ is taken as the definition of the symplectic form $\tilde\omega$.} of the mappings (\ref{sim.21}) and (\ref{sim.19})
\begin{eqnarray}\label{sim.22}
\{ A, B \}=-\tilde\omega(\omega(A), \omega(B))\equiv -\tilde\omega(\vec X_A, \vec X_B).
\end{eqnarray} \par 

\noindent {\sc Exercise 6.2.}  a) Prove that $\{ Q, H \}=c=const$ if and only if  $[ \vec X_Q, \vec X_H]=0$. b) Confirm (\ref{sim.22}).

By analogy with Riemannian geometry, on the symplectic manifold there is a natural possibility of raising and lowering the indices of tensor quantities. It is achieved with use of symplectic tensor and its inverse. For instance, the mapping $U_i=\tilde\omega_{ij}V^j$ and its inversion $V^i=\omega^{ij}U_j$ establish an isomrphism between the spaces of covariant and contravariant vector fields.

\noindent {\bf Affirmation 6.3.} $V^i$ is a Hamiltonian vector field if and only if $U_i=\tilde\omega_{ij}V^j$ obeys the condition
\begin{eqnarray}\label{sim.22.1}
\partial_i  U_j-\partial_j U_i=0.
\end{eqnarray}

\noindent {\bf Proof.} The equation $\partial_i A=\tilde\omega_{ij}V^j\equiv U_i$ for determining of $A$ implies (\ref{sim.22.1}) as a  necessary condition. Conversely, when (\ref{sim.22.1}) is satisfied, the function 
\begin{eqnarray}\label{sim.22.2}
A=\frac{1}{n}\sum_{j=1}^{n}\int U_j(z^k)dz^j, 
\end{eqnarray}
generates the field $V^i$: $V^i=\omega^{ij}\partial_j A$.  $\blacksquare$

As an application of the developed formalism, we mention the following \par

\noindent {\bf Affirmation 6.4.} Consider the Poisson manifold $\mathbb{M}_{2n}$ with non-degenerated Poisson structure: $\det\omega\ne 0$. Let $\tilde\omega$ is the corresponding symplectic form, and $a_i$ is the contravariant vector field defined in (\ref{sim.17}). Then the Hamiltonian equations (\ref{sim.12}) follow from the variational problem
\begin{eqnarray}\label{sim.23}
S_H=\int d\tau\left[a_i(z)\dot z^i-H(z)\right]. 
\end{eqnarray}
{\sc Exercise 6.3.} Prove the affirmation\footnote{While variation of (\ref{sim.23}) formally leads to (\ref{sim.12}), the following point should be taken into account. Formulating a variational problem, we fix two points in phase space and then look for an extremal trajectory between them. The first-order system (\ref{sim.12}) has unique solution for the given initial "position": $z^i(\tau_1)=z^i_1$. This implies, that position at a future instant $\tau_2$ is uniquely determined by the initial position of the system.  So, if we look for the extremal trajectory between two arbitrary chosen points $z^i(\tau_1)=z^i_1$ and $z^i(\tau_2)=z^i_2$, the variational problem (\ref{sim.23}) generally will not have a solution.}.

\subsection{Restriction of symplectic structure to a submanifold and Dirac bracket.}\label{VI B} 

We recall that the  mapping of manifolds $\mathbb{N}=\{x^a\}\rightarrow \mathbb{M}_n=\{z^i\}$, given by  $x^a\rightarrow z^i(x^a)$, induces the mapping $\mathbb{T}^{(0, m)}_{\mathbb{M}}\rightarrow \mathbb{T}^{(0, m)}_{\mathbb{N}}$ of covariant tensor fields, see  (\ref{sim.2.13}). 
Let $\mathbb{M}_n=\{z^k, \tilde\omega_{ij}(z^k)\}$ is a symplectic manifold and $\mathbb{N}_k$ is a submanifold determined by the functions $\Phi^\alpha(z^k)=0$ (see (\ref{sim.2.1.0})),  $n$ and $k$ are even numbers. Consider the embedding $\mathbb{N}_k\rightarrow \mathbb{M}_n$, given by $x^a\rightarrow z^i=(f^\alpha(x^a), x^a)$. Then the induced mapping 
\begin{eqnarray}\label{sim.23.02}
\tilde\omega_{{\bf f} ab}(x^c)=\frac{\partial z^i}{\partial x^a}\frac{\partial z^j}{\partial x^b}\tilde\omega_{ij}(f^\alpha(x^c), x^c), 
\end{eqnarray}
is called restriction of symplectic form $\tilde\omega_{ij}(z^k)$ on $\mathbb{N}_k$. If $\tilde\omega_{{\bf f}}$ obeys the properties (\ref{sim.G5}) and (\ref{sim.H5}), $\mathbb{N}_k$ turns into a symplectic manifold. The inverse matrix then determines a Poisson bracket on $\mathbb{N}_k$. Here we discuss  the necessary and sufficient conditions under which this occurs. We will need the following matrix identity.  

\noindent {\bf Affirmation 6.5.} Consider an invertible antisymmetric matrix 
\begin{eqnarray}\label{apD.6}
A=\left(
\begin{array}{cc}
a & b \\
-b^T& c 
\end{array}\right), \quad \mbox{and its inverse} \quad 
A^{-1}=\left(
\begin{array}{cc}
\alpha & \beta \\
-\beta^T& \gamma 
\end{array}\right).
\end{eqnarray}
Then the matrix $\gamma$ is invertible if and only if $a$ is invertible. In addition, we have
\begin{eqnarray}\label{apD.7}
\gamma^{-1}=c+b^Ta^{-1}b, 
\end{eqnarray}
\begin{eqnarray}\label{apD.7.1}
a^{-1}=\alpha+\beta\gamma^{-1}\beta^T. 
\end{eqnarray}

\noindent {\bf Proof.} Eqs. (\ref{apD.7}) and (\ref{apD.7.1}) immediately follow from the identity $AA^{-1}=1$, written in terms of the blocks. $\blacksquare$  

\noindent {\bf Affirmation 6.6.} The matrix (\ref{sim.23.02}) obeys the properties  (\ref{sim.G5}) and (\ref{sim.H5}) if and only if 
\begin{eqnarray}\label{apD.7.1.10}
\det\{\Phi^\alpha, \Phi^\beta\}\equiv\det\triangle^{\alpha\beta}\ne 0,  \quad \mbox{on} \quad  \mathbb{N}_k.
\end{eqnarray}

\noindent {\bf Proof.}  Consider the problem  in the coordinates of $\mathbb{M}_n$  
\begin{eqnarray}\label{apD.2}
y^k=(y^\alpha, y^a), \quad y^\alpha=\Phi^\alpha(z^\beta, z^b), \quad y^a=z^a, \quad a=1, 2, \ldots, k,
\end{eqnarray}
adapted with the functions $\Phi^\alpha$ (see Sect. \ref{SymNot_C}). Denote $\omega^{ij}(z^k)$ the Poisson tensor of $\mathbb{M}_n$.  Using the transformation law  (\ref{sim.4}), we obtain 
\begin{eqnarray}\label{apD.8}
\omega'^{ij}(y^k)=\left.\left(
\begin{array}{cc}
 \{\Phi^\alpha, \Phi^\beta\}& \{\Phi^\alpha, z^b\} \\
\{z^a, \Phi^\beta\}& \{z^a, z^b\} 
\end{array}\right)\right|_{z^i(y^j)}, \quad \mbox{denote its inverse as} \quad 
\tilde\omega'_{ij} (y^k)=\left(
\begin{array}{cc}
\tilde\omega'_{\alpha\beta}(y^k) & \tilde\omega'_{\alpha b}(y^k) \\
\tilde\omega'_{a\beta}(y^k) & \tilde\omega'_{ab}(y^k) 
\end{array}\right).
\end{eqnarray}

For the latter use we make the following observation. The symplectic matrix $\tilde\omega'_{ij}(y^\alpha, y^a)$ obeys the identity (\ref{sim.H5}). In particular, we have $\partial_a\tilde\omega'_{bc}(y^\alpha, y^a)+cycle=0$ for any fixed $y^\alpha$. Considering $\tilde\omega'_{bc}(y^a, y^\alpha)$ as a function of $y^a$, and applying the Affirmation 6.1, we conclude that its inverse, 
say $\omega_D^{ab}$, 
obeys the identity $\omega_D^{ad}\partial_d \omega_D^{bc}+ cycle(a, b, c)=0$. Using Affirmation 6.5 for  the matrices  (\ref{apD.8}), the explicit form of the inverse matrix is
\begin{eqnarray}\label{apD.80}
\omega_D^{ab}(y^\alpha, y^c)=\left.\left(\{z^a, z^b\}-\{z^a, \Phi^\alpha\}\tilde\triangle_{\alpha\beta}\{\Phi^\beta, z^b\}\right)\right|_{z^i(y^j)}. 
\end{eqnarray}

Let us return back to the proof. In adapted coordinates, the embedding $\mathbb{N}_k\rightarrow \mathbb{M}_n$ is given by $x^a\rightarrow y^i=(y^\alpha, y^a)$, where $y^\alpha=0$ and $y^a=x^a$. The equation (\ref{sim.23.02}) reads 
\begin{eqnarray}\label{apD.80.1}
\tilde\omega_{{\bf f} ab}(x^a)=\tilde\omega'_{ab}(y^\alpha, y^a)|_{y^\alpha=0,  y^a\rightarrow x^a},  
\end{eqnarray}
that is  the restriction of $\tilde\omega'_{ij} (y^k)$ on $\mathbb{N}_k$ reduces to the setting $y^\alpha=0$ in $a, b$\,-block of the matrix $\tilde\omega'_{ij} (y^k)$.  We need to confirm that $\tilde\omega_{{\bf f}}$ is a nondegenerate and closed form. 
The symplectic matrix $\tilde\omega'_{ij}(y^\alpha, y^a)$ obeys the identity (\ref{sim.H5}). In particular, we have $\partial_a\tilde\omega'_{bc}(y^\alpha, y^a)+cycle=0$ for any fixed $y^\alpha$. Taking $y^\alpha=0$, we conclude that $\tilde\omega_{{\bf f}}$ is closed.   Further, using Affirmation 6.5 for 
the matrices  (\ref{apD.8}), we conclude that  the matrix  $\tilde\omega_{{\bf f}}$ is invertible if and only if $\det\{\Phi^\alpha, \Phi^\beta\}\ne 0$. $\blacksquare$

As the restriction (\ref{sim.23.02}) determines a symplectic structure on $\mathbb{N}_k$, its inverse gives a Poisson bracket. Its explicit expression in terms of the original bracket can be obtained using the representation (\ref{apD.80.1}) for $\tilde\omega_{{\bf f} ab}$. Using Affirmation 6.5 for  the matrices  (\ref{apD.8}) and Eq. (\ref{sim.2.2.0}),  we can write for the inverse of $\tilde\omega_{{\bf f} ab}$ the expression 
\begin{eqnarray}
\omega_{{\bf f}}^{ab}(x^c)=\left.\left.\left(\{z^a, z^b\}-\{z^a, \Phi^\alpha\}\tilde\triangle_{\alpha\beta}\{\Phi^\beta, z^b\}\right)\right|_{z^i(y^j)}\right|_{y^\alpha=0, y^a\rightarrow x^a}  \qquad \quad \label{sim.23.03} \\
=\left.\left.\left(\{z^a, z^b\}-\{z^a, \Phi^\alpha\}\tilde\triangle_{\alpha\beta}\{\Phi^\beta, z^b\}\right)
\right|_{z^\alpha\rightarrow f^\alpha(z^a)}\right|_{z^a\rightarrow x^a}.  \label{sim.23.04}
\end{eqnarray}
Thus we obtained  the following result. \par

\noindent {\bf Affirmation 6.7.} Let $\omega^{ij}=\{z^i, z^j\}$ be non degenerate Poisson tensor and $\Phi^\alpha$ be functionally independent functions with $\det \{\Phi^\alpha, \Phi^\beta\}\ne 0$. Then the matrix 
\begin{eqnarray}\label{apD.12}
\omega_{{\bf f}}^{ab}(z^c)=\left.\left(\{z^a, z^b\}-\{z^a, \Phi^\alpha\}\tilde\triangle_{\alpha\beta}\{\Phi^\beta, z^b\}\right)\right|_{z^\alpha=f^\alpha(z^c)}, 
\end{eqnarray}
where $z^\alpha=f^\alpha(z^c)$ are parametric equations of the surface $\Phi^\alpha=0$, obeys the Jacobi identity, and determines a non degenerate Poisson bracket on $\mathbb{N}_k$
\begin{eqnarray}\label{apD.12.1}
\{ A, B \}_{D({\mathbb N})}=\partial_a A \omega^{ab}_{\bf f}\partial_b B. 
\end{eqnarray}

There is a bracket on  $\mathbb{M}_n$ that induces\footnote{With this respect, see the comment at the end of Sect. \ref{I D}.}  the bracket (\ref{apD.12.1}) on $\mathbb{N}_k$ according to Eq.(\ref{dar.6}). The equality (\ref{apD.12}) prompts to consider 
\begin{eqnarray}\label{sim.26}
\{ A, B\}_D=\{A, B\}-\{A, \Phi^\alpha\}\tilde\triangle_{\alpha\beta}\{\Phi^\beta, B\}=
\partial_i A\left[\{z^i, z^j\}-\{z^i, \Phi^\alpha\}\tilde\triangle_{\alpha\beta}\{\Phi^\beta, z^j\}\right]\partial_j B\equiv\partial_i A~\omega_D^{ij}~\partial_j B.
\end{eqnarray}
This is the famous Dirac bracket \cite{Dir_1950}. The tensor $\omega_D^{ij}(z^k)$ obeys the Jacobi identity (see below), and hence turn $\mathbb{M}_n$ into the Poisson manifold $(\mathbb{M}_n, ~ \{{}, {}\}_D)$. The bracket (\ref{apD.12.1}) turns out to be the restriction of 
(\ref{sim.26}) to $\mathbb{N}_k$. To see this, we first note that for any function $A(z^i)$, Eq. (\ref{sim.26})  implies
\begin{eqnarray}\label{sim.27}
\{A, \Phi^\alpha\}_D=0, 
\end{eqnarray}
so $\Phi^\alpha$ are Casimir functions of the Dirac bracket. As we saw in  Sect. \ref{DegPoi}, this implies that all Hamiltonian fields $V_A^i=\omega_D^{ij}\partial_j A$ are tangent to the surfaces $\Phi^\alpha=c^\alpha$, and  we can restrict the Dirac tensor 
$\omega_D^{ij}$  on the submanifold $\mathbb{N}_k$ according to  Eq. (\ref{dar.6}). This gives the Poisson bracket (\ref{apD.12.1}) on $\mathbb{N}_k$  and turns it into a Poisson submanifold of the Poisson manifold $(\mathbb{M}_n, \{{}, {}\}_D)$. 

It remains to prove the Jacobi identity for the Dirac bracket. 

\noindent {\bf Affirmation 6.8.} Consider Poisson manifold $\mathbb{M}_n=\{z^i, ~ \omega^{ij}(z^k)\}$ with a non degenerate tensor $\omega$. Let $\Phi^\alpha(z^k)$ be functionally independent functions which obey the condition (\ref{apD.7.1.10}). Then the Dirac tensor $\omega_D^{ij}(z^k)$, specified in (\ref{sim.26}), satisfies the identity (\ref{sim.E5}), and hence the Dirac bracket (\ref{sim.26})
satisfies the Jacobi identity: $\{A, \{B, C\}_D\}_D+cycle ~ (A, B, C)=0$. 

\noindent {\bf Proof.}    Consider the problem in the coordinates (\ref{apD.2})  adapted with the functions $\Phi^\alpha$. Using Eqs. (\ref{sim.4}) and (\ref{sim.27}), we obtain the Dirac tensor in these coordinates
\begin{eqnarray}\label{apD.3}
\omega_D '^{ij}(y^k)=\left(
\begin{array}{cc}
0 & 0 \\
0& \left.\omega_D^{ab}(z^k)\right|_{z(y)}
\end{array}\right), 
\end{eqnarray}
where $\omega_D^{ab}(z^k)$ is $a, b$\,-block of the Dirac tensor $\omega_D^{ij}(z^k)$ in original coordinates, and then  $\left.\omega_D^{ab}(z^k)\right|_{z(y)}$ is just the expression written in (\ref{apD.80}). The desired Jacobi identity  $\omega_D '^{in}\partial_n \omega_D '^{jk}+cycle=0$ will be fulfilled, if the matrix (\ref{apD.80}) obeys the identity $\omega_D^{ad}\partial_d \omega_D^{bc}+cycle=0$. But this was confirmed above, see the discussion below of Eq. (\ref{apD.8}). $\blacksquare$

The results of this subsection can be summarized  in the form of diagram (\ref{diag}), that relates geometrical structures on the manifold $\mathbb{M}_n$ (top line),  and on its submanifold $\mathbb{N}_k$ (bottom line): 
\begin{eqnarray}\label{diag}
\begin{array}{ccccc}
\tilde\omega_{ij} & \longleftrightarrow & \omega^{ij} & \longrightarrow & \omega_D^{ij}\sim \{ {}, {}\}_D\\
\downarrow & {} & {} & {} & \downarrow \\
\tilde\omega_{{\bf f} ab} & \longleftarrow  & -- &  \longrightarrow & \omega_{\bf f}^{ab}\sim\{ {}, {}\}_{D(\mathbb N)}
\end{array}
\end{eqnarray}
The Dirac bracket appears in the upper right corner of the rectangle, and provides the closure of our diagram. 

Discussion of the Dirac bracket in the coordinate-free language can be found in \cite{Maslov_1993, Cou_1990, Sni_1974, Zam_2009, Bur_2011, Mei_2016}.

\subsection{Dirac's derivation of the Dirac bracket.}\label{VI C}  

Dirac arrived at his bracket in the analysis of a variational problem for singular nondegenerate theories like (\ref{sim.0.5}), (\ref{sim.0.6}). Consider the variational problem 
\begin{eqnarray}\label{sima.12.1}
S=\int d\tau\left[p_a\dot q^a-H_0(q^a, p_b)+\lambda^\alpha\Phi_\alpha(q^a, p_b)\right],  
\end{eqnarray}
for the set of independent dynamical variables $z^i(\tau)\equiv (q^a, p_b)$, $i=(1, 2, \ldots , 2n)$,  and $\lambda^\alpha(\tau)$, $\alpha=(1,  2, \ldots , 2p<2n)$. $H_0$ and $\Phi_\alpha$ are given functions, where  $\Phi_\alpha$ obey the conditon (\ref{apD.7.1.10}).  Variation of the action with respect to $z^i$ and $\lambda^\alpha$ gives the equations of motion\footnote{It is instructive to  compare the systems (\ref{sima.12.2}) and (\ref{sim.16.5}). The constraints $\Phi_\alpha=0$ should not be confused with the first integrals. Indeed, first integrals represent the first-order differential equations which are consequences of a special form of the original equations, $c_{\alpha i}[\dot z^i-\{z^i, H\}]=\frac{d}{d\tau}Q_\alpha(z)=0$, whereas constraints are the algebraic equations. As a consequence, solutions of the systems (\ref{sima.12.2}) and (\ref{sim.16.5}) have very different properties. Solutions of the 
system (\ref{sim.16.5}) pass through any point of $\mathbb{R}$, while all solutions of (\ref{sima.12.2}) live on the submanifold $\Phi_\alpha=0$.}
\begin{eqnarray}\label{sima.12.2}
\dot z^i=\{z^i, H_0\}+\lambda^\alpha\{z^i, \Phi_\alpha\}, \qquad \Phi_\alpha=0, 
\end{eqnarray}
where $\{{}, {}\}$ is the canonical Poisson bracket on $\mathbb{R}_{2n}$. Let $z^i(\tau)$, $\lambda(\tau)$ be a solution of the system. Computing derivative of the identity $\Phi_\alpha(z^i(\tau))=0$ , we obtain the algebraic equations  
\begin{eqnarray}\label{sima.12.3}
\{\Phi_\alpha, H_0\}+\{\Phi_\alpha, \Phi_\beta\}\lambda^\beta=0. 
\end{eqnarray}
that must be satisfied for all solutions, that is they are the consequences of the system. According to this equation, all variables $\lambda^\beta$ are determined algebraically, $\lambda^\beta=-\tilde\triangle^{\beta\alpha}\{\Phi_\alpha, H_0\}$, where $\tilde\triangle$ is the inverse matrix of $\triangle$.  Adding the consequences to the system, we obtain its equivalent form 
\begin{eqnarray}
\dot z^i=\{z^i, H_0\}-\{z^i, \Phi_\alpha\}\tilde\triangle^{\alpha\beta}\{\Phi_\beta, H_0\}\equiv\left[\omega^{ij}-\{z^i, \Phi_\alpha\}\tilde\triangle^{\alpha\beta}\{\Phi_\beta, z^j\}\right]\partial_j H_0, \label{sima.12.4} \\ 
\Phi_\alpha=0, \qquad \lambda^\beta=-\tilde\triangle^{\beta\alpha}\{\Phi_\alpha, H_0\}, \label{sima.12.4.1} \qquad \qquad \qquad \qquad 
\end{eqnarray}
where the sectors $z^i$ and $\lambda^\beta$ turn out to be separated.  The expression appeared on r.h.s. of (\ref{sima.12.4}) suggests to introduce the new bracket on $\mathbb{M}_{2n}$ 
\begin{eqnarray}\label{sima12.5}
\{ A, B\}_D=\{A, B\}-\{A, \Phi_\alpha\}\tilde\triangle^{\alpha\beta}\{\Phi_\beta, B\},
\end{eqnarray}
which is just the Dirac bracket.  Then equations (\ref{sima.12.4}) represent a Hamiltonian system with the Dirac bracket
\begin{eqnarray}\label{sima.12.6}
\dot z^i=\{z^i, H_0\}_D, 
\end{eqnarray}
and with the Hamiltonian being $H_0$.

\section{Poisson manifold and Dirac bracket.}\label{SymDirac}  
\subsection{Jacobi identity for the Dirac bracket.}\label{SymDiracA} 
While our discussion of the Dirac bracket  in previous section was based on a symplectic manifold, the prescription (\ref{sim.26}) can  equally be used to generate a bracket $\{A, B\}_D$ starting from a given degenerate Poisson bracket $\{ A, B\}$. We show that  $\{A, B\}_D$ still satisfies the Jacobi identity. To prove this, we will need the following auxiliary statement. 
 
\noindent {\bf Affirmation 7.1.} Consider the Poisson manifold $\mathbb{M}_{m+n}=\{x^K=(x^{\bar\alpha}, x^i), ~ \omega^{IJ}(x^K), ~ rank ~ \omega=n\}$. Let $K^{\bar\alpha}(x^I)$, $\bar\alpha=1, 2, \ldots , m$  be functionally independent Casimir functions, and $\Phi^\alpha(x^I)$, $\alpha=1, 2, \ldots , p< n$ be functionally independent functions which obey the condition (\ref{apD.7.1.10}). Then 

\noindent {\bf (A)} The $m+p$ functions $K^{\bar\alpha}, \Phi^\beta$ are functionally independent. 

\noindent {\bf (B)} In the coordinates 
\begin{eqnarray}\label{sip.0.1}
z^I=(z^{\bar\alpha}, z^i), \quad \mbox{where} \quad z^{\bar\alpha}=K^{\bar\alpha}(x^I), \quad z^i=x^i, 
\end{eqnarray}
the functions $\Phi^\alpha(z^{\bar\alpha}, z^i)$ obey the condition $rank ~ \frac{\partial \Phi^\alpha}{\partial z^i}=p$. In other words $\Phi^\alpha$, considered as functions of $z^i$, are functionally independent. 

\noindent {\bf Proof.} {\bf (A)} In the coordinates (\ref{sip.0.1}), our functions are $z^{\bar\alpha}$ and $\Phi^\alpha(z^{\bar\alpha}, z^i)$. We will show that functional dependence of the set implies that the matrix $\{\Phi^\alpha, \Phi^\beta\}$ is degenerate. Then  non degeneracy implies functional independence of the set  - the desired result. 

Consider $(m+p)\times (m+n)$ matrix 
\begin{eqnarray}\label{sip.0.2}
J=\frac{\partial(z^{\bar\alpha}, \Phi^\alpha(z^{\bar\alpha}, z^i))}{\partial(z^{\bar\alpha}, z^i)}=\left(
\begin{array}{cc}
{\bf 1}_{m\times m} & {\bf 0} \\
\frac{\partial\Phi^\alpha}{\partial z^{\bar\beta}} & \frac{\partial\Phi^\alpha}{\partial z^i}
\end{array}
\right). 
\end{eqnarray}
If $z^{\bar\alpha}, \Phi^\alpha(z^{\bar \alpha}, z^i)$ are functionally dependent, we have $rank ~ J<m+p$, then some linear combination of rows of the matrix $J$ vanishes: $c_{\bar\alpha}\delta^{\bar\alpha}{}_I+c_\alpha\frac{\partial\Phi^\alpha}{\partial z^I}=0$ for all $I$. This equation, together with explicit expression (\ref{sip.0.2}) for $J$, implies 
\begin{eqnarray}\label{sip.0.3}
c_\alpha\frac{\partial\Phi^\alpha}{\partial z^i}=0, \quad {\vec c}\ne 0.  
\end{eqnarray}
Consider now the matrix $\{\Phi^\alpha, \Phi^\beta\}$ in the coordinates (\ref{sip.0.1}). Using the Poisson tensor 
\begin{eqnarray}\label{sip.0.4}
\omega'^{IJ}(z^K)\equiv
\left(
\begin{array}{cc}
\omega'^{\bar\alpha \bar\beta} & \omega'^{\bar\alpha j} \\ 
\omega'^{i \bar\beta} & \omega'^{ij}
\end{array}
\right)=
\left(
\begin{array}{cc}
{\bf 0}_{m\times m} & {\bf 0} \\
{\bf 0} & \omega^{ij}(x^K)|_{x(z)}
\end{array}
\right), 
\end{eqnarray}
we obtain $\{\Phi^\alpha, \Phi^\beta\}=\partial_i\Phi^\alpha\omega'^{ij}\partial_j\Phi^\beta$. Then (\ref{sip.0.3}) implies the  degeneracy of the matrix: $\{\Phi^\alpha, \Phi^\beta\} c_\beta=0$. 

{\bf (B)} Item (A) implies: $rank ~ J=m+p$. Then from explicit form (\ref{sip.0.2}) for $J$ it follows, that $rank ~ \frac{\partial\Phi^\alpha}{\partial z^i}=p$. $\blacksquare$ 

\noindent {\bf Affirmation 7.2.} The Dirac bracket (\ref{sim.26}), constructed on the base of a degenerate Poisson bracket $\{A, B\}$, satisfies the Jacobi identity. 

\noindent {\bf Proof.} We use the notation specified in Affirmation 7.1. The original Poisson tensor in the coordinates (\ref{sip.0.1}) is  written in Eq. (\ref{sip.0.4}). According to Affiration 4.2, its block $\omega'^{ij}$ is a non degenerate matrix.  $\omega'^{IJ}(z^K)$ satisfies the Jacobi identity, that due to special form (\ref{sip.0.4}) of this tensor reduces to the expression
\begin{eqnarray}\label{sip.0.5}
\omega'^{in}\frac{\partial}{\partial z^n}\omega'^{jk}+cycle=0. 
\end{eqnarray} 
Using the prescription (\ref{sim.26}), we use $\omega'^{IJ}(z^K)$ to write the Dirac tensor  
\begin{eqnarray}\label{sip.0.6}
\omega_D '^{IJ}(z^K)=
\left(
\begin{array}{cc}
{\bf 0}_{m\times m} & {\bf 0} \\
{\bf 0} & \omega_D '^{ij}
\end{array}
\right), \quad \mbox{where} \quad \omega_D '^{ij}=\omega'^{ij}-\omega'^{in}\partial_n\Phi^\alpha\tilde\triangle_{\alpha\beta}\partial_k\Phi^\beta.  
\end{eqnarray}
Jacobi identity for $\omega_D '^{IJ}(z^K)$ will be satisfied, if 
\begin{eqnarray}\label{sip.0.7}
\omega_D '^{in}\frac{\partial}{\partial z^n}\omega_D '^{jk}+cycle=0. 
\end{eqnarray}
Note that $z^{\bar\alpha}$ enter into the expressions (\ref{sip.0.5})-(\ref{sip.0.7}) as the parameters. In particular, the derivative $\frac{\partial}{\partial z^{\bar\alpha}}$ falls out of all these expressions. According to Item (B) of Affirmation 7.1, the functions $\Phi^\alpha(z^{\bar\beta}, z^i)$, considered as functions of $z^i$, are functionally independent. Taking this into account, we can apply the Affirmation 6.8 to the matrices specified by  (\ref{sip.0.5}) and (\ref{sip.0.6}), and conclude that (\ref{sip.0.7}) 
holds. $\blacksquare$

\subsection{Some applications of the Dirac bracket.}\label{SymDiracB} 
 
With given scalar function $A$ we associate the function 
\begin{eqnarray}\label{sim.301}
A_d=A-\{A, \Phi^\alpha\}\tilde\triangle_{\alpha\beta}\Phi^\beta.  
\end{eqnarray}
The two functions coincide on the surface $\Phi^\alpha=0$. There is a remarcable relation between the Dirac bracket of original functions and the Poisson bracket of deformed functions:
\begin{eqnarray}\label{sim.302}
\{A, B\}_D=\{A_d, B_d\}+O(\Phi^\alpha),
\end{eqnarray}
which means that  the two brackets also coincide on the surface. This property can be reformulated in terms of vector fields as follows. Given scalar function $A$, integral lines of Hamiltonian field $V^i=\omega^{ij}\partial_j A_d$ that cross  the surface $\Phi^\alpha=0$, entirely lie on it. 

Below we use the Dirac bracket to analyse some Hamiltonian systems consisting of both dynamical and algebraic equations. 

{\bf 1.} Consider the Hamiltonian system $\dot z^i=\{z^i, H\}_D$ on the Poisson manifold $(\mathbb{M}, \{{}, {}\}_D)$. As $\{\Phi^\alpha, H\}_D=0$, the functions $\Phi^\alpha$ are integrals of motion of the system. According to Affirmation 5.2, all the submanifolds $\mathbb{N}_k^{\vec c}=\{z\in\mathbb{M}_n, \quad \Phi^\alpha(z)=c^\alpha\}$ are invariant submanifolds, that is any trajectory that starts on $\mathbb{N}_k^{\vec c}$ entirely lies on it. In particular, we have  \par

\noindent {\bf Affirmation 7.3.} The equations 
\begin{eqnarray}\label{sim.30}
\dot z^i=\{z^i, H\}_D,  \qquad \Phi^\alpha=0,  
\end{eqnarray}
form a self-consistent system in the sense of Definition 1.1. 

Further, according to Affirmation 5.3, these equations are equivalent to the system: $z^\alpha-f^\alpha(z^a)=0$, $\dot z^b=\left.\{z^b, H(z^i)\}_D\right|_{z^\alpha\rightarrow f^\alpha(z^b)}$. We replace $z^\alpha$ on $f^\alpha(z^b)$ using Eqs. (\ref{dar.5.1}) and (\ref{dar.6}), this gives 
\begin{eqnarray}\label{sim.31}
z^\alpha=f^\alpha(z^b), \quad \dot z^b=\{z^b, \hat H(z^b)\}_{D({\mathbb N})}, 	
\end{eqnarray}
where $\hat H(z^b)=H(z^b, z^\alpha(z^b))$, and $\{ {}, {} \}_{D({\mathbb N})}$ is the bracket (\ref{apD.12}) on ${\mathbb N}$, induced by the Dirac bracket. 
This shows that the variables $z^b$ obey the Hamiltonian equations on the submanifold  $\mathbb{N}$.  

{\bf 2.} Let us rewrite the system (\ref{sim.30}) in terms of original bracket as follows: 
$\dot z^i=\{z^i, H-\Phi^\alpha\tilde\triangle_{\alpha\beta}\{\Phi^\beta, H\}\}+
\Phi^\alpha\{z^i, \tilde\triangle_{\alpha\beta}\{\Phi^\beta, H\}\}$, $\Phi^\alpha=0$, or, equivalently
\begin{eqnarray}\label{sim.32}
\dot z^i=\{z^i, H-\Phi^\alpha\tilde\triangle_{\alpha\beta}\{\Phi^\beta, H\}\}, 
\qquad \Phi^\alpha=0.  
\end{eqnarray}
Note that the functions $\Phi^\alpha$ are not the Casimir functions of the original bracket.  As the systems (\ref{sim.32}) and (\ref{sim.30}) are equivalent, we obtained an example of a self-consistent theory of the type  (\ref{sim.0.5}),  (\ref{sim.0.6}): \par

\noindent {\bf Affirmation 7.4.} Given Poisson manifold $(\mathbb{M}_n, \{{}, {}\})$, let $H$ be given function and $\Phi^\alpha$ is a set of functionally independent functions that obey the condition $\left.\det \{ \Phi^\alpha, \Phi^\beta \}\right|_{\Phi^\alpha=0}\equiv \det \triangle^{\alpha\beta}\ne 0$. Then the equations 
\begin{eqnarray}\label{sim.33}
 \dot z^i=\{z^i, \tilde H\},  \qquad \Phi^\alpha=0.  
\end{eqnarray}
with the Hamiltonian $\tilde H=H-\Phi^\alpha\tilde\triangle_{\alpha\beta}\{\Phi^\beta, H\}$ form a self-consistent system. 

{\bf 3.} {\bf Affirmation 7.5.} Singular nondegenerate theory defined by the equations (\ref{sim.0.5}),  (\ref{sim.0.6}) with the properties   (\ref{sim.0.5.1}) and (\ref{sim.0.5.2})  is a self-consistent, and is equivalent to (\ref{sim.30}).  \par

\noindent {\bf Proof.} Using (\ref{sim.0.5.1}), we rewrite the system (\ref{sim.0.5}),  (\ref{sim.0.6}) in the equivalent form as follows
\begin{eqnarray}\label{sim.34}
\dot z^i=\{z^i, H\}_D+\{z^i, \Phi^\alpha\}\tilde\triangle_{\alpha\beta}\{\Phi^\beta, H\}, 
\qquad \Phi_\alpha=0.  
\end{eqnarray}
Take any point of the submanifold $\Phi^\alpha=0$. According to Affirmation 7.3, there is a solution $z^i(\tau)$ of (\ref{sim.30}) that passes through this point. Due to the condition (\ref{sim.0.5.2}) we have $\left.\{\Phi^\beta, H\}\right|_{z^i(\tau)}=0$, then the direct substitution of $z^i(\tau)$ into (\ref{sim.34}) shows that it is a solution of this system. $\blacksquare$

{\bf 4.} {\bf Example of Hamiltonian reduction.} Let the Hamiltonian system $\dot z^i=\omega^{ij}\partial_j H(z^k)$ with $\det\omega\ne 0$ admits the first integrals $\Phi^\alpha(z^k)$, $\alpha=1, 2, \ldots , n-k$ with the properties 
\begin{eqnarray}\label{sim.34.1}
\{\Phi^\alpha, H\}=0, \qquad \det\{\Phi^\alpha, \Phi^\beta\}\equiv\triangle^{\alpha\beta}\ne 0.  
\end{eqnarray}
Then dynamics can be consistently restricted on any one of invariant surfaces $\mathbb{N}_k^{\vec c}=\{z^k\in\mathbb{M}_n, ~ \Phi^\alpha=c^\alpha \}$. Without loss of generality, we consider reduction on $\mathbb{N}_k^{\vec 0}$, then $\Phi^\alpha=0$ implies $z^\alpha=f^\alpha(z^a)$, while the independent variables obey the equations
\begin{eqnarray}\label{sim.34.2}
\dot z^a=\left.\omega^{aj}\partial_j H\right|_{z^\alpha=f^\alpha(z^a)}. 
\end{eqnarray}
We do the substitution indicated in this equation and show that the result is a Hamiltonian system. Consider the problem in the adapted coordinates (\ref{apD.2}). Then $\Phi^\alpha=0$ turn into $y^\alpha=0$ while instead of (\ref{sim.34.2}) we have
\begin{eqnarray}\label{sim.34.3}
\dot y^a=\left.\omega'^{aj}\partial_j H\right|_{y^\alpha=0}=\left.\omega'^{ab}\partial_b H'\right|_{y^\alpha=0}+\left.\omega'^{a\beta}\partial_\beta H'\right|_{y^\alpha=0},
\end{eqnarray}
where explicit form of $\omega'$ is given by (\ref{apD.8}), and $H'(y^i)=H(z^k(y^i))$.The equation $\{\Phi^\alpha, H\}=0$ in the coordinates $y^k$ gives: $0=\{y^\alpha, H'\}=\omega'^{\alpha a}\partial_a H'+\triangle^{\alpha\beta}\partial_\beta H'$, or $\partial_\beta H'=-\tilde\triangle_{\beta\gamma}\omega'^{\gamma a}\partial_a H'$. Using this expression in (\ref{sim.34.3}) we obtain
\begin{eqnarray}\label{sim.34.4}
\dot y^a=\left.\left.\left[\{z^a, z^b\}-\{z^a, \Phi^\beta\}\tilde\triangle_{\beta\gamma}\{\Phi^\beta, z^b\}\right]\partial_b H(z^i)\right|_{z^i(y^j)}\right|_{y^\alpha=0}.
\end{eqnarray}
Now note that $\left.A(z^i(y^j))\right|_{y^\alpha=0}=\left.A(f^\alpha(z^a),  z^a)\right|_{z^a\rightarrow y^a}$, so the equations of motion reads
\begin{eqnarray}\label{sim.34.5}
\dot z^a=\omega_D^{ab}(f^\alpha(z^a),  z^a)\partial_b H(f^\alpha(z^a), z^a),
\end{eqnarray}
where $\omega_D^{ab}(f^\alpha(z^a) z^a)$ is $(a,b)$\,-block of the Dirac tensor, see  (\ref{apD.12}). According to Affirmation 6.7, it obeys the Jacobi identity, so the equations  (\ref{sim.34.5}) represent a Hamiltonian system, which is equivalent to (\ref {sim.34.2}).

\subsection{Poisson manifold with prescribed Casimir functions.}\label{PresCas}

Let $K^\alpha(z^\beta, z^ b)$ with $\det\frac{\partial K^\alpha}{\partial z^\beta}\ne 0$ scalar functions  in local coordinates $z^i=(z^\beta, z^ b)$ of the manifold $\mathbb{M}_n$, where $\beta=1, 2, \ldots, p$, $b=1, 2, \ldots , n-p$. Without loss of generality, we assume that $n-p$ is an even number: $n-p=2k$. The task is to construct a Poisson bracket on $\mathbb{M}_n$, that has $K_\alpha$ as the Casimir functions.  One possible solution of this task can be found by using of coordinate system where the functions $K_\alpha$ turn into a part of coordinates. 

Introduce the following coordinates on $\mathbb{M}_n$:
\begin{eqnarray}\label{sim.36}
z^{j'}=\varphi^{j'}(z^i)=(K^\alpha(z^i), z^a). 
\end{eqnarray}
Construct the matrix $a$ with elements $a_i{}^{j'}=\frac{\partial z^{j'}}{\partial z^i}=\partial_i\varphi^{j'}$, its inverse is denoted  as $\tilde a\equiv a^{-1}$. In the local coordinates $z^{j'}$, define the bracket 
\begin{eqnarray}\label{sim.35}
\{ A, B\}=\partial_{i'}A ~ W_0^{i'j'}\partial_{j'} B, \qquad 
W_0^{i'j'}(z^{i'})=\left(
\begin{array}{cc}
0_{p\times p}  & 0  \\
0  & \omega_0(z^{j'})
\end{array}
\right),   
\end{eqnarray}
where $\omega_0$ is a $2k\times 2k$ matrix with the elements $\omega_0^{a'b'}(z^{\alpha'}, z^{c'})$ satisfying the identity (\ref{sim.E5}) with respect to $z^{c'}$. As this matrix we can take any known Poisson structure $\omega_0^{a'b'}(z^{c'})$ on the submanifold $K^\alpha(z^\beta, z^ b)=0$. For instance, we could take it in the canonical form
\begin{eqnarray}\label{sim.35.1}
\omega_0^{a'b'}=\left(
\begin{array}{cc}
0_{k\times k} & 1_{k\times k}  \\
-1_{k\times k} & 0_{k\times k} 
\end{array}
\right).   
\end{eqnarray}
According to Eq. (\ref{sim.4}), in the original coordinates $z^i$ the bracket reads
\begin{eqnarray}\label{sim.37}
\{ A, B\}=\partial_i A\omega^{ij}\partial_j B, \qquad \omega^{ij}=\left[\tilde a^T W_0(K^\alpha(z^i), z^a) \tilde a\right]^{ij}. 
\end{eqnarray}
The Affirmation 2.2 guarantees that it satisfies the Jacobi identity, hence it turn $\mathbb{M}_n$ into a Poisson manifold. 

\noindent {\bf Affirmation 7.6.}  $K^\alpha$ are Casimir functions of the bracket  (\ref{sim.37}).  \par 

\noindent {\bf Proof.} Consider, for instance, $\{A, K^1\}=\partial_i A(\tilde a^T W_0 \tilde a)^{ij}\partial_j K^1$. Compute the term: $(W_0 \tilde a)^{ij}\partial_j K^1=(W_0 \tilde a)^{ij}a_j{}^1=W_0^{ij}\delta_j{}^1=W_0^{i1}=0$. $\blacksquare$ 

In resume, the set of functionally independent functions $K^\alpha(z^i)$ turns out to be  the set of Casimir functions of the Poisson manifold  with the bracket (\ref{sim.37}).

Denoting $\partial_\alpha K^\beta=b_{\alpha}{}^{\beta}$, $\partial_a K^\beta=c_a{}^{\beta}$, the Poisson structure (\ref{sim.37}) can be written in the following form:
\begin{eqnarray}\label{sim.38}
\omega=\left(
\begin{array}{cc}
(cb^{-1})^T\omega_0 cb^{-1} & (\omega_0 cb^{-1})^T \\
-\omega_0 cb^{-1} & \omega_0 
\end{array}
\right).   
\end{eqnarray}
Blocks of this matrix can be compared with Eqs. (\ref{dar.4.1}).  We can restrict the bracket (\ref{sim.37}) on the Casimir submanifold, obtaining the bracket (see Eq. (\ref{dar.6}))
\begin{eqnarray}\label{sim.39}
\{A(z^a), B(z^a)\}=\partial_a A\bar\omega^{ab}\partial_b B, \qquad \bar\omega^{ab}=\omega_0^{ab}(f^\alpha(z^a), z^a). 
\end{eqnarray}
In particular, if $\omega_0$ in Eq. (\ref{sim.35}) was originally chosen to be independent of the coordinates $z^\alpha$, we have $\bar\omega^{ab}=\omega_0^{ab}$.  
Casimir submanifold with the bracket (\ref{sim.39}) is the Poisson submanifold of $\mathbb{M}_n$ (\ref{sim.37}) in the sense of definition (\ref{sim15.10}). 

{\sc Example 7.1.}  Consider $\mathbb{M}_3$ and the function $K(z^1, z^2, z^3)$ with $\partial_1 K\ne 0$. Then 
\begin{eqnarray}\label{sim.40}
a=\left(
\begin{array}{ccc}
\partial_1 K & 0 & 0 \\
\partial_2 K & 1 & 0 \\
\partial_3 K & 0 & 1
\end{array}
\right), \qquad \tilde a=\frac{1}{\det a}\left(
\begin{array}{ccc}
1 &0 & 0 \\
-\partial_2 K & 1 & 0 \\
-\partial_3 K & 0 & 1
\end{array}
\right).   
\end{eqnarray}
Taking 
\begin{eqnarray}\label{sim.41}
W=\left(
\begin{array}{ccc}
0 & 0 & 0 \\
0 & 0 & 1 \\
0 & -1 & 0
\end{array}
\right),   
\end{eqnarray}
we obtain the Poisson structure on $\mathbb{M}_3$ that has $K(z)$ as the Casimir function
\begin{eqnarray}\label{sim.42}
\omega=\tilde a^T W\tilde a=\frac{1}{\det a}\left(
\begin{array}{ccc}
0 & \partial_3 K & -\partial_2 K \\
-\partial_3 K & 0 & 1 \\
\partial_2 K & -1 & 0
\end{array}
\right),   \quad \mbox{or} \quad \omega^{ij}=\frac{1}{\det a}\epsilon^{ijk}\partial_k K. 
\end{eqnarray}
If $V_i$ and $U_j$ are contravariant vectors, the quantity $\frac{1}{\det a}\epsilon^{ijk}\partial_k K ~ V_i ~ U_j$ is a scalar function under the diffeomorphisms (\ref{sim.1}). So $\omega^{ij}$ of Eq.  (\ref{sim.42}) is a second-rank covariant tensor, as it should be. Restriction of the bracket (\ref{sim.42}) on the Casimir submanifold $K=0$ gives the canonical Poisson bracket: $\{z^2, z^3\}=1$. 

{\sc Example 7.2.} $SO(3)$ Lie-Poisson bracket. Chosing $K=\frac12[(z^1)^2+(z^2)^2+(z^3)^2]-1$ (see Example 1.2) in the expressions of previous example, we obtain diffeomorphism-covariant form of the Lie-Poisson bracket: 
\begin{eqnarray}\label{sim.43}
\{z^i, z^j \}=\frac{1}{\det a}\epsilon^{ijk}z^k. 
\end{eqnarray}

\section{Discussion.}\label{Concl}

In this short survey, we presented an elementary exposition of the methods of Poisson and symplectic geometry, 
with an emphasis on the construction, geometric meaning and applications of the Dirac bracket. We have traced the role played by the Dirac bracket in the problem of reducing the Poisson structure of a manifold to the submanifold defined by scalar functions, which form the set of second-class constraints. Then the Dirac bracket was applied to the study of the Hamiltonian system (\ref{sim.0.5}), (\ref{sim.0.6}) with second-class constraints (\ref{sim.0.5.1}). Let us briefly describe these results.

Let $\mathbb{M}_n=\{z^k, \omega^{ij}(z^k)\}$ is a nondegenerate Poisson manifold and $\mathbb{N}_k=\{x^a\}$ is a submanifold determined by the equations  $\Phi^\beta(z^k)=0$, $n$ and $k$ are even numbers. Let  $z^\beta=f^\beta(z^a)$ be the solution to these equations. They determine the embedding $\mathbb{N}_k\rightarrow \mathbb{M}_n$, given by $x^a\rightarrow z^i=(f^\beta(x^a), x^a)$. Nondegenerate contravariant tensor $\omega^{ij}$ can not be directly used to unduce the Poisson structure on the submanifold. But we can do this with help of symplectic form $\tilde\omega_{ij}$, corresponding to the Poisson tensor $\omega^{ij}$. For the case of submanifold determined by the second-class constraints,  $\left.\det \{\Phi^\alpha, \Phi^\beta\}_P\right|_{\Phi^\alpha=0}\ne 0$,  the induced mapping (\ref{sim.23.02}) determines the symplectic form $\tilde\omega_{{\bf f}ab}(x^c)$ on $\mathbb{N}_k$. The explicit form of inverse of this matrix is given by  Eq. (\ref{apD.12}), and determines a nondegenerate Poisson bracket $\{ A, B \}_{D({\mathbb N})}=\partial_a A \omega^{ab}_{\bf f}\partial_b B$ on $\mathbb{N}_k$. This solves the reduction problem. 

Next, we may wonder about constructing a degenerate Poisson bracket on $\mathbb{M}_n$ that directly induces the bracket on $\mathbb{N}_k$ with use the Casimir functions, see  Eq. (\ref{dar.6}). The explicit form (\ref{apD.12}) of the Poisson tensor $\omega^{ab}_{\bf f}(x^c)$ immediately prompts the Dirac bracket (\ref{sim.26}) as a solution of this task. The described construction can be resumed in the form of diagram (\ref{diag}). The Dirac bracket appears in the upper right corner of the rectangle, and provides the closure of the diagram.

Consider now the Hamiltonian system (\ref{sim.0.5})-(\ref{sim.0.5.2}) on $\mathbb{M}_n$ and the following Hamiltonian system on $\mathbb{N}_k$: 
\begin{eqnarray}\label{Con1}
\dot x^a=\{x^a, H(x^b)\}_{D({\mathbb N})},  \qquad  H(x^b)=H(f^\alpha(x^b), x^b).   
\end{eqnarray}
Using the Dirac bracket, we demonstrated in Sect. \ref{SymDiracB} that the two systems are equivalent. This implies that the system (\ref{sim.0.5}), (\ref{sim.0.6}) with second-class constraints  (\ref{sim.0.5.1}), (\ref{sim.0.5.2}) is a self-consistent, and its restriction on $\mathbb{N}_k$ is a Hamiltonian system.

\vspace{5mm}

\noindent{\bf Acknowledgments.} 
The work has been supported by the Brazilian foundation CNPq,  and by Tomsk State University Competitiveness Improvement
Program.

\section{Appendices.}\label{App1}

\subsection{Jacobi identity}\label{App A}

\noindent {\bf Affirmation A 1.} Let the bracket (\ref{sim.5}) obeys the Jacobi identity in the coordinates $z^i$. Then the Jacobi identity is satisfied in any other coordinates. \par

\noindent {\bf Proof.} We need to show that the validity of the identity (\ref{sim.B5}) for the bracket (\ref{sim.5}) with $\omega^{ij}(z)$ implies its validity for this bracket with $\omega^{i'j'}$ defined in (\ref{sim.4}). 

Given functions $A(z), B(z), C(z)$, let us consider the auxiliary functions $\tilde A(z)\equiv A(z'(z))$ and so on. Since the Jacobi identity is satisfied in the coordinates $z$, we can write
\begin{eqnarray}\label{sima.11.1}
\partial_i A(z'(z))\omega^{ip}(z)\partial_p\left[\partial_j B(z'(z))\omega^{jk}(z)\partial_k C(z'(z))\right]+cycle(A, B, C)=0.
\end{eqnarray}
Computing the derivatives, we present this identity as follows
\begin{eqnarray}\label{sima.11.2}
\left.\partial_{i'} A\right|_{z'(z)}\frac{\partial z^{i'}}{\partial z^i}\omega^{ip}(z)\partial_p\left[\left.\left[\partial_{j'} B\left.\frac{\partial z^{j'}}{\partial z^j}\right|_{z(z')}\omega^{jk}(z(z')\left.\frac{\partial z^{k'}}{\partial z^k}\right|_{z(z')}\partial_{k'} C\right]\right|_{z'(z)}\right]+cycle(A, B, C)=0. 
\end{eqnarray}
Using the identity $\partial_p\left[\left.D(z')\right|_{z'(z)}\right]=\frac{\partial z^{p'}}{\partial z^p}\left.\partial_{p'} D(z')\right|_{z'(z)}$, we obtain
\begin{eqnarray}\label{sima.11.3}
\left.\left[\partial_{i'} A(z')\omega^{i'p'}(z')\partial_{p'}\left[\partial_{j'} B(z')\omega^{j'k'}(z')\partial_{k'} C(z')\right]+cycle(A, B, C)\right]\right|_{z'(z)}=0, 
\end{eqnarray}
which is just the Jacobi identity for the bracket $\{A(z'), B(z')\}=\partial_{i'} A\omega^{i' j'}(z')\partial_{j'} B$.  $\blacksquare$

\subsection{Darboux theorem}\label{App B}

\noindent {\bf Lemma B1.} (On rectification of a vector field). Let $V^i(z^j)$ be vector field, nonvanishing at the point $z_0 \in \mathbb{M}_n$. Then there are coordinates\footnote{In this section we use the notation $V^i(z^j)$ and $V^i(y^j)$ instead of $V^i$ and $V'{}^i$ to denote components of the vector $\vec V$ in different coordinate systems.} $y^i$ such that $V^i(y^j)=(1, 0, \ldots , 0)$ at all points $y^j$ in some vicinity of $z_0$. The coordinate $y^1$ has simple geometric meaning: its integral lines are just the integral lines of $\vec V$: $V^i(y^j)=\frac{dy^i}{d\tau}$, where $y^i(\tau)=(y^1=\tau, y^a=C^a)$, and $C^2, \ldots , C^n$ are fixed numbers.

\noindent {\bf Proof.} Without loss of generality we take $z_0=(0, 0, \ldots , 0)$ and $V^1(z_0)\ne 0$, then $V^1(z)\ne 0$ in some vicinity of $z_0$. Write the equations for integral lines 
\begin{eqnarray}\label{apB.1}
\frac{d z^i}{d\tau}=V^i(z^j(\tau)), 
\end{eqnarray}
and solve them with the following initial conditions on the hyperplane $z^1=0$:
\begin{eqnarray}\label{apB.2}
z^1(0)=0, \quad z^2(0)=y^2, \quad \ldots , \quad z^n(0)=y^n,  \quad \mbox{where} \quad y^2, \ldots , y^n \quad \mbox{are fixed numbers}.  
\end{eqnarray}
Denote by 
\begin{eqnarray}\label{apB.3}
z^i(\tau)=f^i(\tau, y^2, \ldots , y^n),
\end{eqnarray}
the integral line that at $\tau=0$ passes through the point $(0, y^2, \ldots , y^n)$.  This determines the nondegenerate mapping
\begin{eqnarray}\label{apB.4}
f:  \quad (\tau, y^2, \ldots , y^n) \quad \rightarrow \quad z^i=f^i(\tau, y^2, \ldots , y^n). 
\end{eqnarray}
The nondegeneracy follows from (\ref{apB.2}) and (\ref{apB.4}) as follows: 
\begin{eqnarray}\label{apB.5}
\left.\det\frac{\partial(z^1, \ldots , z^n)}{\partial(\tau, y^2, \ldots , y^n)}\right|_{z_0}=\det
\left(
\begin{array}{cccc}
V^1(z_0) & 0 & \ldots , & 0 \\
V^a(z_0) & {} & \delta^a{}_b & {}
\end{array}
\right)=V^1(z_0) \ne 0,   
\end{eqnarray} 
so we can take the set 
\begin{eqnarray}\label{apB.6}
y^1=\tau, ~ y^2, ~ \ldots , ~ y^n,
\end{eqnarray}
as new coordinates of $\mathbb{M}_n$, and then the transition functions are given by Eq. (\ref {apB.4}). For the latter use, we note that
\begin{eqnarray}\label{apB.6.1}
\left.\frac{\partial z^2}{\partial y^1}\right|_{z_0}=V^2(z_0), \quad 
\left.\frac{\partial z^2}{\partial y^2}\right|_{z_0}=1, \quad \left.\frac{\partial z^2}{\partial y^\alpha}\right|_{z_0}=0, \quad \alpha=3, 4, \ldots , n. 
\end{eqnarray}
According to (\ref {apB.3}), integral line of the field $\vec V$ in the new system is $y^i(\tau)=(y^1=\tau, ~  y^a=C^a)$,  that is it coincides with the coordinate line of $y^1$, then $V^i(y^j)=\frac{dy^i}{d\tau}=(1, 0, \ldots , 0)$. $\blacksquare$

\noindent {\bf Lemma B2.} Let $\mathbb{M}_n=\{z^k, ~ \omega^{ij}(z^k)\}$ be Poisson manifold with $rank ~  \omega(z_0)\ne 0$. Then there is a pair of scalar functions, say $q\in \mathbb{F}_{\mathbb{M}}$ and $p\in \mathbb{F}_{\mathbb{M}}$, with the 
property $\{q, p\}=1$. Their Hamiltonian fields $\vec V_q$ and $\vec U_p$ are linearly independent and have vanishing Lie bracket, $[\vec V_q, \vec U_p]=0$. \par

\noindent {\bf Proof.}  Without loss of generality we take $\omega^{12}\ne 0$. As the function $q(z^i)$ we take the scalar function of the coordinate $z^2$, its representative in the system $z^i$ is 
\begin{eqnarray}\label{apB.7}
q(z^i)=z^2, \quad \mbox{then its Hamiltonian field is} \quad V^i_q(z^k)=\omega^{ik}\frac{\partial}{\partial z^k} z^2=\omega^{i2}=(\omega^{12}, 0, \omega^{32}, \ldots , \omega^{n2}). 
\end{eqnarray}
In particular, $V^1_q=\omega^{12}\ne 0$.  We rectificate this field according to Lemma B1, then its components in the system $y^j$ are\footnote{Compare this discussion with that of around Eq. (\ref {sim.12.00}).}
\begin{eqnarray}\label{apB.8}
V^i_q(y^j)=(1, 0, \ldots , 0).
\end{eqnarray} 
The representative of the function $q$ in the system $y^j$ is $q(y^j)=z^2(y^j)$, so its bracket with any other function reads
\begin{eqnarray}\label{apB.9}
\{q(y^j), B(y^j)\}=V^i_q(y^j)\frac{\partial}{\partial y^i}B=\frac{\partial}{\partial y^1}B. 
\end{eqnarray} 
Taking as the function $p$ the scalar function of the coordinate $y^1$: $p(y^j)=y^1$, we obtain the desired pair of functions
\begin{eqnarray}\label{apB.10}
\{z^2(y^i), y^1\}=1, \quad \mbox{or, in initial coordinates,} \quad \{z^2, y^1(z^j)\}=1. 
\end{eqnarray} 
In the coordinate system $y^j$, the Hamiltonian fields of these functions are
\begin{eqnarray}\label{apB.11}
V_q^i(y^j)=(1, 0, \ldots , 0), \qquad U_p^i(y^j)=\omega'^{ik}\frac{\partial}{\partial y^k}y^1=\omega'^{i1}=(0, \omega'^{21}, \omega'^{31}, \ldots , \omega'^{n1}). 
\end{eqnarray}
From their manifest form they are linearly independent. Besides, as the Hamiltonian field of a constant vanishes, we have $[\vec V_q, \vec U_p]=-\vec W_{\{q, p\}}=-\vec W_1=0$. $\blacksquare$  \par

\noindent {\bf Lemma B3.}  (On existence  of a pair of canonical coordinates). Let $\mathbb{M}_n=\{z^k, ~ \omega^{ij}(z^k)\}$ be Poisson manifold with $rank ~  \omega(z_0)\ne 0$. Then there are coordinates $q, p, \xi^3, \ldots , \xi^n$ with the properties 
\begin{eqnarray}\label{apB.12}
\{ q, p \}=\omega'^{12}=1, \qquad \{q, \xi^\alpha\}=\omega'^{1\alpha}=0, \qquad \{p, \xi^\alpha\}=\omega'^{2\alpha}=0, 
\end{eqnarray}
\begin{eqnarray}\label{apB.13}
\{\xi^\alpha, \xi^\beta \}=\omega'^{\alpha\beta}(\xi^\gamma), \quad \mbox{that is} \quad \omega'^{\alpha\beta} \quad \mbox{do not depend on} \quad q, p. 
\end{eqnarray} 
In addition, Jacobi identity for $\omega^{ij}$ and Eqs. (\ref{apB.12}) and (\ref{apB.13}) imply the Jacobi identity for $\omega'^{\alpha\beta}$: $\omega'^{\alpha\rho}\partial_\rho \omega'^{\beta\gamma}+cycle=0$. 

\noindent {\bf Proof.} {\bf (A)} We take $q(z^i)=z^2$, and rectify the vector field $V_q^i$ using the Lemma B1. In the process, we obtain the coordinates $y^i$,  the components of the field $V_q^i(y^j)=(1, 0, \ldots , 0)$ in these coordinates, and the scalar function $p(y^j)=y^1$ which obeys 
\begin{eqnarray}\label{apB.14.0}
\{q, p\}=1.
\end{eqnarray}
{\bf (B)} Let $U_p^i(y^j)$ be components of Hamiltonian vector field of the function $p$ in the coordinates $y^j$. According the Lemma B2, $\vec V_q$ and $\vec U_p$ are commuting fields, then 
\begin{eqnarray}\label{apB.14}
0=[\vec V_q, \vec U_p]^i=V_q^k\frac{\partial}{\partial y^k}U_p^i-U_p^k\frac{\partial}{\partial y^k}V_q^i=\frac{\partial U_p^i}{\partial y^1}, \quad  \mbox{implies} \quad U^i_p=U^i_p(y^2, \ldots , y^n),
\end{eqnarray}
that is $\vec U_p$ does not depend on $q$ and $p$. Consider integral lines of the field $\vec U_p$. Taking into account that $U_p^1(y^j)=0$, we have
\begin{eqnarray}\label{apB.15}
\frac{d y^1}{d\lambda}=0, \quad  \mbox{then } \quad y^1=C=const,
\end{eqnarray}
\begin{eqnarray}\label{apB.16}
\frac{d y^a}{d\lambda}=U_p^a(y^2(\lambda), \ldots , y^n(\lambda)). 
\end{eqnarray}
For definiteness, we assume $U_p^2(z_0)\ne 0$. We apply the Lemma B2 to the field $ U_p^a(y^b)$, with $a, b=2, 3, \ldots , n$, that is we solve Eqs. (\ref{apB.16}) with initial conditions on the surface $y^2=0$:
\begin{eqnarray}\label{apB.17}
y^2(0)=0, \quad y^3(0)=\xi^3, \quad \ldots , \quad y^n(0)=\xi^n,  
\end{eqnarray}
Denote solution of the problem as 
\begin{eqnarray}\label{apB.18}
y^a(\lambda)=g^a(\lambda, \xi^3, \ldots , \xi^n), \qquad a=2, 3, \ldots , n.
\end{eqnarray}
These equations are invertible, since (\ref{apB.16})-(\ref{apB.18}) imply (here $\alpha, \beta=3, 4, \ldots , n$)
\begin{eqnarray}\label{apB.19}
\left.\det\frac{\partial(y^2, \ldots , y^n)}{\partial(\lambda, \xi^3, \ldots , \xi^n)}\right|_{z_0}=U_p^2(z_0)\det\frac{\partial y^\alpha(\lambda=0)}{\partial\xi^\beta}=U_p^2(z_0)\det\boldsymbol{1}=U_p^2(z_0)\ne 0. 
\end{eqnarray}
We denote the inverse formulas as follows: 
\begin{eqnarray}\label{apB.20}
\lambda=\tilde g(y^2, \ldots , y^n), \quad \xi^3=\tilde g^3(y^2, \ldots , y^n), \quad \ldots , \quad \xi^n=\tilde g^n(y^2, \ldots , y^n),
\end{eqnarray}
and introduce the new coordinates
\begin{eqnarray}\label{apB.21}
(y^1, y^2, y^3, \ldots , y^n) \quad \rightarrow \quad ( y^1, \lambda(y^a), \xi^\alpha(y^a)), \qquad a=2, 3, \ldots n, \quad \alpha=3, 4, \ldots , n,
\end{eqnarray}
with the transition functions (\ref{apB.20}).  Integral lines of the fields $U$ and $V$ in the new coordinates are ($C, \lambda, \xi^3, \ldots \xi^n)$ and $(y^1=\tau, \tilde g(y^2, \ldots , y^n), \tilde g^\alpha(y^2, \ldots , y^n)$. Along the integral lines of $U$ only the second coordinate $\lambda$ changes.  Along the integral lines of $V$ changes the first coordinate, $y^1=\tau$, while $\lambda$ and $\xi^\alpha$, being  functions of $y^2, \ldots , y^n$, remain constants. Therefore, in these coordinates both fields are straightened: $V^i_q=(1, 0, 0, \ldots , 0)$, $U^i_p=(0, 1, 0, \ldots , 0)$. \par 

{\bf (C)} The Poisson brackets of $q$ and $p$ with scalar functions of the coordinates $\xi^\alpha$, $\alpha=3, 4, \ldots , n$ vanish
\begin{eqnarray}\label{apB.22}
\{q, \xi^\alpha\}=V_q(\xi^\alpha)=\frac{\partial\xi^\alpha}{\partial\tau}=0, \qquad 
\{p, \xi^\alpha\}=V_p(\xi^\alpha)=\frac{\partial\xi^\alpha}{\partial\lambda}=0.
\end{eqnarray}
So, the functions $q$ $p$, and $\xi^\alpha$ obey the equation (\ref{apB.12}). \par 

{\bf (D)} The last step is to introduce the mapping 
\begin{eqnarray}\label{apB.23}
(y^1, y^2, y^3, \ldots , y^n) \quad \rightarrow \quad ( q=z^2(y^j), ~ p=y^1, ~ \xi^\alpha=\tilde g(y^2, \ldots , y^n)). 
\end{eqnarray}
Its invertibility follows from the direct computation 
\begin{eqnarray}\label{apB.24}
\left.\det\frac{\partial(q, p, \xi^3, \ldots , \xi^n)}{\partial(y^1, y^2, y^3, \ldots , y^n)}\right|_{z_0}=\left.\det
\left(
\begin{array}{ccccc}
\frac{\partial z^2}{\partial y^1} & \frac{\partial z^2}{\partial y^2} & \ldots & {} &\frac{\partial z^2}{\partial y^n}  \\
\frac{\partial y^1}{\partial y^1} & \frac{\partial y^1}{\partial y^2} & \ldots & {} & \frac{\partial y^1}{\partial y^n}  \\
\ldots & \ldots &  {} & {} & {} \\
\frac{\partial\xi^\alpha}{\partial y^1}  & \frac{\partial\xi^\alpha}{\partial y^2}  & {} & \frac{\partial\xi^\alpha}{\partial y^\beta} & {} \\
\ldots & \ldots &  {} & {} & {} \\
\end{array}
\right)\right|_{z_0}=\det
\left(
\begin{array}{ccccc}
0 & 1 & 0 &\ldots  & 0 \\
1 & 0 & 0 &\ldots  & 0 \\
\ldots & \ldots &  {} & {} & {} \\
0  & \frac{\partial\xi^\alpha}{\partial y^2}  & {} & {\boldsymbol 1} & {} \\
\ldots & \ldots &  {} & {} & {} \\
\end{array}
\right)=-1. 
\end{eqnarray} 
In the computation we used the equations (\ref{apB.6.1}), (\ref{apB.7}), (\ref{apB.20}) and (\ref{apB.19}). In particular: $\left.\frac{\partial z^2}{\partial y^1}\right|_{z_0}=\left.\frac{\partial f^2(\tau, y^2, \ldots , y^n)}{\partial\tau} \right|_{z_0}=V^2_q|_{z_0}=\omega^{22}=0$.  
Therefore we can take $q, p, \xi^\alpha$ as a coordinate system on $\mathbb{M}_n$. As we saw above, the coordinates obey the desired property (\ref{apB.12}). To confirm (\ref{apB.13}), we use $\{q, \xi^\alpha\}=0$ in the Jacobi identity, obtaining 
\begin{eqnarray}\label{apB.25}
\{q, \{\xi^\alpha, \xi^\beta \}\}=-\{\xi^\alpha, \{\xi^\beta, q \}\}-\{\xi^\beta , \{q, \xi^\alpha \}\}=0,  \quad \mbox{or} \quad \frac{\partial}{\partial p} \{\xi^\alpha, \xi^\beta \}=0. 
\end{eqnarray}
So $\{\xi^\alpha, \xi^\beta \}\equiv\omega'^{\alpha\beta}$ does not depend on $p$. Similar computation of  $\{p, \{\xi^\alpha, \xi^\beta \}\}$ implies, that $\omega'^{\alpha\beta}$ does not depend on $q$. $\blacksquare$ \par 

If $rank ~ \omega'^{\alpha\beta}(\xi^\gamma)\ne 0$,  the manifold $\mathbb{M}_{n-2}=\{\xi^\gamma, ~ \omega'^{\alpha\beta}(\xi^\gamma)\}$, in turn, satisfies the conditions of the Lemma B3.

\noindent {\bf Generalized Darboux theorem.} Let $\mathbb{M}_n=\{z^k, ~ \omega^{ij}(z^k)\}$ be Poisson manifold with 
$rank~\omega=2k$ at the point $z_0^i\in\mathbb{M}_n$. Then there are local coordinates, where $\omega$ has the form: 
\begin{eqnarray}\label{apB.26}
\omega'=\left(
\begin{array}{ccc}
0_{p\times p}  & 0 & 0 \\
0  & 0_{k\times k} & 1_{k\times k} \\
0  & -1_{k\times k} & 0_{k\times k} 
\end{array}\right), \qquad p=n-2k, 
\end{eqnarray} 
at all points in some vicinity of $z^i_0$.

\noindent {\bf Proof.} The proof is carried out by induction on the pairs of canonical coordinates constructed in Lemma B3. After $k$ steps, we get the coordinates $\xi^\alpha, q^b, p^c$, $\alpha=1, 2, \ldots , n-2k$, $b, c=1, 2, \ldots , k$, in which the tensor $\omega$ has the block-diagonal form
\begin{eqnarray}\label{apB.27}
\omega'=\left(
\begin{array}{ccc}
\omega'^{\alpha\beta}  & 0 & 0 \\
0  & 0_{k\times k} & 1_{k\times k} \\
0  & -1_{k\times k} & 0_{k\times k} 
\end{array}\right), 
\end{eqnarray}
and $\omega'^{\alpha\beta}=\{\xi^\alpha, \xi^\beta\}$. From the rank condition and from the manifest form (\ref{apB.27}) of the matrix $\omega'$, we have $2k=rank ~ \omega'=rank ~ \omega'^{\alpha\beta}+2k$, or $rank ~ \omega'^{\alpha\beta}=0$. This 
implies $\omega'^{\alpha\beta}=0$ for all $\alpha$ and $\beta$. $\blacksquare$

{\bf Affirmation B1.} Let $Q(z^i)$ be first integral of the Hamiltonian system $\dot z^i=\omega^{ij}\partial_j H$ with a non degenerate tensor $\omega^{ij}$. Then solution of this system of $n$ equations reduces to the solution of a Hamiltonian system composed by $n-2$ equations. 

\noindent {\bf Proof.}  Introduce the coordinates $z'^i$: $z'^1=z^1, z'^2=Q(z^i), z'^3=z^3, \ldots , z'^n=z^n$, thus turning $Q$ into the second coordinate of the new system. Applying the Lemmas B2 and B3, we construct the coordinates $q, p, \xi^\alpha$ with $q=Q$. Poisson tensor in these coordinates has the form 
\begin{eqnarray}\label{apB.28}
\omega'=\left(
\begin{array}{ccc}
 0 & 1 & 0 \\
-1  & 0 & 0 \\
0  & 0 & \omega'^{\alpha\beta}(\xi^\gamma) 
\end{array}\right).  
\end{eqnarray}
Consider our Hamiltonian equations in these coordinates. The equation $\dot q=\partial_p H'$ together with $q=c_2=const$ implies that $H'$ does not depend on $p$: $H'=H'(q, \xi^\gamma)$. Then on the surface $q=c_1=const$, the original system is equivalent to 
\begin{eqnarray}\label{apB.29}
\dot p=-\left.\partial_q H'(q, \xi^\gamma)\right|_{q=c_2}, 
\end{eqnarray}
\begin{eqnarray}\label{apB.30}
\dot \xi^\alpha=\omega'^{\alpha\beta}\partial_\beta H'(c_2, \xi^\gamma),  \qquad \alpha=3, 4, \ldots , n. 
\end{eqnarray}
The $n-2$ Hamiltonian equations (\ref{apB.30}) can be solved separately from (\ref{apB.29}), let $\xi^\alpha(\tau, c_2, \ldots , c_n)$ be their general solution.  Using these functions in Eq. (\ref{apB.29}), the latter is solved by direct integration: $p=-\int d\tau\left.\partial_q H'(q, \xi^\gamma)\right|_{q=c_2, \xi=\xi(\tau, c_2, \ldots , c_n)}$. $\blacksquare$

It should be noted that the range of applicability of this affirmation in applications is rather restricted. Indeed, to find manifest form of the equations (\ref{apB.30}),  we need to rectify two vector fields. And for this, it is necessary to solve twice a system of equations like the original system!

\subsection{Frobenius theorem}\label{App C}

The equation $\partial_x X(x, y, z)=0$ has two functionally independent solutions: $X_1=y$ and $X_2=z$. Frobenius theorem can be thought as a generalization of this result to the case of the system of first-order partial differential equations 
$A_a^i(z^k)\partial_i X(z^k)=0$.  The theorem can also be reformulated in a purely geometric language, see the end of this section. 

We will need some properties of vector fields and their integral lines on a smooth  manifold $\mathbb{M}_n=\{z^i, ~ i=1, 2, \ldots , n\}$. We recall that integral line of the vector field $V^i(z^k)$ on $\mathbb{M}_n$ is a solution $z^i(\tau)$ to  $\frac{d z^i(\tau)}{d\tau}=V^i(z^k(\tau))$.  As before, we assume that through each point of the manifold passes unique integral line of $\vec V$. By $\{ \mathbb{N}_k^{\vec c}, ~ \vec c\in \mathbb{R}^{n-k}\}$ we denote a foliation of $\mathbb{M}_n$ (see Sect. \ref{SymNot_C}), with the leaves
\begin{eqnarray}\label{apC.0.3}
\mathbb{N}_k^{\vec c}=\{z^i\in \mathbb{M}_n, ~  F^\alpha(z^i)=c^\alpha \}. 
\end{eqnarray}
 
\noindent {\bf Affirmation C1.} Let $\vec V(z^k)$ be vector field on $\mathbb{M}_n$ and $F(z^k)$ be  scalar function with non vanishing gradient. The following two conditions are equivalent: \par

\noindent {\bf (A)} $\vec V$ touches  the surfaces $F(z^k)=c=const$: $V^i\partial_i F=0$ at each point $z^k\in \mathbb{M}_n$.  

\noindent {\bf (B)} $\vec V$ is tangent\footnote{See the definition of a vector field tangent to a submanifold on page \pageref{tang}.}  to the surfaces $F(z^k)=c=const$, that is integral lines of $\vec V$ lie on the surfaces. 

\noindent {\bf Proof.} Let $z^k(\tau)$ be an integral line of $\vec V$. Then the Affirmation follows immediately from the equality
\begin{eqnarray}\label{apC.0.5}
\frac{d}{d\tau}F(z^k(\tau))=
\left.V^i(z^k)\partial_i F(z^k)\right|_{z(\tau)}.  \qquad \blacksquare
\end{eqnarray} 
Evidently, the same is true for a set of vector fields: \par 

\noindent {\bf Affirmation C2.}  Let 
$\vec A_1(z^k), \ldots , \vec A_k(z^k)$ be vector fields on $\mathbb{M}_n$, linearly independent at each point $z\in {\mathbb M}_n$, and  $\{ \mathbb{N}_k^{\vec c}, ~ \vec c\in \mathbb{R}^{n-k} \}$ be a foliation of $\mathbb{M}_n$.  The following two conditions are equivalent: \par 

\noindent {\bf (A)}  The vectors $\vec A_a$ touch  ${\mathbb N}_k^{\vec c}$  at each point $z^k\in \mathbb{M}_n$: $A_a^i\partial_i F^\alpha=0$ at each point $z\in {\mathbb M}_n$.

\noindent {\bf (B)}  The vectors $\vec A_a$ are tangent to ${\mathbb N}_k^{\vec c}$, that is each integral line of each $\vec A_a$ lies in one of the submanifolds $\mathbb{N}_k^{\vec c}$ (hence  $\vec A_a(z^k)$ form a basis of $\mathbb{T}_{\mathbb{N}_k^c}(z^k)$.

\noindent {\bf Lemma C1.} There is a set of $k$ linearly independent vector fields $\vec U_a(z^k)$ on $\mathbb{M}_n$ with the following properties. \par

\noindent {\bf (A)} For any $z^k\in\mathbb{M}_n$, the vectors  $\vec U_a(z^k)$ touch the submanifold $\mathbb{N}_k^{\vec c}$ that passes through this point: 
\begin{eqnarray}\label{apC.0.4}
U_a^i\partial_i F^\alpha=0. 
\end{eqnarray} 
At each point they form a basis of tangent space to the submanifold. \par

\noindent {\bf (B)} Integral lines of $\vec U_a$ that pass through $z^k\in\mathbb{M}_n$, lye in $\mathbb{N}_k^{\vec c}$ that passes through this point. \par

\noindent {\bf (C)} $\vec U_a$ are commuting fields
\begin{eqnarray}\label{apC.0.5}
[\vec U_a, \vec U_b]=0.
\end{eqnarray} 

\noindent {\bf Proof.} Introduce the coordinates, adapted with the foliation: $z^k\rightarrow y^k=(y^\alpha, y^a)$, with the transition functions $y^a=z^a$, $y^\alpha=F^\alpha(z^\beta, z^b)$.  In these coordinates the sumanifolds $\mathbb{N}_k^{\vec c}$ look like hyperplanes:
\begin{eqnarray}\label{apC.0.6}
\mathbb{N}_k^{\vec c}=\{y^i\in\mathbb{M}_n, ~ y^\alpha=c^\alpha \}, 
\end{eqnarray} 
and $y^a$ can be taken as local coordinates of $\mathbb{N}_k^{\vec c}$. Consider the vector fields $\vec U_a$ on $\mathbb{M}_n$, which in the system $y^k$ have the following components: $U_a^i(y^k)=\delta_a{}^i$. Their integral lines are just lines of the coordinates $y^a$ of the submanifolds $\mathbb{N}_k^{\vec c}$. Evidently, the fields obey the conditions (A)-(C) of the Lemma. Their explicit form in the original coordinates is as follows:
\begin{eqnarray}\label{apC.0.7}
U_a^i(z^k)=\left.\left[\frac{\partial z^i}{\partial y^j}U_a^j(y^k)\right]\right|_{y(z)}= (U_a^b, ~ U_a^\beta)=\left(\delta_a{}^b, ~ \left.\frac{\partial f^\beta(z^c, y^\gamma)}{\partial z^a}\right|_{y^\gamma\rightarrow F^\gamma(z^b, z^\beta)}\right),
\end{eqnarray} 
where $f^\beta(z^c, y^\gamma)$ is solution to the system $F^\beta(f^\beta, z^c)=y^\gamma$. 
Since (\ref{apC.0.4}) and (\ref{apC.0.5}) are covariant equations, the fields (\ref{apC.0.7}) satisfy them in the original coordinates $z^k$. $\blacksquare$

\noindent {\bf Lemma C2.} An invertible linear combination of vector fields with closed algebra also form a closed algebra:
\begin{eqnarray}\label{apC.0.1}
\mbox{if} \quad \vec V_a=b_a{}^b\vec U_b, \quad \det b\ne 0, \quad \mbox{and} \quad  [\vec U_a, \vec U_b]=c_{ab}{}^c\vec U_c, \quad \mbox{then} \quad  [\vec V_a, \vec V_b]=\gamma_{ab}{}^c\vec V_c. 
\end{eqnarray}

\noindent {\bf Proof.} This  follows from direct alculation, that also implies
\begin{eqnarray}\label{apC.0.2}
\gamma_{ab}{}^c=b_a{}^d b_b{}^e c_{de}{}^f \tilde b_f{}^c+V_a^j(\partial_j b_b{}^f)\tilde b_f{}^c-(a\leftrightarrow b),  
\end{eqnarray} 
where $\tilde b$ is inverse for $b$. $\blacksquare$

\noindent {\bf Lemma C3.} Let $\vec A_1, \ldots , \vec A_k$ is a set of linealy independent vector fields on $\mathbb{M}_n$,  with closed algebra of commutators
\begin{eqnarray}\label{apC.1}
[ \vec A_a, \vec A_b]^i=c_{ab}{}^d(z^k)A_d^i. 
\end{eqnarray}
Then there is a set of linearly independent fields $\vec V_a$, which are linear combinations of $\vec A_a$ and have vanishing commutators
\begin{eqnarray}\label{apC.2}
[ \vec V_a, \vec V_b]=0. 
\end{eqnarray} \par

\noindent {\bf Proof.} The components $A_a^i=(a_a{}^b, b_a{}^\beta)$ of  linearly independent fields form $k\times n$ matrix with rank equal $k$. Without loss of generality we assume $\det a_a{}^b\ne 0$, and let $\tilde a_a{}^b$ be the inverse matrix. We show that 
$\vec V_a\equiv \tilde a_a{}^b\vec A_b$ are the desired fields.

The expressions (\ref{apC.1}) with components $i=c$ can be solved with respect to $c_{ab}{}^d$ as follows: $[ \vec A_a, \vec A_b]^c=c_{ab}{}^d a_d{}^c$ implies $c_{ab}{}^c=[ \vec A_a, \vec A_b]^d\tilde a_d{}^c$. Using this equality, we exclude $c_{ab}{}^c$ from the expressions  (\ref{apC.1}) with $i=\beta$, obtaining $[ \vec A_a, \vec A_b]^\beta=[ \vec A_a, \vec A_b]^d\tilde a_d{}^c b_c{}^\beta$. In more detail, this reads
\begin{eqnarray}\label{apC.3}
A_a{}^i\partial_i b_b^\beta-(a\leftrightarrow b)=A_a{}^i(\partial_i a_b{}^d)\tilde a_d{}^c b_c{}^\beta-(a\leftrightarrow b)= \qquad \qquad \qquad \qquad \cr 
A_a{}^i\partial_i(a_b{}^d\tilde a_d{}^c b_c{}^\beta)-A_a{}^i a_b{}^d\partial_i(\tilde a_d{}^c b_c{}^\beta)-(a\leftrightarrow b)= 
A_a{}^i\partial_i b_b^\beta-A_a{}^i a_b{}^d\partial_i(\tilde a_d{}^c b_c{}^\beta)-(a\leftrightarrow b), 
\end{eqnarray}
which implies $A_a{}^i a_b{}^d\partial_i(\tilde a_d{}^c b_c{}^\beta)-(a\leftrightarrow b)=0$. Contraction of this equality with $\tilde a_e{}^a\tilde a_f{}^b$ gives the following relation between components of the fields with closed commutator algebra: 
\begin{eqnarray}\label{apC.4}
\tilde a_a{}^c A_c{}^i\partial_i(\tilde a_b{}^d b_d{}^\beta)-(a\leftrightarrow b)=0. 
\end{eqnarray}
Now, the fields $\vec V_a\equiv \tilde a_a{}^c\vec A_c$ with the components $V_a{}^i=(V_a{}^b, ~  V_a{}^\beta)=(\delta_a{}^b, ~ \tilde a_a{}^c b_c{}^\beta)$ satisfy the conditions of the Lemma. Indeed, $[\vec V_a, \vec V_b]^c=V_a{}^i\partial_i\delta_b^c-(a\leftrightarrow b)=0$, and 
$[\vec V_a, \vec V_b]^\beta=V_a{}^i\partial_i(\tilde a_b{}^d b_d{}^\beta)-(a\leftrightarrow b)=0$ due to (\ref{apC.4}). $\blacksquare$

Given vector field $V^i(z^k)$, let us denote $\varphi^i(\tau, z_0)$ the unique solution to the problem 
\begin{eqnarray}\label{apC.5}
\frac{d z^i}{d\tau}=V^i(z^k(\tau)),  \qquad z^i(0)=z^i_0. 
\end{eqnarray}
For any fixed value of $\tau$, the integral lines $\varphi^i(\tau, z), ~ z\in\mathbb{M}_n $ determine the transformation
\begin{eqnarray}\label{apC.6}
\varphi_\tau: \mathbb{M}_n \rightarrow \mathbb{M}_n,  \qquad z^i \rightarrow \varphi^i(\tau, z^k).
\end{eqnarray}
Sometimes we will also use the coordinate-free notation $\varphi_\tau(z)$ for the integral line $\varphi^i(\tau, z^k)$.  Composition of two transformations has the property
\begin{eqnarray}\label{apC.7}
\varphi_\tau\circ \varphi_s=\varphi_{\tau+s}.
\end{eqnarray}
Indeed, $\varphi^i(\tau, \varphi^j(s, z^k))$ and $\varphi^i(\tau+s, z^k)$ as functions of $\tau$ obey the problem (\ref{apC.5}) with $z^i_0=\varphi^i(s, z^k)$. Since the problem has unique solution, they coincide. So, the set of transformations $\{ \varphi_\tau, ~  \tau\in \mathbb{R}\}$ is a one-parametric Lie group with the group product being the composition law (\ref{apC.7}). 

Let $\varphi_\tau$ and $\psi_\lambda$ be the one-parametric groups criated by linearly independent fields $V^i(z^k)$ and $U^i(z^k)$. There is a remarkable relation between commutativity of the transformations and of the vector fields.

\noindent {\bf Lemma C4.} The following two conditions are equivalent: {\bf (A)} $\varphi_\tau\circ \psi_\lambda(z^k)=\psi_\lambda\circ \varphi_\tau(z^k)$ for all $\tau$, $\lambda$ and $z^k$. {\bf (B)} $[\vec V(z), \vec U(z) ]=0$ for all $z$. 

\noindent {\bf Proof.} (A) $\rightarrow$  (B). Expanding in series of Taylor, we obtain $[\varphi_\tau\circ \psi_\lambda(z^k)-\psi_\lambda\circ \varphi_\tau(z^k)]^i=\varphi^i(\tau, \psi^j(\lambda, z^k))-\psi^i(\lambda, \varphi^j(\tau, z^k))= [\vec V(z), \vec U(z)]^i\tau\lambda+O^2(\tau)+O^2(\lambda)+O^3(\tau, \lambda)$. Since l.h.s. vanishes for any $\tau$ and $\lambda$, we conclude $[\vec V(z), \vec U(z) ]=0$. 

(B) $\rightarrow$  (A). Consider the fields $\vec V$ and $\vec U$ in the coordinates $y^k$ of the Lemma B1. Then $\vec V(y^k)=(1, 0, \ldots , 0)$ and its integral line through the point $y^k$ is 
\begin{eqnarray}\label{apC.8.1}
\varphi^i(\tau, y^k)=(y^1+\tau, y^2,  \ldots , y^n). 
\end{eqnarray}
Besides, the condition (B) reads $0=[ \vec V, \vec U]^i=\frac{\partial U^i}{\partial y^1}$, 
that is the field $\vec U$ does not depend on $y^1$. Consider $\psi_\lambda\circ \varphi_\tau(z^k)$ and $\varphi_\tau\circ \psi_\lambda(z^k)$ in the system $y^k$ as functions of $\lambda$. Using (\ref {apC.8.1}), we can write
\begin{eqnarray}\label{apC.9}
\psi^i(\lambda, \varphi^j(\tau, y^k))=(\psi^1, ~ \psi^2, \ldots , ~ \psi^n), \quad \mbox{then} \quad 
\psi^i(0, \varphi^j(\tau, y^k))=\varphi^i(\tau, y^k)=(y^1+\tau, ~ y^2, ~  \ldots , ~ y^n), 
\end{eqnarray}
\begin{eqnarray}\label{apC.10}
\varphi^i(\tau, \psi^j(\lambda, y^k))=(\psi^1+\tau, ~ \psi^2, \ldots , ~ \psi^n), \quad \mbox{then} \quad 
\varphi^i(\tau, \psi^j(0, y^k))=(y^1+\tau, ~ y^2, ~  \ldots , ~ y^n).  
\end{eqnarray}
By construction, $\psi^i(\lambda)$ satisfy the equation
\begin{eqnarray}\label{apC.11}
\frac{d x^i}{d\lambda}=U^i(x^2, x^3, \ldots , x^n).
\end{eqnarray}
As the r.h.s. of this equation does not depend on $x^1$, the function $\varphi^i(\lambda)$ also satisfy this equation. Besides, $\psi^i(\lambda)$ and $\varphi^i(\lambda)$ satisfy the same initial conditions, see (\ref{apC.9}) and (\ref{apC.10}).  Hence they coincide. $\blacksquare$

Any set of coordinate lines, say the lines of the coordinates $z^1, z^2, \ldots , z^k$, can be used to construct a set of commuting vector fields. They are the tangent fields to the coordinate lines. The following Lemma is an inversion of this statement. It also generalizes the Lemma B1 to the case of several fields.

\noindent {\bf Lemma C5.} (On rectification of the commuting vector fields). Let $\vec V_1, \vec V_2, \ldots , \vec V_k$ be linearly independent and commuting vector fields in vicinity of  $z_0\in\mathbb{M}_n$: $[ \vec V_a, \vec V_b]=0$. Then: \par 

\noindent  There are coordinates $y^i=(y^a, y^\alpha)$, $\alpha=k+1, \ldots , n$, where the fields $\vec V_a$ are tangent to the coordinate lines $y^a$: $V_a^i(y^j)=\delta_a^i$, $a=1, 2, \ldots, k$.   \par

Notice the immediate consequences of the Lemma: through each point $z_1\in \mathbb{M}_n$ passes a surface $\mathbb{N}_k$ such that $\vec V_1(z), \vec V_2(z), \ldots , \vec V_k(z)$ form a basis of the tangent spaces $\mathbb{T}_{\mathbb{N}}(z)$ at any point $z\in\mathbb{N}_k$. Integral lines of the fields $\vec V_a$, that cross $\mathbb{N}_k$, entirely lie in $\mathbb{N}_k$. Evidently, in the coordinates $y^k$ these surfaces are given by the equations $y^\alpha=c^\alpha=const$.  

\noindent {\bf Proof.}  Without loss of generality we assume that the point $z_0$ has null coordinates. Selecting the appropriate $n-k$ vectors among the basic vectors of coordinate lines, say $\vec e_\alpha$, with coordinates $e_\alpha^i=\delta_\alpha^i$, $\alpha=k+1,  \ldots, n$, we complete the vectors $\vec V_a(z_0)$ up to a basis of $\mathbb{T}_{\mathbb{M}}(z_0)$. Then determinant of the matrix composed from components of the basic vectors is not equal to zero at $z_0$
\begin{eqnarray}\label{apC.12}
\det(\vec V_1, \ldots , \vec V_k, \vec e_{k+1}, \ldots , \vec e_n)|_{z_0=0}\ne 0. 
\end{eqnarray}

Denote $\varphi_{\tau_a}$ the one-parametric group (\ref{apC.6}) criated by the field $\vec V_a$.  Consider the mapping $h: O(\vec 0)\in\mathbb{R}^n\rightarrow\mathbb{M}_n$ defined according the rule 
\begin{eqnarray}\label{apC.13}
z=h(\tau_1, \ldots \tau_k, y^1, \ldots, y^{n-k})=\varphi_{\tau_1}\circ \ldots \circ \varphi_{\tau_k}(0, \ldots , 0, y^{k+1}, \ldots y^n).
\end{eqnarray}
Derivatives of this function at the point $\tau_a=y^\alpha=0$ are $\left.\frac{dh}{d\tau_a}\right|_{0}=\left.\frac{d}{d\tau_a}\varphi_{\tau_a}(0, \ldots , 0, 0,  \ldots , 0)\right|_{\tau_a=0}=\vec V_a(0)$ and $\left.\frac{dh}{dy^\alpha}\right|_{0}=\left.\frac{d}{dz^\alpha}(0, \ldots , 0, 0, \ldots , y^\alpha, \ldots , 0)\right|_{z^\alpha=0}=(0, \ldots , 0, 0, \ldots , 1, \ldots , 0)=\vec e_\alpha$. Then $\left.\det\frac{\partial(z^1, z^2, \ldots , z^n)}{\partial(\tau_1, \ldots \tau_k, y^{k+1}, \ldots, y^n)}\right|_0=\det(\vec V_1, \ldots , \vec V_k, \vec e_{k+1}, \ldots , \vec e_n)|_{0}\ne 0$, see (\ref{apC.12}). So the mapping (\ref{apC.13}) is invertible, and we can take $y^i\equiv(\tau_a, y^\alpha)$ as a coordinate system on $\mathbb{M}_n$. The transition functions are given by Eq. (\ref{apC.13}). 

Consider the integral line $\varphi_{s_a}(z)$ of the field $\vec V_a$ through some point $z$.  According to Lemma C4, commutativity of the fields implies the commutativity of their one-parametric groups, so we have 
\begin{eqnarray}\label{apC.14}
\varphi_{s_a}(z)=\varphi_{s_a}\circ\varphi_{\tau_1}\circ \ldots\varphi_{\tau_a}\ldots\circ\varphi_{\tau_k}(0, y^\alpha)=
\varphi_{\tau_1}\circ \ldots\varphi_{\tau_a+s_a}\ldots\circ\varphi_{\tau_k}(0, y^\alpha)=h(\tau_1, \ldots , \tau_a+s_a , \ldots , \tau_k, y^\alpha). 
\end{eqnarray}
This shows that integral lines of $\vec V_a$ are the coordinate lines of $y^a$\,-coordinate of the new system. Hence the integral lines lie in the submanifolds  $\mathbb{N}_k=\{ y^k\in\mathbb{M}_n, ~  y^\alpha=c^\alpha=const \}$. 

To find equations of these surfaces in the original coordinates, denote $\tilde h$ the inverse mapping of (\ref{apC.13}). Let the point $z_1$ has coordinates 
$\tau_1, \ldots , \tau_k, c^{k+1}, \ldots c^n$ in the system $y^i$. Then the submanifold is 
$\mathbb{N}_k=\{z\in\mathbb{M}_n, ~  \tilde h^\alpha(z^i)=c^\alpha\}$.   $\blacksquare$

\noindent {\bf Frobenius theorem.} Let $A_a^i(z^k)$, $a=1, 2, \ldots k$ be a set of functions with $rank ~ A=k$. The system of first-order partial differential equations 
\begin{eqnarray}\label{apC.15}
A_a^i(z^k)\partial_i X(z^k)=0, 
\end{eqnarray}
has $n-k$ functionally independent solutions, if and only if the vectors $\vec A_a$ form a set with closed algebra
\begin{eqnarray}\label{apC.16}
[\vec A_a(z^k),  \vec A_b(z^k)]=c_{ab}{}^c(z^k)\vec A_a(z^k). 
\end{eqnarray}

\noindent {\bf Proof.} Let the functions $F^\alpha(z^k)$, $\alpha=1, 2 , \ldots , n-k$ represent the solutions:
\begin{eqnarray}\label{apC.17}
A_a^i(z^k)\partial_i F^\alpha(z^k)=0.  
\end{eqnarray}
Consider the foliation $\{ \mathbb{N}_k^{\vec c}, ~ \vec c\in \mathbb{R}^{n-k} \}$ determined by $F^\alpha$ according to Eq. (\ref{apC.0.3}), and let $\vec U_a(z^k)$ be vector fields described in Lemma C1. 

Denoting $z_a^i(\tau)$ integral lines of $\vec A_a(z^k)$, we have $\frac{d }{d\tau}F^\alpha(z_a^i(\tau))=\left.A_a^i(z^k)\partial_i F^\alpha(z^k)\right|_{z_a^i(\tau)}=0$ according to (\ref{apC.17}). Then $F^\alpha(z_a^i(\tau))=c^\alpha=const$, that is integral lines of $\vec A_a(z^k)$ lie in $\mathbb{N}_k^{\vec c}$, and $\vec A_a(z^k)$ are tangent vectors to this submanifold at each point. Then we can present them 
through the basic vectors $\vec U_b$: $\vec A_a=b_a{}^b\vec U_b$ of Lemma C1. According to Lemma C1, $[\vec U_a, \vec U_b]=0$. According to Lemma C2, this implies (\ref{apC.16}). 

Let (\ref{apC.16}) is satisfied. Assuming $A_a^i=(a_a{}^b, b_a{}^\beta)$ with $\det a\ne 0$ (see Lemma C3), we write the system (\ref{apC.15}) in the equivalent form $\tilde a_a{}^b A_b^i(z^k)\partial_i X(z^k)\equiv V_a^i(z^k)\partial_i X(z^k)=0$. According to Lemma C3, we have $[\vec V_a, \vec V_b]=0$. According to Lemma C5, there are coordinates $y^k$ where $V_a^i(y^k)=\delta_a{}^i$. In these coordinates our system acquires the form $\frac{\partial}{\partial y^a}X'(y^\beta, y^b)=0$. The functions $F^\beta(y^\beta, y^b)=y^\beta$ give $n-k$ functionally independent solutions. $\blacksquare$

\noindent {\bf Frobenius theorem, geometric formulation.} Let 
$\vec A_1(z^k), \ldots , \vec A_k(z^k)$ be linearly independent vector fields on $\mathbb{M}_n$.
The following two conditions are equivalent: \par 

\noindent {\bf (A)} The  fields $\vec A_a$  form closed  algebra: 
\begin{eqnarray}\label{apC.18}
[\vec A_a(z^i), \vec A_b(z^i)]=c_{ab}{}^c(z^i) \vec A_c(z^i). 
\end{eqnarray}
\noindent {\bf (B)}  There is a foliation $\{ \mathbb{N}_k^{\vec c}, ~ \vec c\in \mathbb{R}^{n-k} \}$ of $\mathbb{M}_n$ such that the fields $\vec A_a(z^k)$ touch the leaf $\mathbb{N}_k^{\vec c}$ (see Eq. (\ref{apC.0.3}) ) at each point $z^k\in \mathbb{M}_n$ (hence $\vec A_a$ form a basis of $\mathbb{T}_{\mathbb{N}_k^c}(z^k)$, see the Affirmation C2).  

\noindent {\bf Proof.}  $(B)\rightarrow (A)$. Consider $z_0\in\mathbb{M}_n$ and let $z_0\in\mathbb{N}_k^{\vec c}$, where $\mathbb{N}_k^{\vec c}$ is one of submanifolds  specified in (B). Let $z^i(\tau)$ be integral line of the field $[\vec A_a, \vec A_b]^i$, which at $\tau=0$ passes through $z_0$. We get 
\begin{eqnarray}\label{apC.19}
\frac {d}{d\tau}F^\alpha(z^i(\tau))=\left. [\vec A_a, \vec A_b]^i\partial_i F^\alpha\right|_{z^i(\tau)}=\left.\left[\vec A_a(\vec A_b(F^\alpha))-(a\leftrightarrow b)\right)\right|_{z^i(\tau)}=0,
\end{eqnarray}
since $\vec A_b(F^\alpha)=A_b^i\partial_i F^\alpha=0$. The equality  (\ref{apC.19}) implies that integral line of the field $[\vec A_a, \vec A_b]^i$  through $z_0$ entirely lies in $\mathbb{N}_k^{\vec c}$, so the vector $[\vec A_a, \vec A_b]^i(z_0)$ is tangent to ${\mathbb N}_k^{\vec c}(z_0)$.  Hence it can be presented through the basic vectors $\vec A_a$, which gives the desired result (\ref{apC.18}). 

$(A)\rightarrow (B)$. Let (\ref{apC.18}) is satisfied. Using Lemma C3, we construct $k$ linearly independent and commuting fields $\vec V_a$. According to Lemma C5, there are coordinates $y^k$ where $V_a^i(y^k)=\delta_a{}^i$. Consider the foliation $\{ \mathbb{N}_k^{\vec c}, ~ \vec c\in \mathbb{R}^{n-k} \}$ where $\mathbb{N}_k^{\vec c}=\{ z^k\in \mathbb{M}_n, y^\alpha=c^\alpha=const \}$. By construction, $\vec V_a\in \mathbb{T}_{\mathbb{N}_k^c}$  and form a basis of $\mathbb{T}_{\mathbb{N}_k^c}$ at each point $z^k\in \mathbb{M}_n$. According to Lemma C3, the linearly independent vectors $\vec A$ are linear combinations of $\vec V_a$, so they also form a basis of  $\mathbb{T}_{\mathbb{N}_k^c}$ at each point $z^k\in \mathbb{M}_n$.  $\blacksquare$

\end{document}